\documentclass[reqno]{amsart}
\theoremstyle{plain}
\usepackage{latexsym}
\usepackage{amssymb}
\usepackage{mathrsfs}
\usepackage{wasysym}
\usepackage{epsfig}
\def\Xint#1{\mathchoice
   {\XXint\displaystyle\textstyle{#1}}%
   {\XXint\textstyle\scriptstyle{#1}}%
   {\XXint\scriptstyle\scriptscriptstyle{#1}}%
   {\XXint\scriptscriptstyle\scriptscriptstyle{#1}}%
   \!\int}

\def\XXint#1#2#3{{\setbox0=\hbox{$#1{#2#3}{\int}$}
     \vcenter{\hbox{$#2#3$}}\kern-.5\wd0}}

\newcommand{\cha}{1\!\!1}

\newcommand{\na}{\nabla}
\newcommand{\OO}{\Omega}

\newcommand{\Vb}{\mathbb{V}}

\newcommand{\ui}{\tilde{u}}
\newcommand{\ai}{\tilde{a}}
\newcommand{\bi}{\tilde{b}}

\newcommand{\pssi}{\varphi}

\newcommand{\lt}{\left}
\newcommand{\rt}{\right}
\newcommand{\nl}{\newline}

\newcommand{\nn}{\nonumber}

\newcommand{\lm}{\lambda}

\newcommand{\qd}{\quad}
\newcommand{\spt}{\mathrm{Spt}}
\newcommand{\veps}{\varepsilon}
\newcommand{\vthe}{\vartheta}
\newcommand{\vs}{\varsigma}

\newcommand{\vp}{\varphi}
\newcommand{\ep}{\epsilon}
\newcommand{\thea}{\varphi}

\newcommand{\wt}{\widetilde}

\newcommand{\mm}{\ti{m}}

\newcommand{\AII}{\mathcal{A}}

\newcommand{\GI}{\mathcal{G}}
\newcommand{\EI}{\mathcal{E}}

\newcommand{\MI}{\mathcal{M}}

\newcommand{\KI}{\mathcal{K}}
\newcommand{\DI}{\mathcal{D}}
\newcommand{\CI}{\mathcal{C}}

\newcommand{\UI}{\mathcal{U}}
\newcommand{\WI}{\mathcal{W}}

\newcommand{\BI}{\mathcal{B}}

\newcommand{\ti}{\tilde}
\newcommand{\la}{\langle}

\newcommand{\mi}{\tilde{m}}
\newcommand{\ra}{\rangle}

\newcommand{\R}{\mathrm {I\!R}}

\newcommand{\ca}[1]{\mathrm{Card}\lt(#1\rt)}
\newcommand{\dia}{\diamondsuit}

\newtheorem{a1}{Lemma}
\newtheorem{a2}{Theorem}

\newtheorem{a5}{Proposition}
\newtheorem{a6}{Corollary}

\theoremstyle{remark}

\begin{document}
\title[A quantitative characterisation  of functions with low Aviles Giga]
{A quantitative characterisation of functions with low Aviles Giga energy on convex domains}
\author{Andrew Lorent}
\address{Mathematics Department\\University of Cincinnati\\ 2600 Clifton Ave.\\ Cincinnati\\ Ohio 45221}
\email{lorentaw@uc.edu}

\maketitle

\begin{abstract}
Given a connected Lipschitz domain $\Omega$ we let $\Lambda(\Omega)$ be the subset of functions
in $W^{2,2}(\Omega)$ with $u=0$ on $\partial \Omega$ and whose gradient (in the sense of trace) satisfies $\na u(x)\cdot \eta_x=1$ where 
$\eta_x$ is the inward pointing unit normal to $\partial \Omega$ at $x$.
The functional $I_{\ep}(u)=\frac{1}{2}\int_{\Omega} \ep^{-1}\lt|1-\lt|\na u\rt|^2\rt|^2+\ep\lt|\na^2 u\rt|^2 dz$
minimised over $\Lambda(\Omega)$ serves as a model in connection with problems in liquid crystals and thin film blisters, it is also
the most natural higher order generalisation of the Modica Mortola functional. In \cite{otjab} Jabin, Otto, Perthame characterised a class of functions which includes all limits of sequences $u_n\in\Lambda\lt(\Omega\rt)$ with $I_{\ep_n}(u_n)\rightarrow 0$
as $\ep_n\rightarrow 0$. A corollary to their work is that if there exists such a sequence $(u_n)$ for a bounded
domain $\Omega$, then $\Omega$ must be a ball and (up to change of sign) $u:=\lim_{n\rightarrow \infty} u_n
=\mathrm{dist}(\cdot,\partial\Omega)$. We prove a quantitative generalisation of this corollary for the class of bounded convex sets.

There exists positive constant $\gamma_1$ such that if $\Omega$ is a convex set of diameter $2$ and 
$u\in \Lambda(\Omega)$ with $I_{\ep}(u)=\beta$ then $\lt|B_1(x)\triangle \Omega\rt|\leq
c\beta^{\gamma_1}$ for some $x$ and
$$
\int_{\Omega} \lt|\na u(z)+\frac{z-x}{\lt|z-x\rt|} \rt|^2 dz\leq c\beta^{\gamma_1}.
$$

A corollary of this result is that there exists positive constant $\gamma_2<\gamma_1$ such that if 
$\Omega$ is convex with diameter $2$ and $C^2$ boundary with curvature bounded by $\ep^{-\frac{1}{2}}$, then for any minimiser $v$ of $I_{\ep}$ over $\Lambda(\Omega)$, 
$$
\|v-\zeta\|_{W^{1,2}(\Omega)}\leq c  (\ep+\inf_{y}\lt|\Omega\triangle B_1(y)\rt|)^{\gamma_2}
$$
where $\zeta(z)=\mathrm{dist}(z,\partial \Omega)$. Neither of the constants $\gamma_1$ or $\gamma_2$ are 
optimal. 

\end{abstract}

\section{Introduction}
\label{sec1}

We consider the following functional 
\begin{equation}
\label{eqe1}
I_{\ep}(u)=\frac{1}{2}\int_{\Omega} \ep^{-1}\lt|1-\lt|\na u\rt|^2\rt|^2+\ep\lt|\na^2 u\rt|^2 dz
\end{equation}
the study of which arises from a number of sources, one of the earliest and most important
is the article by Aviles, Giga \cite{avg1}. We will refer to the quantity $I_{\ep}(u)$ as the Aviles-Giga
energy of functional $u$.
Functional $I_{\ep}$ is usually minimised over the space of functions $u\in W^{2,2}(\Omega)$ where
$u(x)=0$ and $\na u(x)\cdot\eta_x=1$ on $\partial\Omega$ (in the sense of trace) where $\eta_x$ is the inward pointing
unit normal, we will denote this space of functions by $\Lambda(\Omega)$.

Aviles, Giga raised the problem of the study of the limiting behavior of $I_{\ep}$ as $\ep\rightarrow 0$
in connection with the theory of smectic liquid crystals \cite{avg1}. In \cite{gior} Gioia, Ortiz studied
$I_{\ep}$ as a model for thin film blisters. Jin, Kohn \cite{jin} introduced the by now classic method of 
estimating the energy by `divergence of vectorfields'. 
A related functional arising from micromagnetics was studied by Riviere, Serfaty \cite{riv2}, in this
case the functional acts on vector fields $m$ (in two dimensions) satisfying $\lt|m\rt|=1$ in $\Omega$
and the functional is given by $M_{\ep}(m)=\ep\int_{\Omega} \lt|\na m\rt|^2
+\ep^{-1}\int_{\R^2} \lt|\na^{-1} \mathrm{div} \ti{m}\rt|^2$ where $\ti{m}$ is vectorfield $m$ extended
trivially by $0$ outside $\Omega$. For the Aviles Giga functional we minimise
over curl free vector fields and the functional forces the norm of the vector field to be close to $1$ with weighting
$\ep^{-1}$ while constraining an $\ep$ multiple of the $L^2$ norm (squared) of the gradient, on the other hand
the micromagnetics functional is minimised over vectorfields whose norm is taken to be $1$ from the outset and
the functional forces the vector field to be divergence free with weighting $\ep^{-1}$
\footnote{the term $\int_{R^2} \lt|\na^{-1} \mathrm{div} m\rt|^2$ is the $L^2$ norm of the Hodge projection
onto curl free vector fields} while again
constraining an $\ep$ multiple of the $L^2$ norm (squared) of the gradient. Functional $M_{\ep}$ is much more
rigid and very much stronger results are known for it than for $I_{\ep}$, see \cite{riv3},\cite{riv2},\cite{ambkir1}, \cite{amler}.

Roughly speaking, the conjecture is that as $\ep\rightarrow 0$ the energy of minimisers of $I_{\ep}$
will converge to a collection of curves on which the gradient of the minimisers make a jump of order $O(1)$
perpendicularly across the curve. This has already been proved for functional $M_{\ep}$ \cite{riv2}. A way to think
about this is the following, given a connected Lipschitz domain $\Omega$ let $w$ be the distance from
$\partial \Omega$ and let $v_{\ep}$ be $w$ convolved by a convolution kernel of diameter $\ep$, the regions
where $\lt|\na v_{\ep}\rt|\not\sim 1$ will be exactly the $\ep$ neighborhoods of the curves on which $\na w$ has a
jump discontinuity. If $\Omega$ is a ball $\na w$ will have a discontinuity only at one point, in
all other cases there will be non trivial curves of singularities and for the specific function
$v_{\ep}$, it is exactly in an $\ep$ neighborhood of these curves that the energy
will concentrate. The conjecture is that what we can observe directly for $v_{\ep}$ will 
hold true for the minimisers of $I_{\ep}$.

The most natural way to study
these questions is within the frame work of $\Gamma$ convergence. One of the earliest successes of $\Gamma$
convergence was the characterisation of the $\Gamma$ limit of the so called Modica Mortola functional
$A_{\ep}(w)=\int_{\Omega} \ep\lt|\na w\rt|^2+\ep^{-1}\lt|1-\lt|w\rt|^2\rt|^2$ which is minimised
over scalar functions $w$ satisfying an integral condition of the form $\int_{\Omega} w=0$. It was
shown by Modica, Mortola \cite{mod} (confirming a conjecture of DeGiorgi)
that the $\Gamma$ limit of $A_{\ep}$ is a constant multiple of the $H^{n-1}$ measure of the
jump set $J_w$ minimised over the space of functions
$w\in\lt\{v\in BV:v\in\lt\{1,-1\rt\}\;a.e.\text{ and }\int v=0\rt\}$. Given the elementary inequality
\begin{equation}
\label{eqd1}
\ep\lt|\na w\rt|^2+\ep^{-1}\lt|1-\lt|w\rt|^2\rt|^2\geq
\lt|\na w\rt|\lt|1-\lt|w\rt|^2\rt|
\end{equation}
we have that for any sequence $(w_n)$ of equibounded $A_{\ep_n}$ energy (for some
subsequence $\ep_n\rightarrow 0$) has a uniform $L^1$ control of $\na \lt(w_n-\frac{w_n^3}{3}\rt)$ and the measure we
obtain as the limit of this $L^1$ sequence of gradients will naturally be supported on the jump set of the
limiting function. In some sense the nature of the $\Gamma$ limit of $A_{\ep}$ could be anticipated from 
(\ref{eqd1}).

Functional $I_{\ep}$ is the most natural higher order
generalisation of $A_{\ep}$, in the case of $I_{\ep}$ the conjectured $\Gamma$ limit is
surprising, this is part of the reason that functional $I_{\ep}$ has received so much attention. The
first works on identifying the $\Gamma$ limit are by Aviles, Giga \cite{avg1} and Jin, Kohn \cite{jin}, 
later these
ideas were developed by Ambrosio, DeLellis, Mantegazza \cite{amb2}, roughly speaking the limiting function space
is conjectured to have a structure similar to the space of functions whose gradient is $BV$ and
the limiting energy is conjectured to have the form $\int_{J_{\na u}} \lt|\na u^{+}-\na u^{-}\rt|^3 dH^1$.
Much progress has been made on this conjecture, particularly equi-coercivity of $I_{\ep}$ has been
shown independently in \cite{amb2} and in the work of Desimone, Kohn, Muller, Otto \cite{otto1}.
A proposed limiting function space $AG(\Omega)$ and limiting functional $I$
as been suggested in \cite{amb2} and it was shown
that all limits of sequences of functions $(u_n)$ with $\sup_{n} I_{\ep_n}(u_n)<\infty$
are such that $u_n\overset{W^{1,3}}{\rightarrow} u\in AG(\Omega)$ and $\liminf I_{\ep_n}(\na u_n)\geq
I(u)$. The compactness proofs provided by \cite{otto1} and  \cite{amb2} are different but
share some common ideas. The proof by \cite{otto1} identifies the set
of all smooth functions $\Phi:\R^2\rightarrow \R^2$ for which there exists smooth $\Psi:\R^2\rightarrow \R^2$ such that 
\begin{equation}
\label{eqd1.6}
\int \lt|\mathrm{div}\lt[\Phi(\na u)\rt]\rt|\leq c\int \lt|\Psi(\na u)\cdot \na \lt(1-\lt|\na u\rt|^2\rt)\rt|\text{ for any }C^2\text{ function  }u,
\end{equation}
influenced by ideas of Tartar and Murat on compensated compactness \cite{tar} \cite{mur} the authors are able to
prove that this set of $\Phi$ is sufficiently rich so as to force $\na u_n$ to converge strongly. In
\cite{avg1} the authors (building on work of Jin Kohn \cite{jin}) found two third order
polynomial vector fields $\Sigma_1:\R^2\rightarrow \R^2$ and $\Sigma_2:\R^2\rightarrow \R^2$ such that
\begin{equation}
\label{eqd20}
\int \lt|\mathrm{div}\lt[\Sigma_i(\na u)\rt]\rt|\leq c\int \lt|\na^2 u\rt|\lt|1-\lt|\na u\rt|^2\rt|
\text{ for any }C^2\text{ function  }u,\text{ for }i=1,2.
\end{equation}
Using some elementary and surprising identities satisfied by $\Sigma_1(\na u), \Sigma_2(\na u)$ a different approach to
compactness was found. Rather naturally considering (\ref{eqd20}), the function space $AG(\Omega)$ proposed
by \cite{amb2} is given by the set of functions $v$ for which $\mathrm{div}(\Sigma_i(\na v))$ forms a Radon measure for $i=1,2$
and the limiting energy functional $I(v)$ is given by the total absolute value of this measure on $\Omega$.

Given vector field $w$ let $\chi(\xi,w):=\cha_{\lt\{\xi\cdot w>0\rt\}}$, Jabin, Perthame \cite{jab2} showed
that gradients of sequences of bounded Aviles-Giga energy (in fact their method extends to more general
functionals) are compact and the limit $\na u$ satisfies a kinetic equation of the form $\xi\cdot \na_x \chi(\xi,R(\na u))=q$
where $q$ is the distribution derivative with respect to $\xi$ of some measure on $\R^2_{\xi}\times \R^2_{x}$ and $R$ is the rotation given by $R(x,y)=(-y,x)$.
By application of kinetic averaging lemmas \cite{lion1} this
leads to some regularity; $\na u\in W^{s,q}$ for all $0\leq s<\frac{1}{5}$,
$q<\frac{5}{3}$ and using the kinetic equation a different proof of compactness was found.
The kinetic equation deduced by \cite{jab2} was motivated by the characterisation of the set of $\Phi$
satisfying (\ref{eqd1.6}) given in \cite{otto1}, indeed defining $\ti{\Phi}(z)=\lt|z\rt|^2 e$ for $z\cdot e>0$
and $0$ otherwise,  in \cite{otto1} it was shown that a sequence $\Phi_n$ satisfying (\ref{eqd1.6}) could
be found that approximates $\ti{\Phi}$ pointwise. Using the kinetic equation deduced in \cite{jab2},
Jabin, Otto, Perthame \cite{otjab} were able to characterise zero energy limits (and the domains that
allow them) for $I_{\ep}$, in fact their result is stronger, they showed that if a divergence free
vector field $m$ satisfies the kinetic equation $\xi\cdot \na \chi(m,\xi)=0$, $\lt|m(x)\rt|=1$ a.e.\ in
$\Omega$ and $m(x)\cdot \eta_x=0$ on $\partial \Omega$ then either $\Omega$ is a strip and $m$ is a
constant or $\Omega=B_r(x)$ for some $r>0$, $x\in \R^2$ and $m(z)=\lt(\frac{z-x}{\lt|z-x\rt|}\rt)^{\perp}$ or $m(z)=-\lt(\frac{z-x}{\lt|z-x\rt|}\rt)^{\perp}$. An analogous result for zero energy limits of $M_{\ep}$ is stated in  \cite{leri} and is a consequence of
the main theorem of \cite{amler}.

As a corollary, given a sequence $u_n\in \Lambda(\Omega)$ and $\ep_n\rightarrow 0$ such that $I_{\ep_n}(u_n)\rightarrow 0$
as $n\rightarrow \infty$, letting $u$ be the limit of this sequence, the vector field $R(\na u)$ satisfies
the hypothesis stated and hence we have (up to a sign) a complete description of $\na u$.

The main theorem of this paper is a quantitative generalisation of the corollary to Jabin, Otto, Perthame
theorem over the class of bounded convex sets.
\begin{a2}
\label{T0}
Let $\ep>0$ and $\Omega$ be a convex domain with diameter $2$. Let $u\in W^{2,2}(\Omega)$
with $u=0$ on $\partial \Omega$ and 
$\na u(x)\cdot \eta_x=1$ of $\partial \Omega$ (in the sense of trace) where $\eta_x$ is the inward
pointing unit normal. Then there exists positive constants $\CI>1$ and $\gamma<1$ such that for some 
$x\in\Omega$, 
$$
\lt|\Omega\triangle B_1(x)\rt|\leq \CI\lt(I_{\ep}(u)\rt)^{\gamma}
$$
and
$$
\int_{\Omega} \lt|\na u(z)+\frac{z-x}{\lt|z-x\rt|}\rt|^2 dz\leq \CI\lt(I_{\ep}(u)\rt)^{\gamma}.
$$
\end{a2}
\begin{a6}
\label{CC1}
Let $\ep>0$ and $\Omega$ be a convex set of diameter $2$ and with $C^2$ boundary and curvature bounded above by $\ep^{-\frac{1}{2}}$. Let $\Lambda(\Omega):=\lt\{u\in W^{2,2}(\Omega):u=0\text{ on }\partial \Omega\text{ and }\na u(z)\cdot \eta_z=1\text{ for }z\in \partial \Omega\rt\}$. There exists 
positive constants $\CI=\CI(\Omega)>1$ and $\lambda<1$ such that if $u$ is a minimiser of $I_{\ep}$ over $\Lambda(\Omega)$, then 
\begin{equation}
\label{uz9}
\|u-\zeta\|_{W^{1,2}(\Omega)}\leq \CI \lt(\ep+\inf_{y\in \Omega}\lt|\Omega\triangle B_1(y)\rt|\rt)^{\lambda}
\end{equation}
where $\zeta(z)=\mathrm{dist}(z,\partial \Omega)$. 
\end{a6}

In Theorem \ref{T0} we take $\gamma=512^{-1}$ and in Corollary \ref{CC1}, $\lambda=5462^{-1}$. Neither constant is 
optimal. Corollary \ref{CC1} requires a fair amount of technical work establishing an upper bound for the minimizer of $I_{\ep}$ in terms 
of the `eccentricity' $\inf_{y\in \Omega, r>0} \lt|\Omega\triangle B_r(y)\rt|$. For the reader primarily interested in the asymptotic behavior 
of minimizers as $\ep\rightarrow 0$ recent powerful results on $\Gamma$-convergence upper bound of $I_{\ep}$ 
(in the case where the function $u$ being approximated satisfies  $\na u\in BV(\Omega:S^1)$) by Conti, DeLellis \cite{conti-del} and Poliakovsky \cite{poli} do much 
of the work for us and we can give a relatively shorter proof of the following corollary to Theorem \ref{T0}. Note that Corollary \ref{CC3} stated 
below is a corollary to Corollary \ref{CC1}.

\begin{a6}
\label{CC3}
Let $\Omega$ be a convex set of diameter $2$ with $C^2$ boundary. Let $\Lambda(\Omega)$ be 
as defined in Corollary \ref{CC1}. There exists 
positive constants $\CI=\CI(\Omega)>1$ and $\lambda<1$ such that if $u^{\ep}$ is a minimiser of $I_{\ep}$ over $\Lambda(\Omega)$, then 
\begin{equation}
\label{uz9.8}
\limsup_{\ep\rightarrow 0}\|u^{\ep}-\zeta\|_{W^{1,2}(\Omega)}\leq \CI \lt(\inf_{y\in \Omega}\lt|\Omega\triangle B_1(y)\rt|\rt)^{\lambda}
\end{equation}
where $\zeta(z)=\mathrm{dist}(z,\partial \Omega)$. 
\end{a6}
\bf Plan of paper. \rm After the introduction in Section \ref{sec1} we sketch the proof 
of the main theorem in Section \ref{sec2}. In Section \ref{sec3} we prove the main theorem. In Section \ref{sec4} we establish Corollary \ref{CC3}, the additional lemmas needed to establish Corollary \ref{CC1} are given in Section \ref{sec5}.  

\subsection{Background}

Given a sequence $\ep_n\rightarrow 0$ and $u_n\in \Lambda(\Omega)$ with $\limsup I_{\ep_n}(u_n)<\infty$, let
$u$ be the limit of $u_n$, the vector valued measure given by
$\nu_u:=(\mathrm{div}\lt[\Sigma_1(\na u)\rt],\mathrm{div}\lt[\Sigma_2(\na u)\rt])$ (where $\Sigma_1,\Sigma_2$ are the third order
polynomial vector fields that satisfy (\ref{eqd20})) gives us the expression of the limiting energy, i.e.\
$I(u)=\|\nu_u\|(\Omega)$. If we consider
the $1$-dimensional part of the measure
$$
\Gamma:=\lt\{x:\limsup_{r\rightarrow 0} \frac{\|\nu_u(B_r(x))\|}{r}>0\rt\}
$$
it has been shown that $\Gamma$ is $1$-rectifiable \cite{delo} (see also \cite{delo2})
and an analogous result has been shown for $M_{\ep}$ \cite{ambkir1}. It was also
shown $\na u$ has jump discontinuous across the rectifiable set $\Gamma$ exactly as would be the
case if $\na u$ was $BV$ and its jump set was given by $\Gamma$. However it is not known (even if $u_n$ are the minimisers of $I_{\ep_n}$) if measure $\|\nu_u\|$ is even singular with respect to Lebesgue measure. Note that 
for the function $M_{\ep}$ the minimiser of the limiting energy is known to be rectifiable 
\cite{amler}, for a sequence with only equibounded energy the measure is not known to be singular.

The original motivation for Theorem \ref{T0} was to prove a version of it for $\Omega=B_1(0)$ without
boundary conditions, under the hypotheses $\int_{B_1} \lt|1-\lt|\na u\rt|^2\rt|\lt|\na^2 u\rt| dz=\beta$,
$\int_{B_1} \lt|1-\lt|\na u\rt|^2\rt| dz\leq \ep$ and
$\sup\lt\{\|u-A\|_{L^{\infty}\lt(B_1\lt(0\rt)\rt)}:A\text{ is affine with }\lt|\na A\rt|=1\rt\}\leq 1000^{-1}$,
the conclusion in this case would be that there exists a smooth function $\psi$ with $\lt|\na \psi\rt|=1$
everywhere such that $\|\na u-\na \psi\|_{L^2\lt(B_{2^{-1}}\lt(0\rt)\rt)}\leq c\beta^{\gamma}$ for some
$\gamma>0$. This is a kind of quantitative version of the main proposition required to prove compactness in 
\cite{amb2}, (see Proposition 4.6). The hope is to use such a quantitative result to show $\|\nu_u\|$ is singular, or at least that $\na u$ is continuous at $H^1$ a.e. point outside $\Gamma$, we will address these issues in a forthcoming
paper \cite{lorft}.

The many strong results about measure $\|\nu_u\|$ (and the measure that gives the limiting functional
for the micromagnetics function) have been achieved by characterising various
kinds of \em blow up \rm of the measure and understanding well the absolute (i.e.\ non quantitative)
situation in the limit \cite{ambkir1}, \cite{delo}, \cite{delo2}, \cite{otjab}, \cite{amler}. In some
sense there are only two possibilities, to take a limit and have an absolute situation
and to understand the measure from this, or to stop before the limit and have a non-absolute situation
and try and understand something about it with a quantitative theorem. Our primary motivation
in proving a quantitative version of Jabin-Otto-Perthame Theorem was so as to obtain
a result that could be used for the latter approach.

By Poincare's inequality it is easy to see $\inf_{\Lambda(\Omega)} I_{\ep}\geq c\ep$ and so 
Theorem \ref{T0} follows from the following slightly more general result.
\begin{a2}
\label{T1}
Let $\Omega$ be a convex body centered on $0$ with $\mathrm{diam}(\Omega)=2$. Let $\beta>0$, suppose
$u:W^{2,2}(\Omega)\rightarrow \R$ is a function satisfying
\begin{equation}
\label{eq1}
\int_{\Omega} \lt|1-\lt|\na u\rt|^2\rt|\lt|\na^2 u\rt| dz\leq \beta
\end{equation}
and
\begin{equation}
\label{eq2}
\int_{\Omega} \lt|1-\lt|\na u\rt|^2\rt|^2 dz\leq \beta^2
\end{equation}
and in addition $u$ satisfies $u=0$ on $\partial \Omega$ and $\na u(z)\cdot \eta_z=1$ on $\partial \Omega$ in the sense of trace 
where $\eta_z$ is the inward pointing unit normal to $\partial\Omega$ at $z$.

Then there exists positive constant $\CI_1>0$ such that
$\lt|B_{1}\lt(0\rt)\triangle\Omega\rt|< \CI_1 \beta^{\frac{1}{512}}$ and 
\begin{equation}
\label{eq70.6}
\int_{\Omega} \lt|\na u(z)+\frac{z}{\lt|z\rt|}\rt|^2 dz\leq \CI_1\beta^{\frac{1}{512}}.
\end{equation}
\end{a2}

\bf Acknowledgments. \rm  Part of this paper was written while the author was 
the Emma e Giovanni Sansone Junior Visitor at Centro
di Ricerca Matematica Ennio De Giorgi, Pisa. 
The hospitality and support this institute is gratefully acknowledged. I would also like to express my great thanks to the referee for numerous suggestions, simplifications and 
improvements. The quality of the paper has been substantially increased by the input of the referee.

\section{Sketch of the proof}
\label{sec2}

\subsection{Sketch of the proof of Theorem \ref{T1}}

While the proof for convex domains is slightly involved, there are only a couple
of ideas that are really central. We will sketch the proof for the case $\Omega=B_1(0)$, ignoring
(without comment) many technicalities in order to give an impression of the basic skeleton.

The real engine of the proof is the characterisation in \cite{otto1} of the set of $\Phi$ such that (\ref{eqd1.6}) is satisfied.
As mentioned in the introduction, as consequence of the characterisation it was shown there exists a sequence
of $\Phi_n$ satisfying (\ref{eqd1.6}) that converge pointwise to the function $\ti{\Phi}(z)=\lt|z\rt|^2 e$ for
$z\cdot e>0$ and $0$ otherwise. Following closely the proof of this
it is possible to extract the existence of functions $\Phi_{\theta}$ and $\Psi_{\theta}$ with
$\|\na \Phi_{\theta}\|\leq
c\beta^{-\frac{1}{4}}$, $\|\Psi_{\theta}\|\leq c\beta^{-\frac{1}{4}}$, $\|\na\Psi_{\theta}\|\leq c\beta^{-\frac{1}{2}}$
such that the following two inequalities hold.

\begin{figure}[h]
\centerline{\psfig{figure=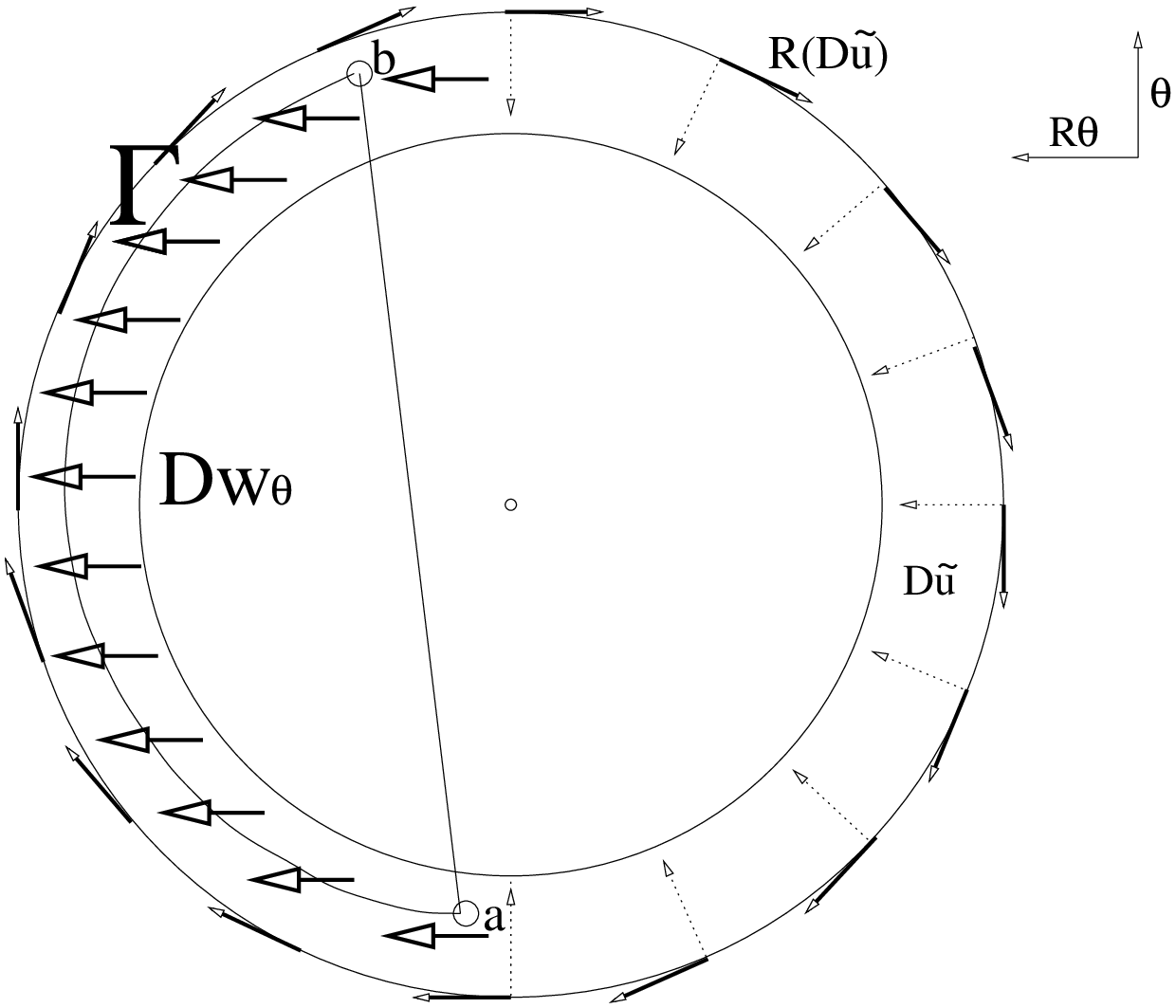,height=9cm,angle=0}}
\caption{}\label{fig10}
\end{figure}

Let $\Lambda_{\theta}(z):=\theta$ for $z\cdot \theta>0$ and $0$ otherwise,
\begin{equation}
\label{eqe1.8}
\lt|\Phi_{\theta}\lt(z\rt)-\Lambda_{\theta}\lt(z\rt)\rt|\leq c\beta^{\frac{1}{4}}
\text{ for }z\in N_{\sqrt{\beta}}(S^1)\backslash B_{2\beta^{\frac{1}{4}}}(\theta)
\end{equation}
and (letting $R(z_1,z_2)=(-z_2,z_1)$ be the anti-clockwise rotation)
\begin{equation}
\label{eqe2}
\mathrm{div}\lt[\Phi_{\theta}\lt(R(\na w)\rt)-\Psi_{\theta}\lt(R(\na w)\rt)
\lt(1-\lt|R(\na w)\rt|^2\rt)\rt]\leq c\beta^{-\frac{1}{2}}\lt|1-\lt|\na w\rt|^2\rt|\lt|\na^2 w\rt|
\text{ for any }w\in W^{2,1}.
\end{equation}
Recall, for simplicity we have taken $\Omega=B_1(0)$, as $\na u(z)=-\frac{z}{\lt|z\rt|}$ on $\partial B_1(0)$ then 
we can extend $u$ to a function $\ui:B_{11/10}(0)\rightarrow \R$ such that
$$
\int_{B_{11/10}(0)} \lt|1-\lt|\na \ui\rt|^2\rt|\lt|\na^2 \ui\rt| dz\leq c\beta, 
\int_{B_{11/10}(0)} \lt|1-\lt|\na \ui\rt|^2\rt|^2 dz\leq c\beta^2
$$ 
and
\begin{equation}
\label{eqe3}
\na\ui(z)=-\frac{z}{\lt|z\rt|}\text{ for any }z\in B_{11/10}(0).
\end{equation}
It is more convenient to work with vectorfields that are \em almost \rm curl free instead of
\em almost \rm divergence free. So notice that (\ref{eqe1.8}) can be rewritten as
\begin{equation}
\label{eqe4}
\lt|R\lt(\Phi_{\theta}\lt(z\rt)\rt)-R\lt(\Lambda_{\theta}\lt(z\rt)\rt)\rt|\leq c\beta^{\frac{1}{4}}
\text{ for }z\in N_{\sqrt{\beta}}(S^1)\backslash B_{2\beta^{\frac{1}{4}}}(\theta)
\end{equation}
and we have
$\int_{B_{11/10}(0)}\lt|\mathrm{curl}\lt[R\lt(\Phi_{\theta}\lt(R\lt(\na \ui\rt)\rt)\rt)-R\lt(\Psi_{\theta}\lt(R\lt(\na\ui\rt)\rt)\rt)
\lt(1-\lt|\na \ui\rt|^2\rt)\rt]\rt|\leq c\sqrt{\beta}$.
By the quantitative Hodge decomposition type theorem from \cite{amb2} (Theorem 4.3) we can find a scalar valued
function $w_{\theta}$ such that
\begin{equation}
\label{eqe5}
\int_{B_{11/10}(0)}\lt|\na w_{\theta}-\lt(R\lt(\Phi_{\theta}\lt(R\lt(\na \ui\rt)\rt)\rt)
-R\lt(\Psi_{\theta}\lt(\na R\lt(\na\ui\rt)\rt)\rt)
\lt(1-\lt|\na \ui\rt|^2\rt)\rt) \rt| dz\leq c\sqrt{\beta}.
\end{equation}
The real power of (\ref{eqe5}) is that on the annulus $\AII:=B_{11/10}(0)\backslash B_{1}(0)$ we
know that $\na\ui(z)=-\frac{z}{\lt|z\rt|}$ and hence given inequality (\ref{eqe4}) (and the fact that
$\lt|\na \ui\rt|=1$ on $\AII$) we have a that $\Phi_{\theta}\lt(R\lt(\na \ui(z)\rt)\rt)\in
N_{\beta^{\frac{1}{4}}}\lt(\theta\rt)$ for any $z\in \AII\cap H\lt(R\theta,0\rt)$, see figure \ref{fig10}.

In much the same way in the ball $B_1(0)$, by inequalities (\ref{eqe4}), (\ref{eqe5}) and $\int_{B_1(0)}
\lt|1-\lt|\na \ui\rt|^2\rt|^2\leq \beta^2$ we have that there exists a large set $\GI\subset B_1(0)\cap
H(0,R\theta)$, with $\lt| B_1(0)\backslash \GI\rt|\leq \sqrt{\beta}$ such that if $z\in\GI$ then
$\na w_{\theta}(z)\in B_{\beta^{\frac{1}{4}}}(R\theta)$ or $\na w_{\theta}(z)\in B_{\beta{\frac{1}{4}}}(0)$
depending on whether $R(\na u(z))\cdot \theta>0$ or $R(\na u(z))\cdot \theta\leq 0$.

It is not hard to see we can find points $a,b\in N_{\beta^{\frac{1}{8}}}(\la \theta\ra\cap \partial B_1(0))$ with
$\lt|a-b\rt|\sim 2$,  $\theta\cdot \frac{b-a}{\lt|b-a\rt|}>0$, the angle between $\frac{b-a}{\lt|b-a\rt|}$ and $\theta$ is at least
$\beta^{\frac{1}{8}}$ and $H^1(\lt[a,b\rt]\backslash \GI)\leq \beta^{\frac{1}{4}}$. Let $\GI_1=\lt\{x\in\GI:\na u(z)\cdot R^{-1}\lt(\theta\rt)>0\rt\}$ and $\GI_2=\GI\backslash \GI_1$. As can be seen
from figure \ref{fig10} we can connect $a$ to $b$ with a path $\Gamma\subset \AII$ so
\begin{eqnarray}
\label{eqe40}
\lt|w_{\theta}(b)-w_{\theta}(a)\rt|&=&\lt|\int_{\Gamma} \na w_{\theta}(z)t_z dH^1 z\rt|\geq \lt|
R\theta\cdot \lt(\int_{\Gamma} t_z dH^1 z\rt)\rt|-c\beta^{\frac{1}{4}}\nn\\
&=&\lt|R\theta\cdot \frac{b-a}{\lt|b-a\rt|}\rt|\lt|b-a\rt|-c\beta^{\frac{1}{4}}.
\end{eqnarray}
On the other hand
\begin{eqnarray}
\label{eqe41}
\lt|w_{\theta}(b)-w_{\theta}(a)\rt|&=&\lt|\int_{\lt[a,b\rt]} \na w_{\theta}(z)\frac{b-a}{\lt|b-a\rt|} dH^1 z\rt|
\leq \lt|\int_{\lt[a,b\rt]\cap \GI_1} \na w_{\theta}(z)\frac{b-a}{\lt|b-a\rt|} dH^1 z\rt|+c\beta^{\frac{1}{4}}\nn\\
&\leq&\lt|\int_{\lt[a,b\rt]\cap \GI_1} R\theta\cdot\frac{b-a}{\lt|b-a\rt|} dH^1 z\rt|
+c\beta^{\frac{1}{4}}\nn\\
&=& \lt|R\theta\cdot\frac{b-a}{\lt|b-a\rt|}\rt|H^1(\lt[a,b\rt]\cap \GI_1)+c\beta^{\frac{1}{4}}
\end{eqnarray}
and since $\lt|R\theta\cdot\frac{b-a}{\lt|b-a\rt|}\rt|\geq \beta^{\frac{1}{8}}$ so putting
(\ref{eqe40}) and (\ref{eqe41}) together
$$
\lt|a-b\rt|\leq H^1\lt(\lt[a,b\rt]\cap \GI_1\rt)+\frac{c\beta^{\frac{1}{4}}}{\lt|R\theta\cdot\frac{b-a}{\lt|b-a\rt|}\rt|}
\leq H^1\lt(\lt[a,b\rt]\cap \GI_1\rt)+c\beta^{\frac{1}{8}}.
$$
So by arguing in the same way for lines parallel to $\lt[a,b\rt]$ by Fubini's theorem we can show
$\lt|H\lt(\frac{a+b}{2},R\lt(\frac{b-a}{\lt|b-a\rt|}\rt)\rt)\backslash \GI_1\rt|\leq c\beta^{\frac{1}{8}}$.
Thus all but $\beta^{\frac{1}{8}}$ points $z\in B_1(0)\cap H(0,R(\theta))$ are such that
$\na u(z)\cdot R^{-1}(\theta)>0$.  As $\theta$ is arbitrary we can rephrase this the following way.
Given $\phi\in S^1$ for all but $\beta^{\frac{1}{8}}$ points $z\in B_1(0)\cap H(0,\phi)$ are such that
$\na u(z)\cdot (-\phi)>0$.

Now take $\psi=\lt(\begin{matrix}  \cos \beta^{\frac{1}{16}}\\ \sin\beta^{\frac{1}{16}}\end{matrix}\rt)$. For all
but $\beta^{\frac{1}{8}}$ points in $H(0,e_1)\cap H(0,-\psi)\cap H(0,-e_2)$ we have that
$\na u(z)\cdot (-e_1)>0$ and $\na u(z)\cdot \psi>0$, it is not hard to show this implies $\lt|\na u(z)\cdot e_1\rt|\leq
c\beta^{\frac{1}{16}}$ and since $\na u(z)\cdot e_2>0$ and $\lt|\na u(z)\rt|\sim 1$ we have $\na u(z)\in
B_{c\beta^{\frac{1}{16}}}(e_2)$ with an exceptional set of measure less than $c\beta^{\frac{1}{8}}$.
So integrating a carefully chosen line inside
$H(0,e_1)\cap H(0,-\psi)\cap H(0,-e_2)$ and using the fact that $u=0$ on $\partial B_1(0)$ we can show
$\lt|u(0)-1\rt|\leq c\beta^{\frac{1}{16}}$.

Now since $\lt|\na u\rt|$ is mostly very close to $1$ and we have zero boundary condition, so avoiding technicalities
assuming the coarea formula we have $\int_{\theta\in S^1} \int_{\R_{+}\theta\cap B_1(0)} \lt|\lt|\na u(z)\rt|^2-1\rt| dH^1 z dH^1 \theta\leq c\sqrt{\beta}$. Note also that for any $\theta\in S^1$, $u(\theta)=0$ so by the fundamental theorem of 
Calculus 
\begin{eqnarray}
\lt|\int_{\R_{+}\theta\cap B_1(0)} \na u(z)\cdot (-\theta) dH^1 z-1\rt|&\leq&\lt|(u(0)-u(\theta))-1\rt|\nn\\
&\leq&c\beta^{\frac{1}{16}}\nn
\end{eqnarray}
so
\begin{eqnarray}
&~&\int_{\theta\in S^1} \int_{\R_{+}\theta\cap B_1(0)} \lt|\na u(z)+\theta\rt|^2 dH^1 z \theta\nn\\
&~&\qd\qd\qd\qd\qd=
\int_{\theta\in S^1} \int_{\R_{+}\theta\cap B_1(0)} \lt|\na u(z)\rt|^2+2\na u(z)\cdot \theta+\lt|\theta\rt|^2 dH^1 z dH^1 \theta\nn\\
&~&\qd\qd\qd\qd\qd \leq c\beta^{\frac{1}{16}}.
\end{eqnarray}
This concludes the sketch of the proof of Theorem \ref{T1}. 

%

\subsection{Sketch of the proof of Corollary \ref{CC1} and Corollary \ref{CC3}}

In order to deduce Corollary \ref{CC1} we need to apply Theorem \ref{T0} to the minimizer 
of $I_{\ep}$ over $\Lambda(\Omega)$. We can only do this if the minimizer has small energy (and from 
Theorem \ref{T0} we know it can only have small energy if $\Omega$ is close to a ball). For this reason it is 
necessary to construct a function in $\Lambda(\Omega)$ with this property. 
It turns out this is a surprisingly delicate task, it is achieved in 
Section \ref{sec4} and Section \ref{sec5} of the paper.  

The obvious way to attempt the construction is to make some adaption of the function 
$\zeta(z)=\mathrm{dist}(z,\partial \Omega)$, this function clearly satisfies the correct boundary condition. The 
first problem is that $\na \zeta$ will have its gradient in BV and it is easy to construct examples of convex domains that are 
close to balls for which the singular part of $\na \zeta$ is widely spread over the domain. So 
it is necessary to convolve $\zeta$, let $\psi$ denote the convolution of 
$\zeta$ with a convolution kernel of support size $\sim \ep$.

We need to check that the function $\psi$ we obtain 
by convolving $\zeta$ will have small energy. By recent results of \cite{amdel} we 
have that $\na \zeta\in SBV(\Omega:S^1)$. So by Poincare inequality if for most balls the gradient of $\na \zeta$ is not too concentrated in balls of sized $\ep$ then 
we would have $\int_{\Omega} \lt|1-\lt|\na \psi\rt|^2\rt|^2 dz$ is small. Now 
assuming $\Omega$ is close to a ball, then for $x$ not too close to the center of 
$\Omega$ (which we assume is $0$) it is not hard to show that $\lt|\na\zeta(z)+\frac{z}{\lt|z\rt|}\rt|$ is small. By convexity of $\Omega$, if $\Phi^t$ is a 
parameterization of $\zeta^{-1}(t)$ then 
$h\rightarrow \na \zeta(\Phi^t(h))$ will 
be a monotonic parameterization of $S^1$. So the total variation of $\na \zeta$ 
can be explicitly bounded above. The closer $\Omega$ is to a ball the better the estimate on $\lt|\na \zeta(z)+\frac{z}{\lt|z\rt|}\rt|$ holds but near the center it breaks down. To overcome this we do the following. Let 
$\beta=\lt|\Omega\triangle B_1(0)\rt|$ and let 
$\eta(z):=1-\beta^{\frac{3}{32}}+\lt|z\rt|$, so $\Pi:=\lt\{z:\eta(z)\leq\zeta(z)\rt\}$ is roughly a ball centered on $0$ of radius $\beta^{\frac{3}{32}}$. 
So defining $w:=\min\lt\{\zeta,\eta\rt\}$ we have 
$\lt|\na w\rt|=1$ a.e.\ and $\na w\in SBV$. Notice that $\int_{J_{\na w}\cap \Omega} \lt|\na w^{+}-\na w^{-}\rt|^3 dH^1\leq \int_{J_{\na \zeta}\backslash \Pi} \lt|\na \zeta^{+}-\na \zeta^{-}\rt|^3 dH^1+8 H^1(\Gamma)$. Now $\Pi$ is a convex set of diameter approximately 
$\beta^{\frac{3}{32}}$ so $H^1(\Gamma)\sim \beta^{\frac{3}{32}}$. So we have the estimate $\lt|\na\zeta(z)+\frac{z}{\lt|z\rt|}\rt|\leq c\beta^{\frac{3}{32}}$ so 
$\lt|\na \zeta^{-}(z)-\na \zeta^{+}(z)\rt|\leq c\beta^{\frac{3}{32}}$ for any $z\in J_{\na \zeta}\backslash \Pi$. Now by convexity of $\Omega$ and hence monotonicity 
of the gradient along the level set $\zeta^{-1}(t)$ we can prove an explicit 
upper bound $V(\na \zeta,\Omega\backslash  \Pi)\leq 8\pi$. So 
we can estimate 
\begin{eqnarray}
\int_{J_{\na \zeta}\backslash \Pi} \lt|\na \zeta^{+}-\na \zeta^{-}\rt|^3 dH^1
&\leq& 
\sup_{J_{\na \zeta}\backslash \Pi}\lt|\na \zeta^{+}-\na \zeta^{-}\rt|^2\int_{J_{\na \zeta}\backslash \Pi} \lt|\na \zeta^{+}-\na \zeta^{-}\rt| dH^1   \nn\\
&\leq& \sup_{J_{\na \zeta}\backslash \Pi}\lt|\na \zeta^{+}-\na \zeta^{-}\rt|^2  
V(\na \zeta, \Omega\backslash \Pi)\leq 8\pi\beta^{\frac{3}{16}}.
\end{eqnarray}
Putting these things together we have $\int_{J_{\na w}\cap \Omega} \lt|\na w^{+}-\na w^{-}\rt|^3 dH^1\leq c\beta^{\frac{3}{32}}$. This allows us to apply recent results on $\Gamma$-upper bounds of functions whose 
gradient belongs to $SBV$ by \cite{conti-del}, \cite{poli}. These results give the existence of a sequence $u^{\ep}$ with 
the same boundary conditions as $w$ and with the property that 
$\limsup_{\ep\rightarrow 0} I_{\ep}(u^{\ep})\leq c\beta^{\frac{3}{32}}$. This energy bound allows us to apply Theorem \ref{T0} and hence to establish 
Corollary \ref{CC3}. 

To establish Corollary \ref{CC1} requires us to construct a Sobolev function 
by adapting $w$ with `our own hands'. Function $\psi$ we obtained by convolving 
$\zeta$ has a problem in that the convolution will destroy the boundary condition. To circumvent 
this obstacle, in an $\sqrt{\ep}$ neighborhood of the $\partial \Omega$ we convolve the $\zeta$ with a convolution kernel who support 
decreases in proportion to the distance to the boundary. Let 
the new function be denoted by $\pssi$. We make the assumption that $\partial \Omega$ is 
$C^2$ with curvature bounded above by $\ep^{-\frac{1}{2}}$ and this allows us estimate the various error 
terms involved in differentiating a function that is convolved with a kernel of varying support. Clearly the goal is to show 
that $\int_{\Omega} \ep^{-1}\lt|1-\lt|\na \pssi\rt|^2\rt| dz\leq \beta^{\frac{3}{32}}$ 
and $\ep\int_{\Omega} \lt|\na^2 \pssi\rt|^2 dz\leq \beta^{\frac{3}{32}}$.  Establishing the upper bounds required in $\Omega\backslash \lt(N_{\sqrt{\ep}}(\partial \Omega)\cup N_{\ep}(\Pi)\rt)$ can be achieved by Poincare inequalities and the estimate $V(\Omega\backslash \Pi,\na \zeta)\leq 8\pi$. 
Establishing the upper bounds on $N_{\sqrt{\ep}}(\partial \Omega)$ can be achieved by very precise estimates on $\na \pssi$ and $\na^2 \pssi$ which are made due to the fact 
that the curvature conditions on $\partial \Omega$ implies $\na \zeta$ has no singular points in this neighborhood.  The length of $\partial \Pi$ is less than $c \beta^{\frac{3}{32}}$ so as 
$\|\na \pssi\|_{\infty}<c$ we know $\int_{N_{\ep}(\partial \Pi)} \ep^{-1}\lt|1-\lt|\na \pssi\rt|^2\rt| dz\leq c\beta^{\frac{3}{32}}$. Similarly as for 
$z\in \Omega\backslash N_{\sqrt{\ep}}(\partial \Omega)$, $\|\na^2 \pssi\|_{\infty}\leq c\ep^{-1}$ so $\ep\int_{N_{\ep}(\partial \Pi)} \lt|\na^2 \pssi\rt|^2 dz\leq c\beta^{\frac{3}{32}}$. The energy of $\pssi$ in $\Pi\backslash N_{\ep}(\partial \Pi)$ can easily be estimated and shown to be negligible so putting these 
things together gives that $I_{\ep}(\pssi)\leq c\beta^{\frac{3}{32}}$. This upper bound allows us to apply Theorem \ref{T0} and hence to establish 
Corollary \ref{CC1}.

%
%

\section{Proof of Theorem}
\label{sec3}

It should be re-emphasized that the main calculations that makes this lemma work 
(specifically equation (\ref{eq9})) are very minor adoptions of the calculations in 
\cite{otto1}. 
\begin{a1}
\label{L1}
Let $\Omega$ be a convex body centered on $0$ with $\mathrm{diam}(\Omega)\leq 2$. Suppose
$u:W^{2,1}(\Omega)\rightarrow \R$ satisfies (\ref{eq1}) and (\ref{eq2}).
For each $\theta\in S^1$ define $\Lambda_{\theta}:\R^2\rightarrow S^1$ be defined by
\begin{equation}
\label{eq750}
\Lambda_{\theta}\lt(z\rt)=
\begin{cases} \theta & \text{if $z\cdot \theta>0$,}\\
0 &\text{if $z\cdot\theta\leq 0$.}
\end{cases}
\end{equation}
Let $R\in SO(2)$ be the anti-clockwise rotation defined by  $R(z_1,z_2)=(-z_2,z_1)$ and let
$m=R(\na u)$, we will show there exists a set $\Gamma\subset S^1$ with $H^1(S^1\backslash \Gamma)\leq 40\pi\beta^{\frac{1}{8}}$ and $-\Gamma=\Gamma$ such that for any $\theta\in\Gamma$ we can find function $w_{\theta}:\Omega\rightarrow \R$ with the property
\begin{equation}
\label{eq3}
\int_{\Omega} \lt|\na w_{\theta}-R\lt(\Lambda_{\theta}\lt(m\rt)\rt)\rt|\leq c\beta^{\frac{1}{8}}.
\end{equation}
\end{a1}

\em Proof of Lemma \ref{L1}. \rm Let $M=2\lt[\frac{\beta^{-\frac{1}{4}}}{8}\rt]$, we divide $S^1$ into
$M$ disjoint connected subsets of length $\frac{2\pi}{M}$, denote them $A_1,A_2,\dots A_M$. We assume
they have been ordered sequentially, i.e.\ $\overline{A_i}\cap \overline{A_{i+1}}\not=\emptyset$ for $i=1,2,\dots
M-1$. Also assume they have been ordered so that $-A_i=A_{i+\frac{M}{2}}$ for $i=1,2,\dots \frac{M}{2}$. 
Let 
$$
\BI=\lt\{k\in\lt\{1,2,\dots \frac{M}{2}\rt\}:\lt|\lt\{x\in\Omega:\frac{\na u\lt(x\rt)}{\lt|\na u\lt(x\rt)\rt|}\in
\overline{A_k}\cup \overline{A_{k+\frac{M}{2}}}\rt\}\rt|\geq \beta^{\frac{1}{8}}\rt\}.
$$
Since $\ca{\BI}\beta^{\frac{1}{8}}\leq \lt|\Omega\rt|\leq 4\pi$ we have that
$\ca{\BI}\leq 4\pi\beta^{-\frac{1}{8}}$.

Let $\DI:=\lt\{k\in\lt\{2,3,\dots \frac{M}{2}-1\rt\}:\lt\{k-1,k,k+1\rt\}\cap \BI\not=\emptyset\rt\}$. A simple covering argument
shows that $\ca{\DI}\leq 20\pi\beta^{-\frac{1}{8}}$.

Let 
$\Gamma=\lt\{\theta\in S^1: \theta\in\bigcup_{k\in \lt\{2,3,\dots \frac{M}{2}-1\rt\}\backslash D} \overline{A_k}\cup 
\overline{A_{k+\frac{M}{2}}}\rt\}$. Note that for any $\theta\in \Gamma$ we have
\begin{equation}
\label{eqb10}
\lt|\lt\{x\in\Omega:\frac{\na u\lt(x\rt)}{\lt|\na u\lt(x\rt)\rt|}\in
B_{2\beta^{\frac{1}{4}}}\lt(\theta\rt)\cup B_{2\beta^{\frac{1}{4}}}\lt(-\theta\rt)\rt\}\rt|\leq 3\beta^{\frac{1}{8}}.
\end{equation}
So pick $\theta\in\Gamma$ without loss of generality we can assume $\theta=e_1$.
Let $s:\R\rightarrow \R_{+}$ be a smooth monotone function where $s(x)=0$ if $x\leq 0$ and
$s(x)=x$ if $x>\beta^{\frac{1}{4}}$ and $\|\na^2 s\|_{L^{\infty}}\leq \beta^{-\frac{1}{4}}$ and
$\|\na^3 s\|_{L^{\infty}}\leq \beta^{-\frac{1}{2}}$, it is clear such a function exists.

Let $\vp(z)=s(z\cdot e_1)=s(z_1)$.
Define $\Phi:\R^2\rightarrow \R^2$ by
\begin{eqnarray}
\label{eq608}
\Phi(z)&:=&\vp(z)\lt(\begin{matrix}  z_1\\ z_2\end{matrix}\rt)+\lt(\na\vp(z)\cdot \lt(\begin{matrix}  -z_2\\ z_1\end{matrix}\rt) \rt)
\lt(\begin{matrix}  -z_2\\ z_1\end{matrix}\rt)\nn\\
&=& \lt(\begin{matrix}  \vp\lt(z\rt) z_1+z_2^2 \vp_{,1}\lt(z\rt)\\  \vp\lt(z\rt) z_2-z_2 z_1 \vp_{,1}\lt(z\rt)  \end{matrix}\rt).
\end{eqnarray}
Define
\begin{equation}
\label{eqqz1}
\Psi\lt(z\rt)=\lt(\begin{matrix}\Psi_1\lt(z\rt)\\ \Psi_2\lt(z\rt)\end{matrix}\rt):=\lt(\begin{matrix} -\vp_{,1}\lt(z\rt)\\
\frac{z_2}{2} \vp_{,11}\lt(z\rt)\end{matrix}\rt).
\end{equation}
Recall $m\lt(z\rt):=R\lt(\na u\lt(z\rt)\rt)$
so $m$ is divergence free. Note (using the fact $\vp_{,2}\equiv 0$ and  $\vp_{,12}\equiv 0$ and
$\mathrm{div}m\equiv 0$ for the third inequality, and using $\mathrm{div}m=0$ for the last inequality)
\begin{eqnarray}
\label{eq7}
\mathrm{div}\lt[\Phi\lt(m\rt)\rt]&=&\mathrm{div}\lt(\begin{matrix}  \vp\lt(m\rt) m_1+m_2^2 \vp_{,1}\lt(m\rt)\\  \vp\lt(m\rt) m_2-m_2 m_1 \vp_{,1}\lt(m\rt)  \end{matrix}\rt)\nn\\
&=&(\vp_{,1}(m)m_{1,1}+\vp_{,2}(m)m_{2,1})m_1+\vp(m)m_{1,1}+2m_2 m_{2,1} \vp_{,1}(m)\nn\\
&~&\qd
+m_2^2 (\vp_{,11}(m) m_{1,1}+\vp_{,12}(m)m_{2,1})
+(\vp_{,1}(m)m_{1,2}+\vp_{,2}(m)m_{2,2})m_2\nn\\
&~&\qd+\vp(m)m_{2,2}
-((m_{1,2}m_2+m_1 m_{2,2})\vp_{,1}(m)\nn\\
&~&\qd+m_1m_2(\vp_{,11}(m)m_{1,2}+\vp_{,12}(m)m_{2,2}))\nn\\
&=&m_1\vp_{,1}(m)m_{1,1}+2 m_2 m_{2,1} \vp_{,1}(m)+m_2^2  m_{1,1} \vp_{,11}(m)+m_2 m_{1,2}\vp_{,1}(m)\nn\\
&~&\qd
-((m_{1,2} m_2+m_1 m_{2,2}) \vp_{,1}(m)+m_1 m_2 m_{1,2}\vp_{,11}(m)\nn\\
&=&2\vp_{,1}(m)(m_1 m_{1,1}+m_2 m_{2,1})-\vp_{,11}(m)m_2(m_1 m_{1,2}+m_2 m_{2,2}).
\end{eqnarray}
Note also that
\begin{eqnarray}
\Psi(m)\cdot \na(1-\lt|m\rt|^2)&=&-\Psi(m)\cdot \lt(\begin{matrix} 2(m_1 m_{1,1}+m_2 m_{2,1})\\
 2(m_1 m_{1,2}+m_2 m_{2,2})\end{matrix}\rt)\nn\\
&=&2\vp_{,1}(m)(m_1 m_{1,1}+m_2 m_{2,1})-m_2 \vp_{,11}(m)(m_1 m_{1,2}+m_2 m_{2,2})\nn
\end{eqnarray}
so by (\ref{eq7}) we have
\begin{equation}
\label{eq9}
\mathrm{div}\lt[\Phi\lt(m\rt)\rt]=\Psi(m)\cdot \na(1-\lt|m\rt|^2).
\end{equation}
Let $\ti{\Phi}:=R\lt(\Phi\rt)$ and $\ti{\Psi}:=R\lt(\Psi\rt)$ note $\mathrm{curl}\lt[\ti{\Phi}(m)\rt]\overset{(\ref{eq9})}{=}\mathrm{div}\lt[\Phi(m)\rt]=\Psi(m)\cdot \na(1-\lt|m\rt|^2)$. So
\begin{eqnarray}
\label{eq9.8}
\mathrm{curl}\lt[\ti{\Psi}(m)(1-\lt|m\rt|^2)\rt]&=&\mathrm{div}[\Psi(m)](1-\lt|m\rt|^2)+\Psi(m)\cdot \na(1-\lt|m\rt|^2)\nn\\
&=&\mathrm{div}\lt[\Psi\lt(m\rt)\rt](1-\lt|m\rt|^2)+\mathrm{curl}\lt[\ti{\Phi}\lt(m\rt)\rt].
\end{eqnarray}
Thus using the fact that $\lt|\na \Psi(z)\rt|\leq c\lt|z\rt|\|\na^3 \vp\|_{L^{\infty}\lt(\R^2\rt)}\leq 
c\beta^{-\frac{1}{2}}\lt|z\rt|$ we have 
\begin{eqnarray}
\label{eqe60}
&~&\mathrm{curl}\lt[\ti{\Phi}(m)-\ti{\Psi}(m)(1-\lt|m\rt|^2)\rt]\nn\\
&~&\qd\qd\qd\qd\overset{(\ref{eq9.8})}{=}-\mathrm{div}[\Psi(m)](1-\lt|m\rt|^2)\nn\\
&~&\qd\qd\qd\qd=-(\Psi_{1,1}(m)m_{1,1}+\Psi_{1,2}(m)m_{2,1}+\Psi_{2,1}(m)m_{1,2}+\Psi_{2,2}(m)m_{2,2})(1-\lt|m\rt|^2)\nn\\
&~&\qd\qd\qd\qd \leq c\beta^{-\frac{1}{2}}\lt|m\rt|\lt|1-\lt|m\rt|^2\rt|\lt|\na m\rt|.
\end{eqnarray}
Hence 
\begin{eqnarray}
\label{eq609}
\int_{\Omega} \lt|\mathrm{curl}\lt[\ti{\Phi}(m)-\ti{\Psi}(m)(1-\lt|m\rt|^2)\rt]\rt|&\leq&
c \beta^{-\frac{1}{2}} \int_{\Omega} \lt|m\rt|\lt|1-\lt|m\rt|^2\rt|\lt|\na m\rt|\nn
\end{eqnarray}

So using (\ref{eqe60}), note that if $x$ is such that $\lt|m(x)\rt|\geq 2$ then 
for $J(x):=\lt|m(x)\rt|^3$ we have $\lt|\na J(x)\rt|\leq c\lt|1-\lt|m\rt|^2\rt|\lt|\na m\rt|$ 
and so 
$$
\int_{\lt\{x:2\leq \lt|m(x)\rt|\leq 4\rt\}} \lt|\na J(x)\rt| dx\leq
c\int_{\Omega}  \lt|1-\lt|m\rt|^2\rt|\lt|\na m\rt|\leq c\beta
$$
so applying the Co-area formula we know $\int_{8}^{64} H^1(J^{-1}(s)) ds\leq c\beta$ thus we must be able to find $t\in [8,64]$ such that $H^1(J^{-1}(t))\leq c\beta$. Let 
\begin{equation}
\label{ezz1}
\GI:=\lt\{x\in \Omega: J(x)<t\rt\}
\end{equation}
so define $w:\Omega\rightarrow \R$ by 
\begin{equation}
\label{uleq2}
w(x)=\lt\{\begin{array}{ll} \ti{\Phi}(m)-\ti{\Psi}(m)(1-\lt|m\rt|^2)
&\text{ for } x\in \GI\\
0 &\text{ for }  x\in \Omega\backslash \GI \end{array} \rt.
\end{equation}
So if $x\in\GI$, 
\begin{eqnarray}
\label{eqzz11}
\mathrm{curl}(w)&=&
\mathrm{curl}\lt(\ti{\Phi}(m)-\ti{\Psi}(m)(1-\lt|m\rt|^2)\rt)\nn\\
&\overset{(\ref{eqe60}),(\ref{ezz1})}{\leq}& 
c\beta^{-\frac{1}{2}}\lt|1-\lt|m\rt|\rt|\lt|\na m\rt|.
\end{eqnarray}
So if $x\in \mathrm{int}\lt(\Omega\backslash \GI\rt)$, $\mathrm{curl}\lt(\ti{\Phi}(m)-\ti{\Psi}(m)(1-\lt|m\rt|^2)\rt)=0$. 

Since $m\in W^{1,1}(\Omega)$ and $\ti{\Phi}(x)-\ti{\Psi}(x)(1-\lt|x\rt|^2)$ is $C^1$ so the vector field 
$\ti{\Phi}(m)-\ti{\Psi}(m)(1-\lt|m\rt|^2)$ is BV by Theorem 3.94 \cite{amb}. So 
by Theorem 3.83 \cite{amb} we have that $w$ is also BV and the singular part of $\na w$, 
which we denote by $\lt[\na w\rt]_s$, is 
supported on $J^{-1}(t)\cap \Omega$ and as $\lt|\ti{\Phi}(m(x))\rt|\leq c\lt|m(x)\rt|^2$ and $\lt|\ti{\Psi}(m(x))\rt|\leq c\beta^{-\frac{1}{4}}\lt|m(x)\rt|$ we have that 
$$
\text{ess sup}_{J^{-1}(t)\cap \Omega}  
\lt|\ti{\Phi}(m(x))-\ti{\Psi}(m(x))(1-\lt|m(x)\rt|^2)\rt|\leq c\beta^{-\frac{1}{4}}
$$
and thus $\|\lt[\na w\rt]_s\|(S)\leq c\beta^{-\frac{1}{4}}H^1(J^{-1}(t)\cap \Omega)\leq c\beta^{\frac{3}{4}}$. 
Now we know that for any set $S\subset \Omega$, 
$$
\|\mathrm{curl} w\|(S)\leq c\|\na w\|(S)
$$ 
and so in particular 
\begin{equation}
\label{uleq1}
\|\mathrm{curl} w\|(J^{-1}(t))\leq c\|\na w\|(J^{-1}(t))\leq c\beta^{\frac{3}{4}}. 
\end{equation}
Thus 
\begin{eqnarray}
\label{uleq3}
\|\mathrm{curl} w\|(\Omega)&\leq& \|\mathrm{curl} w\|(J^{-1}(t))+
\|\mathrm{curl} w\|(\GI)\nn\\
&~&+\|\mathrm{curl} w\|(\mathrm{int}(\Omega\backslash \GI))\nn\\
&\overset{(\ref{eqzz11}),(\ref{uleq1})}{\leq}&c\beta^{\frac{3}{4}}+c\beta^{-\frac{1}{2}} \int_{\GI}\lt|1-\lt|m\rt|\rt|\lt|\na m\rt|\nn\\
&\overset{(\ref{eq1})}{\leq}&c\sqrt{\beta}. 
\end{eqnarray}

Now we try and understand the nature of vector field $\ti{\Phi}(m(x))-\ti{\Psi}(m(x))(1-\lt|m(x)\rt|^2)$. 
Note that if $z\in N_{\sqrt{\beta}}\lt(S^1\rt)\cap\lt\{z_1>0\rt\}\backslash 
\lt(B_{2\beta^{\frac{1}{4}}}\lt(e_2\rt)\cup B_{2\beta^{\frac{1}{4}}}\lt(-e_2\rt)\rt)$ then $\vp(z)=z_1$, $\vp_{,1}(z)=1$ and so $\Phi(z)\overset{(\ref{eq608})}{=}\lt(\begin{matrix} z_1^2+z_2^2\\0\end{matrix}\rt)$
on the other hand if $z\in N_{\sqrt{\beta}}(S^1)\cap\lt\{z_1\leq 0\rt\}\backslash \lt(B_{2\beta^{\frac{1}{4}}}\lt(e_2\rt)\cup B_{2\beta^{\frac{1}{4}}}\lt(-e_2\rt)\rt)$
then $\vp(z)=\vp_{,1}(z)=0$ and so
$\Phi(z)=\lt(\begin{matrix} 0\\0\end{matrix}\rt)$. 

Now if $z\in N_{\sqrt{\beta}}(S^1)\cap \lt\{z_1> 0\rt\}\backslash \lt(B_{2\beta^{\frac{1}{4}}}\lt(e_2\rt)\cup B_{2\beta^{\frac{1}{4}}}\lt(-e_2\rt)\rt)$ we have
\begin{eqnarray}
\label{eqb2}
\lt|(\ti{\Phi}(z)-\ti{\Psi}(z)(1-\lt|z\rt|^2))-R\lt(\Lambda_{e_1}(z)\rt)\rt|&\leq&
\lt|\ti{\Phi}(z)-R\lt(\Lambda_{e_1}(z)\rt)\rt|+c\sqrt{\beta} \sup_{z\in N_{\sqrt{\beta}}(S^1)}\lt|\ti{\Psi}(z)\rt|\nn\\
&=&\lt|R\lt(\begin{matrix} z_1^2+z_2^2\\ 0\end{matrix}\rt)-R\lt(\begin{matrix} 1 \\ 0\end{matrix}\rt)\rt|+c\beta^{\frac{1}{4}} \nn\\
&\leq&c\beta^{\frac{1}{4}}.
\end{eqnarray}
And if we have $z\in N_{\sqrt{\beta}}\lt(S^1\rt)\cap\lt\{z_1\leq 0\rt\}\backslash 
\lt(B_{2\beta^{\frac{1}{4}}}\lt(e_2\rt)\cup B_{2\beta^{\frac{1}{4}}}\lt(-e_2\rt)\rt)$
arguing in the same way we can conclude
\begin{equation}
\label{eq12.6}
\lt|(\ti{\Phi}(z)-\ti{\Psi}(z)(1-\lt|z\rt|^2))-R\lt(\Lambda_{e_1}(z)\rt)\rt|\leq c\beta^{\frac{1}{4}}.
\end{equation}
Let $\Pi:=\lt\{z\in\Omega:\lt|m(z)\rt|\in (1-\sqrt{\beta},1+\sqrt{\beta})\rt\}$ and let
\begin{equation}
\label{eqc33}
\EI:=\lt\{x\in\Omega:\frac{\na u\lt(x\rt)}{\lt|\na u\lt(x\rt)\rt|}\in 
B_{2\beta^{\frac{1}{4}}}\lt(e_1\rt)\cup B_{2\beta^{\frac{1}{4}}}\lt(-e_1\rt)\rt\},
\end{equation}
note from (\ref{eqb10}) we know $\lt|\EI\rt|\leq 3\beta^{\frac{1}{8}}$. Note also 
$\sqrt{\beta}\lt|\Omega\backslash \Pi\rt|\leq 
c\int_{\Omega\backslash \Pi} \lt|1-\lt|\na u\rt|^2\rt|\overset{(\ref{eq2})}{\leq} \beta$ thus 
\begin{equation}
\label{ulazu1}
\lt|\Omega\backslash \Pi\rt|\leq c\sqrt{\beta}.
\end{equation}
Now from (\ref{eqb2}) and (\ref{eq12.6})
\begin{equation}
\label{eq15.7}
\lt|\int_{\Pi\backslash \EI} (\ti{\Phi}(m)-\ti{\Psi}(m)(1-\lt|m\rt|^2))-R\lt(\Lambda_{e_1}(m)\rt)  dz\rt|\leq c\beta^{\frac{1}{4}}
\end{equation}
on the other hand recalling the fact that $\lt|\ti{\Psi}(z)\rt|\leq \beta^{-\frac{1}{4}}\lt|z\rt|$, 
$\lt|\ti{\Phi}(z)\rt|\leq c\lt|z\rt|^2$ and using the definition of $\GI$ (see (\ref{ezz1})) we have 
\begin{eqnarray}
\label{eq173}
&~&\lt|\int_{\GI\backslash \Pi}   \lt((\ti{\Phi}(m)-\ti{\Psi}(m)(1-\lt|m\rt|^2))-R\lt(\Lambda_{e_1}(m)\rt)\rt) dz\rt|\nn\\
&~&\qd\qd\qd\leq c\lt|\GI\backslash \Pi\rt|\nn\\
&~&\qd\qd\qd \overset{(\ref{ulazu1})}{\leq} c\sqrt{\beta}.
\end{eqnarray}
Thus applying (\ref{eq15.7}) to (\ref{eq173}) gives 
\begin{equation}
\label{eq19}
\lt|\int_{\GI\backslash \EI} \lt((\ti{\Phi}(m)-\ti{\Psi}(m)(1-\lt|m\rt|^2))-
R\lt(\Lambda_{e_1}(m)\rt)\rt) dz \rt|\leq c\beta^{\frac{1}{4}}.
\end{equation}
Recall we have $\lt|\EI\rt|\leq 3\beta^{\frac{1}{8}}$ so 
\begin{eqnarray}
\label{eqc56}
\lt|\int_{\EI\cap \GI} \lt((\ti{\Phi}(m)-\ti{\Psi}(m)(1-\lt|m\rt|^2))-R\lt(\Lambda_{e_1}(m)\rt)\rt) dz \rt|
&\leq&c\lt|\EI\rt|\nn\\
&\leq&c\beta^{\frac{1}{8}}.\nn
\end{eqnarray}
Putting this inequality together with (\ref{eq19}) gives
\begin{equation}
\label{eqzz10}
\lt|\int_{\GI} \lt((\ti{\Phi}(m)-\ti{\Psi}(m)(1-\lt|m\rt|^2))-R\lt(\Lambda_{e_1}(m)\rt)\rt) dz\rt|\leq c\beta^{\frac{1}{8}}.
\end{equation}
So by definition of $w$ (see (\ref{uleq2})) we have that 
\begin{eqnarray}
\label{ulaz5}
\lt|\int_{\Omega} w-R(\Lambda_{e_1}(m)) dz\rt|&\overset{(\ref{eqzz10})}{\leq}&c\beta^{\frac{1}{8}}+\lt|\int_{\Omega\backslash \GI} 
R(\Lambda_{e_1}(m)) dz\rt|\nn\\
&\leq&c\beta^{\frac{1}{8}}+\lt|\Omega\backslash \GI\rt|\nn\\
&\overset{(\ref{ulazu1})}{\leq}&c\beta^{\frac{1}{8}}.
\end{eqnarray}

Now from (\ref{uleq3}) applying Theorem 4.3 from (\cite{amb2}) there exists
$w_{e_1}\in W^{1,1}\lt(\Omega\rt)$ such that
\begin{equation}
\label{eq20}
\int_{\Omega} \lt|\na w_{e_1}-w\rt| dz\leq c\beta^{\frac{1}{8}}
\end{equation}
thus putting this together with (\ref{ulaz5}) and gives (\ref{eq3}). $\Box$

%
%
%
%

\begin{a1}
\label{L2}
Let $\Omega$ be a convex body centered on $0$ and let
$u:W^{2,2}(\Omega)\rightarrow \R$ be a function satisfying (\ref{eq1}) and (\ref{eq2}) and $u=0$ on 
$\partial \Omega$ and $\na u(z)\cdot\eta_z=1$ on $\partial \Omega$ in the sense of trace, where $\eta_z$ is the inward pointing
unit normal to $\partial \Omega$ at $z$.

For any $r>0$ define $\OO_r:=N_{r}(\Omega)$, we will show we can construct a function
$\ui:W^{2,1}(\OO_r)\rightarrow \R$ satisfying
\begin{equation}
\label{eq24}
\int_{\OO_r} \lt|1-\lt|\na \ui\rt|^2\rt|\lt|\na^2 \ui\rt| dz\leq \beta,\;\;\int_{\OO_r}
\lt|1-\lt|\na \ui\rt|^2\rt| dz\leq \beta,
\end{equation}
and
\begin{equation}
\label{eq27}
\ui\lt(z\rt)=
\begin{cases} u(z)+r & \text{for }z\in \overline{\Omega}\\
 r-d(z,\Omega) &\text{if } z\in \OO_r\backslash \Omega
\end{cases}
\end{equation}
\end{a1}
\em Proof of Lemma \ref{L2}. \rm

\em Step 1. \rm We will show $\na u(x)=\eta_x$ for $H^1$ a.e.\ $x\in\partial \Omega$

\em Proof of Step 1. \rm Recall $\na u\in W^{1,1}(\Omega)$ and $\na u$ is defined on $\partial \Omega$ 
in the sense of trace, as the trace operator is bounded we know 
$\int_{\partial \Omega} \lt|\na u\rt| dH^1 <\infty$. 

We define 
\begin{equation}
\label{eq27.6}
v\lt(z\rt)=
\begin{cases} u(z)& \text{for }z\in \overline{\Omega}\\
 0 &\text{if } z\in \Omega_r\backslash \Omega
\end{cases}
\end{equation}

So note the vector field $\na v(z)$ is equal to $\na u(z)$ inside $\Omega$ and is zero outside, so by Theorem 3.8 \cite{amb} $\na v\in BV(\Omega_r)$ and hence by Theorem 3.76 \cite{amb} and Theorem 2, 
Section 5.3 \cite{evans2} for $H^1$ a.e.\ $x\in \partial \Omega$ the following limits exist 
\begin{equation}
\label{ula40}
\lim_{\rho\rightarrow 0} \Xint{-}_{B_{\rho}(x)\cap \lt\{z:(z-x)\cdot \eta_x>0\rt\}} \lt|\na v(z)-\na u(x)\rt| dz=0
\end{equation}
and 
\begin{equation}
\label{ula41.7}
\lim_{\rho\rightarrow 0} \Xint{-}_{B_{\rho}(x)\cap \lt\{z:(z-x)\cdot \eta_x\leq 0\rt\}} \lt|\na v(z)\rt| dz=0.
\end{equation}

Let $w_x^{\rho}(z)=\frac{v(x+\rho z)}{\rho}$, by (\ref{ula40}) and (\ref{ula41.7}) for any sequence 
$\rho_n\rightarrow 0$ we have $w_x^{\rho_n}(z)\overset{W^{1,1}}{\rightarrow} w_x$ as $n\rightarrow \infty$ 
where 
\begin{equation}
\label{ula50}
w_x\lt(z\rt)=
\begin{cases} \na u(x)\cdot z& \text{ for }z\in H(0,\eta_x)\\
 0 &\text{ for } z\in H(0,-\eta_x)
\end{cases}
\end{equation}
however $\na w_x$ would not be curl free unless $\na u(x)=\lm\eta_x$ for some $\lm\in \R$. As we 
know $\na u(x)\cdot \eta_x=1$ this implies $\na u(x)=\eta_x$ for $H^1$ a.e.\ $x\in \partial \Omega$. 
This completes the proof of Step 1. \nl

%
%
%
%
%
%

\em Step 2. \rm  For any $z\in \OO_r\backslash\Omega$, $\ui(z)=d(z,\partial \OO_r)$.

\em Proof of Step 2. \rm Note that $\|\na \ui\|_{L^{\infty}(\Omega_r\backslash \Omega)}\leq 1$.
Let $x\in\partial\OO_r$, let $q(x)$ be the metric projection onto a convex set $\Omega$, i.e.\ the unique
point for which $\lt|x-q(x)\rt|=d(x,\Omega)$. Since $x\in\partial\OO_r=\partial(N_r(\Omega))=
\lt\{x\in \Omega^c:d(x,\Omega)=r\rt\}$ so $\lt|x-q(x)\rt|=r$.

Since $\ui(x)=0$ and $\ui(q(x))=r$ and as $\ui$ is $1$-Lipschitz on $\OO_r\backslash \Omega$
this implies $\ui((1-\alpha)x+\alpha q(x))=\alpha r$ for any $\alpha\in [0,1]$.

Now let $Q(z):=d(z,\partial \OO_r)$. For every $x\in\partial \OO_r$, $Q(q(x))\leq \lt|q(x)-x\rt|=r$.
As $\partial \OO_r=\partial( N_{r}(\Omega))$ so we know $Q(q(x))\geq r$ and thus have $Q(q(x))=r$. We also know
$Q$ is $1$-Lipschitz and $Q(x)=0$, thus in the same way as before
$Q((1-\alpha)x+\alpha q(x))=\alpha r$ for any $\alpha\in \lt[0,1\rt]$. Therefor
$Q(z)=\ui(z)$ for any $z\in[x,q(x)]$, $x\in\partial\OO_r$ and this completes the proof of Step 2. \nl

\em Step 3. \rm We will show that $\ui\in W^{2,1}(\OO_r)$ and that $\ui$ satisfies (\ref{eq24}).

\em Proof of Step 3. \rm First we claim that $\ui\in W^{2,1}(\Omega_r\backslash \Omega)$ and 
\begin{equation}
\label{eq18.5}
\int_{\OO_r\backslash \Omega} \lt|\na^2 \ui\rt| dz\leq c.
\end{equation}
%

Note that $\ui(z)=\mathrm{dist}(z,\partial \Omega_r)$ in $\Omega_r\backslash \Omega$. By Corollary 1.4 \cite{amdel} 
for any compact subset $\Omega'\subset\subset \Omega_r$ we have $\na \ui\in SBV(\Omega'\backslash \Omega)$. Also 
as $\ui(z)=r-\mathrm{dist}(z,\Omega)$ for any $z\in \Omega_r\backslash \Omega$ again by  Corollary 1.4 \cite{amdel} 
for any compact subset $\Omega''\subset\subset \R^2\backslash \overline{\Omega}$ we have $\na \ui\in SBV\lt((\Omega_r\backslash \Omega)\cap 
\Omega''\rt)$. Putting these thing together we have $\na \ui\in SBV(\Omega_r\backslash \Omega)$. 
Recall $\ui(x)=r-d(z,\Omega)$ for $z\in \Omega_r\backslash \Omega$, so 
as $\Omega$ is convex for every $z\in \Omega_r\backslash \Omega$ there is a unique point $b(z)\in 
\partial\Omega$ such that $d(z,\Omega)=\lt|b(z)-z\rt|$ and 
$\na \ui(z)=\frac{b(z)-z}{\lt|b(z)-z\rt|}$, since $b$ is a continuous function this shows that 
$\na\ui$ is continuous on $\Omega_r\backslash\overline{\Omega}$, hence 
$S_{\na \ui}\cap \OO_r\backslash \overline{\Omega}=\emptyset$ (recall Definition 3.63 \cite{amb}). 
So by equation (4.2) of Section 4.1 \cite{amb} we 
have that $\na\ui\in W^{1,1}(\Omega_r\backslash \Omega)$. So in particular (\ref{eq18.5}) holds true.

Since $\Omega$ is an extension domain by Theorem 1, Section 4.4 \cite{evans2} there exists
a function $p:W^{1,2}(\R^2)\rightarrow \R^2$ such that $p(z)=\na\ui(z)$ on $\Omega$ and $\spt p$ is
compact. Similarly as $\OO_r\backslash \Omega$ is an extension domain there exists a function
$q:W^{1,1}(\R^2)\rightarrow \R^2$ such that $q(z)=\na \ui(z)$ on $\OO_r\backslash \Omega$ and
$\spt q$ is compact. We define $w:\OO_r\rightarrow \R^2$ by $w:=p\cha_{\Omega}+q\cha_{\OO_r\backslash \Omega}$,
by Theorem 3.83 \cite{amb} $w\in BV(\OO_r:\R^2)$ and since $p$ and $q$ agree on $\partial \Omega$
we have that $\na w$ as a measure is absolutely continuous with respect to Lebesgue measure (and hence $w\in W^{1,1}(\OO_r:\R^2)$) and
$\na w=\na p\cha_{\Omega}+\na q\cha_{\OO_r\backslash \Omega}$. Now as $w=\na \ui$ a.e.\ in $\OO_r$ we have that
$\na \ui\in W^{1,1}(\OO_r)$.

Since $\na^2 \ui\in L^1$ we know
\begin{eqnarray}
\int_{\OO_r} \lt|1-\lt|\na \ui\rt|^2\rt|\lt|\na^2 \ui\rt| dz&=&
\int_{\Omega} \lt|1-\lt|\na \ui\rt|^2\rt|\lt|\na^2 \ui\rt| dz
+\int_{\OO_r\backslash\Omega} \lt|1-\lt|\na \ui\rt|^2\rt|\lt|\na^2 \ui\rt| dz\nn\\
&=&\int_{\Omega} \lt|1-\lt|\na \ui\rt|^2\rt|\lt|\na^2 \ui\rt| dz\nn\\
&\leq&\beta.\nn
\end{eqnarray}
Similarly $\int_{\OO_r} \lt|1-\lt|\na \ui\rt|^2\rt| dz=\int_{\Omega} \lt|1-\lt|\na \ui\rt|^2\rt| dz\leq \beta$. $\Box$

%
%
%
%

\begin{a1}
\label{L3}
Let $\Omega$ be a convex body with $\mathrm{diam}(\Omega)=2$. Let
$u:W^{2,2}(\Omega)\rightarrow \R$ be a function satisfying (\ref{eq1}) and (\ref{eq2}) and $u=0$ on 
$\partial \Omega$ and $\na u(z)\cdot \eta_z=1$ on $\partial \Omega$ in the sense of trace where $\eta_z$ is the
inward pointing unit normal to $\partial\Omega$ at $z$. For any $x,v\in \R^2$ let 
$H(x,v):=\lt\{z\in \R^2:(z-x)\cdot v>0\rt\}$.

Let $\Gamma\subset S^1$ be the set constructed in Lemma \ref{L1}.
Let $\UI:=\OO_{1/10}$ be the convex body and $\ui:W^{2,1}(\UI)\rightarrow \R$ be the
function constructed in Lemma \ref{L2}. Let $R$ be the anti-clockwise rotation defined by 
$R(z_1,z_2)=(-z_2,z_1)$. Let $R_0\in\lt\{R^{-1},R\rt\}$. There exists a set $\wt{\Gamma}\subset \Gamma$ with 
$H^1(\Gamma\backslash \wt{\Gamma})=0$ such that for every $\theta\in \wt{\Gamma}$ there exists unique points $a_{\theta},b_{\theta}\in \partial \UI$ with
$\eta_{a_{\theta}}=\theta$ and $\eta_{b_{\theta}}=-\theta$ with the property that if we define
$\GI^{R_0}_{\theta}:=\lt\{z\in\UI:\na \ui(z)\cdot R_0^{-1}\theta> 0\rt\}$,  
\begin{equation}
\label{eq104}
\lt|\UI\cap H\lt(\frac{a_{\theta}+b_{\theta}}{2}, R_0\lt(\frac{b_{\theta}-a_{\theta}}
{\lt|b_{\theta}-a_{\theta}\rt|}\rt)\rt)\backslash \GI^{R_0}_{\theta}\rt|\leq c\beta^{\frac{1}{24}}.
\end{equation}
\end{a1}

\em Proof of Lemma \ref{L3}. \rm Without loss of generality assume $\Omega$ is centered on $0$, i.e.\ 
$\int_{\Omega} z dz=0$. Since $\partial \UI$ is smooth and $\UI$ is convex there exists a set $\Xi\subset S^1$ with 
$H^1(S^1\backslash \Xi)=0$ with the following property,
\begin{equation}
\label{ulzq1}
\exists\;\text{unique}\; a_{\thea}\in\partial \UI\text{ with }\eta_{a_{\thea}}=\thea\text{ and a unique }b_{\thea}\in\partial \UI\text{ with }\eta_{b_{\thea}}=-\thea \text{ for all }\thea\in \Xi.
\end{equation}

Now by Lemma \ref{L2}, (\ref{eq24}) function $\ui$ satisfies (\ref{eq1}) and (\ref{eq2}) so by 
Lemma \ref{L1} there 
exists $\Gamma\subset S^1$ with $H^1(S^1\backslash \Gamma)\leq 40\pi\beta^{\frac{1}{8}}$ satisfying 
(\ref{eq3}) for every $\theta\in \Gamma$. Define $\wt{\Gamma}:=\Gamma\cap \Xi$. Pick $\theta\in \wt{\Gamma}$ and 
let $\thea:=R R_0^{-1} \theta$ so note that $\thea=\theta$ or $\thea=-\theta$ depending
on whether $R_0=R$ or $R_0=R^{-1}$.

Note since $\Omega$ is convex $\Omega\subset \overline{H(a_{\thea},\thea)}$ we also know that 
$b_{\thea}\in H(a_{\thea},\thea)$  (since otherwise given that $\partial \Omega$ is smooth it would 
not be possible that $\eta_{b_{\thea}}=-\thea$), hence defining $\tau_{\thea}=\frac{b_{\thea}-a_{\thea}}{\lt|b_{\thea}-a_{\thea}\rt|}$ we have $\tau_{\thea}\cdot \thea>0$.

\begin{figure}[h]
\centerline{\psfig{figure=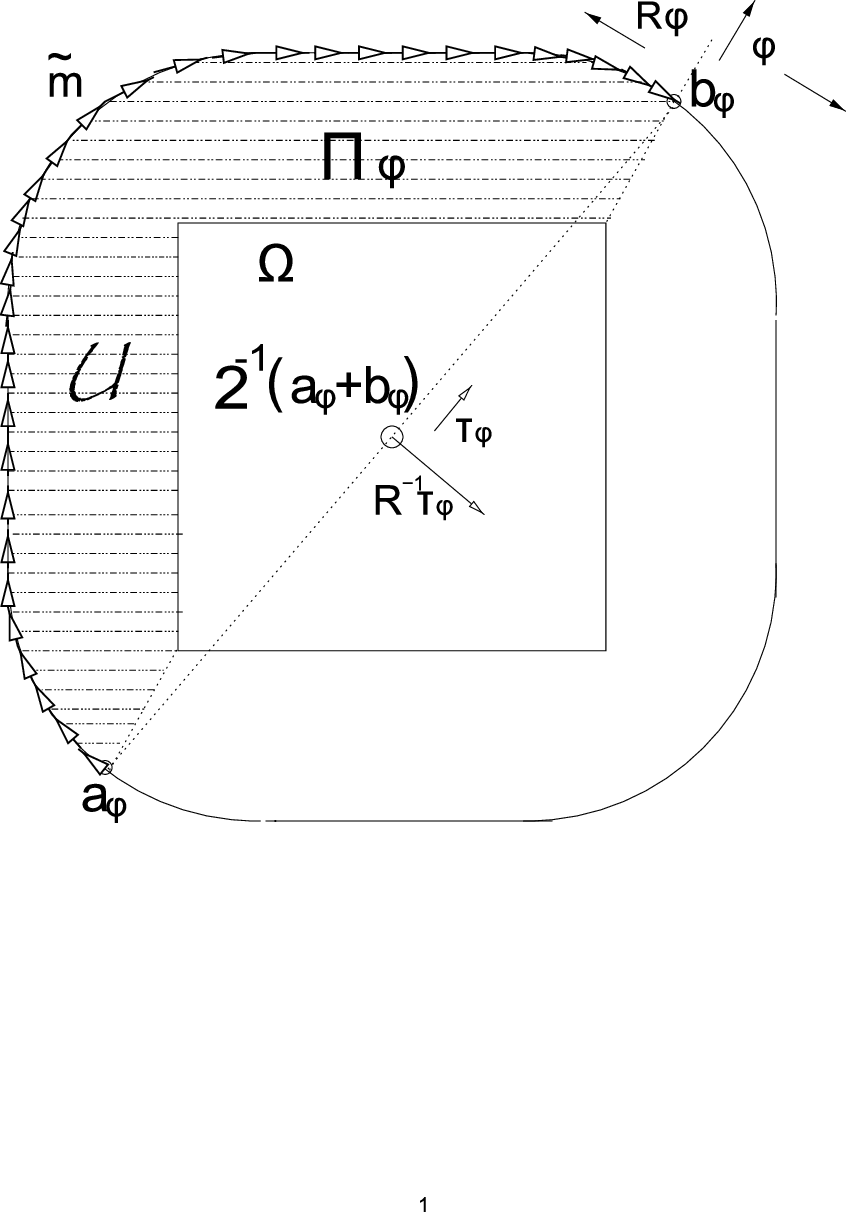,height=15cm,angle=0}}
\caption{}\label{fig5}
\end{figure}

Let
$\mi=R(\na \ui)$, it is easy to see that
\begin{equation}
\label{eq300}
\Pi_{\thea}:=\lt\{z\in\UI\backslash \Omega:\mi(z)\cdot \thea> 0\rt\}
=\lt\{z\in\UI\backslash \Omega:\na u\lt(z\rt)\cdot R^{-1}\thea>0\rt\}
\end{equation}
forms a connected set whose boundary is contained in $\partial\UI$ and $\partial \Omega$ and in two lines parallel
to $\thea$, see figure \ref{fig5}, also note the endpoints of
$\partial \UI\cap \overline{\Pi_{\thea}}$ are given by $a_{\thea}$ and $b_{\thea}$.

Since either $\thea=\theta\in \wt{\Gamma}$ or $\thea=-\theta\in \wt{\Gamma}$ so we can apply Lemma \ref{L1}, to $\mi$ and thus there exists function $w_{\thea}:\UI\rightarrow \R$ such that
\begin{equation}
\label{eqe34}
\int_{\UI} \lt|\na w_{\thea}-R\lt(\Lambda_{\thea}\lt(\ti{m}\rt)\rt)\rt| dx\leq
c\beta^{\frac{1}{8}}.
\end{equation}
By the Co-area formula and Chebyshev's inequality there exists
a set $H\subset \lt[0,1/10\rt]$ such that $H^1( \lt[0,1/10\rt]\backslash H)\leq
c\beta^{\frac{1}{24}}$ where
\begin{equation}
\label{eq303}
\int_{\ui^{-1}\lt(t\rt)} \lt|\na w_{\thea}-R\lt(\Lambda_{\thea}\lt(\mm\rt)\rt)\rt| dH^1 \leq
c\beta^{\frac{1}{12}}\text{ for all }r\in H.
\end{equation}
Pick $s_0\in\lt[1/10-c\beta^{\frac{1}{24}},1/10\rt]\cap H$.
Recall $\tau_{\thea}=\frac{b_{\thea}-a_{\thea}}{\lt|b_{\thea}-a_{\thea}\rt|}$ and define 
\begin{equation}
\label{eq751}
\WI_{\thea}:=\overline{\UI}\cap H\lt(\frac{a_{\thea}+b_{\thea}}{2},R\tau_{\thea}\rt).
\end{equation}
We claim that
\begin{equation}
\label{eq301.5}
\partial \UI\cap \overline{\Pi_{\thea}}= \partial\UI\cap \overline{\WI_{\thea}}.
\end{equation}
Since the endpoints of $\partial \UI\cap \overline{\Pi_{\thea}}$ are the
same as the endpoints of $\partial\UI\cap \overline{\WI_{\thea}}$ it is sufficient to show
$H^1\lt(\partial\UI\cap\overline{\Pi_{\thea}}\cap\overline{\WI_{\thea}}\rt)>0$. Let
$$
\Lambda=\sup\lt\{\lambda>0:\lt(\frac{a_{\thea}+b_{\thea}}{2}
+\lm R\tau_{\thea}+\la\tau_{\thea}\ra\rt)\cap\partial\UI\not=\emptyset\rt\}
$$
then let $c_{\thea}$ be the point given by
$\lt(\frac{a_{\thea}+b_{\thea}}{2}+\Lambda R\tau_{\thea}+\la\tau_{\thea}\ra\rt)\cap\partial\UI$,
since $\partial \UI$ is smooth $\eta_{c_{\thea}}=R^{-1}\tau_{\thea}$, so
$\na u(c_{\thea})=R^{-1}\tau_{\thea}$ and thus
$\na u\lt(c_{\thea}\rt)\cdot R^{-1}\thea=R^{-1}\tau_{\thea}\cdot R^{-1}\thea=\tau_{\thea}\cdot \thea>0$. As 
this inequality is strict, in a neighborhood of $c_{\thea}$ the same
inequality will be satisfied. Thus we have $H^1\lt(\partial\UI\cap\overline{\Pi_{\thea}}\cap\overline{\WI_{\thea}}\rt)>0$
and so we have established (\ref{eq301.5}).

By the construction of $\Pi_{\thea}$, $\WI_{\thea}$ and by (\ref{eq301.5}) and 
the choice of $s_0\in \lt[\frac{1}{10}-c\beta^{\frac{1}{24}},\frac{1}{10}\rt]$ we have
\begin{equation}
\label{eq301}
H^1\lt(\partial\Omega_{s_0}\cap \overline{\Pi_{\thea}}\triangle \overline{\WI_{\thea}} \rt)\leq c\beta^{\frac{1}{24}}.
\end{equation}

There must exist
$\psi\in(0,2\beta^{\frac{1}{24}})$ such that defining $Q:=\lt(\begin{matrix}  \cos \psi & -\sin \psi\\ \sin\psi & \cos\psi\end{matrix}\rt)$ we have 
\begin{equation}
\label{uz8}
\lt|R\thea\cdot Q \tau_{\thea}\rt|>\beta^{\frac{1}{24}}.
\end{equation}
Let $\zeta_{\thea}:= \frac{a_{\thea}+b_{\thea}}{2}+\CI_2\beta^{\frac{1}{24}}R\tau_{\thea}$. From the
construction it is clear that we can chose constant $\CI_2$ large enough so that
$$
\ca{\partial \Omega_{s_0}\cap H\lt(\frac{a_{\thea}+b_{\thea}}{2},R\tau_{\thea}\rt)
\cap\lt\{\zeta_{\thea}+\la Q\tau_{\thea}\ra\rt\}}=2.
$$
Let 
\begin{equation}
\label{ghg1}
\mathfrak{A}:=\sup\lt\{t>0: \partial \Omega_{s_0}\cap \lt\{\zeta_{\thea}+tR\tau_{\thea}+\la Q\tau_{\thea}\ra \rt\}\not=\emptyset  \rt\}.
\end{equation}
For $t\in (0,\mathfrak{A})$ let 
$\varrho^1_t,\varrho^2_t$ be the points defined by 
$\lt\{\varrho^1_t,\varrho^2_t\rt\}=\partial \Omega_{s_0}\cap
\lt\{\zeta_{\thea}+t R\tau_{\thea}+\la Q\tau_{\thea}\ra\rt\}$ and
$\varrho^2_t\cdot Q\tau_{\thea}\geq \varrho^1_t\cdot Q\tau_{\thea}$. 
By (\ref{eq301}) we can assume constant $\CI_2$ was chosen large enough so that
$\varrho^1_t,\varrho^2_t\in\Pi_{\thea}$. Let $\Sigma_t$ be the connected component
of $\partial \Omega_{s_0}\backslash \lt\{\varrho^1_t,\varrho^2_t\rt\}$ that lies inside $\Pi_{\thea}$. 
Thus
\begin{eqnarray}
\label{eq305.6}
\lt|(w_{\thea}(\varrho^2_t)-w_{\thea}(\varrho^1_t))-(\varrho^2_t-\varrho^1_t)\cdot R\thea\rt|
&=&\lt|\int_{\Sigma_t} \na w_{\thea}(z)\cdot t_z dH^1 z- \int_{\Sigma_t} R\thea\cdot t_z dH^1 z   \rt|\nn\\
&=&\lt|\int_{\Sigma_t} (\na w_{\thea}(z)-R\thea)\cdot t_z dH^1 z   \rt|\nn\\
&\overset{(\ref{eq303})}{\leq}&c\beta^{\frac{1}{12}}.
\end{eqnarray}

Let
\begin{equation}
\label{eqe33}
e_t=\int_{\lt[\varrho^1_t,\varrho^2_t\rt]} \lt|\na w_{\thea}-
R\lt(\Lambda_{\thea}(\ti{m})\rt)\rt| dH^1 x,
\end{equation}
so by the fundamental theorem of Calculus 
$$
\lt|\lt(w_{\thea}(\varrho^2_t)-w_{\thea}(\varrho^1_t)\rt)-\int_{\lt[\varrho^1_t,\varrho^2_t\rt]}
R\lt(\Lambda_{\thea}(\ti{m})\rt)\cdot Q\tau_{\thea} dH^1 x\rt|\leq e_t. 
$$
Thus in combination with (\ref{eq305.6}) we have
\begin{equation}
\label{eq313}
\lt|\lt(\varrho^2_t-\varrho^1_t\rt)\cdot R\thea-\int_{\lt[\varrho^1_t,\varrho^2_t\rt]}
R\lt(\Lambda_{\thea}(\ti{m})\rt)\cdot Q\tau_{\thea} dH^1 x\rt|\leq e_t+c\beta^{\frac{1}{12}}.
\end{equation}

Given the definition of $\Lambda_{\thea}$ (see (\ref{eq750})) and of $\GI_{\theta}^{R_0}$ (see
the statement of Lemma \ref{L3}) so
$$
R(\Lambda_{\thea}(\mi(x)))=R\thea\Leftrightarrow \mi(x)\cdot \thea>0\Leftrightarrow
\na \ui(x)\cdot R^{-1} \thea>0\Leftrightarrow \na \ui(x)\cdot R_0^{-1} \theta>0\Leftrightarrow
x\in \GI_{\theta}^{R_0}.
$$
In exactly the same way $\Lambda_{\thea}(\mi(x))=0\Leftrightarrow x\not\in \GI_{\theta}^{R_0}$. Hence
$$
\int_{\lt[\varrho^1_t,\varrho^2_t\rt]} \Lambda_{\thea}(\mi(x)) dH^1 x=\thea
H^1\lt(\lt[\varrho^1_t,\varrho^2_t\rt]\cap \GI^{R_0}_{\theta}\rt)
$$
which from (\ref{eq313})
$$
\lt| \lt(\varrho^2_t-\varrho^1_t\rt)\cdot R\thea
-Q\tau_{\thea}\cdot R\thea H^1\lt(\lt[\varrho^1_t,\varrho^2_t\rt]\cap \GI^{R_0}_{\theta}\rt) \rt|\leq
e_t+c\beta^{\frac{1}{12}}
$$
since (recall (\ref{uz8})) we chose $Q$ so that $\lt|R\thea\cdot Q\tau_{\thea}\rt|>\beta^{\frac{1}{24}}$ and
since $\frac{\varrho^2_t-\varrho^1_t}{\lt| \varrho^2_t-\varrho^1_t\rt|}=Q\tau_{\thea}$
so
$$
\lt|\lt|\varrho^2_t-\varrho^1_t\rt|-  H^1\lt(\lt[\varrho^1_t,\varrho^2_t\rt]\cap \GI^{R_0}_{\theta}\rt)\rt|
\leq c\beta^{-\frac{1}{24}}e_{t}+c\beta^{\frac{1}{24}}.
$$
Thus (recall definition (\ref{ghg1}) of $\mathfrak{A}$)
\begin{equation}
\label{eq501}
H^1\lt(\lt[\varrho^1_t,\varrho^2_t\rt]\cap \GI^{R_0}_{\theta}\rt)\geq \lt|\varrho^1_t-\varrho^2_t\rt|
-c \beta^{-\frac{1}{24}}e_t-c\beta^{\frac{1}{24}}\text{ for any }t\in\lt[0,\mathfrak{A}\rt].
\end{equation}
So
\begin{eqnarray}
\label{uz10}
\lt|\Omega_{s_0}\cap H\lt(\zeta_{\thea},R\lt(Q\tau_{\thea}\rt)\rt)\cap \GI^{R_0}_{\theta}\rt|&=&
\int_{\lt[0,\mathfrak{A}\rt]}
H^1\lt(\lt[\varrho^1_t,\varrho^2_t\rt]\cap \GI^{R_0}_{\theta}\rt) dt\nn\\
&\overset{(\ref{eq501})}{\geq}&\int_{\lt[0,\mathfrak{A}\rt]}
\lt|\varrho^1_t-\varrho^2_t\rt|-c\beta^{-\frac{1}{24}} e_t-c\beta^{\frac{1}{24}} dt\nn\\
&\overset{(\ref{eqe33})}{\geq}&
\lt|\Omega_{s_0}\cap H\lt(\zeta_{\thea},R\lt(Q\tau_{\thea}\rt)\rt)\rt|
-c\beta^{\frac{1}{24}}\nn\\
&~&-c\beta^{-\frac{1}{24}}\int_{\UI} \lt| \na w_{\thea}-R\lt(\Lambda_{\thea}\lt(\ti{m}\rt)\rt)\rt| dx\nn\\
&\overset{(\ref{eqe34})}{\geq}& \lt|\Omega_{s_0}\cap H\lt(\zeta_{\thea},R\lt(Q\tau_{\thea}\rt)\rt) \rt|
-c\beta^{\frac{1}{24}}.
\end{eqnarray}
Note $\lt|\UI\backslash \Omega_{s_0}\rt|\leq c\beta^{\frac{1}{24}}$ and 
by definition of $\WI_{\thea}$ (see (\ref{eq751})) 
$\lt|\WI_{\thea}\backslash H\lt(\zeta_{\thea},R\lt(Q\tau_{\thea}\rt)\rt)\rt|\leq c\beta^{\frac{1}{24}}$ 
this together with (\ref{uz10}) gives $\lt|\WI_{\thea}\backslash \GI_{\theta}^{R_0}\rt|\leq
c\beta^{\frac{1}{24}}$. Now if $R_0=R$ and so $\thea=\theta$, it is imediate that
$\tau_{\thea}=\frac{b_{\theta}-a_{\theta}}{\lt|b_{\theta}-a_{\theta}\rt|}$ and so (again recalling 
definition (\ref{eq751})) (\ref{eq104}) follows.
On the other hand if $R_0=R^{-1}$ then $\thea=-\theta$ and so $a_{\thea}=b_{\theta}$, $b_{\thea}=a_{\theta}$,
which implies $\tau_{\thea}=-\frac{b_{\theta}-a_{\theta}}{\lt|b_{\theta}-a_{\theta}\rt|}$ so
$R\tau_{\thea}=R\lt(-\frac{b_{\theta}-a_{\theta}}{\lt|b_{\theta}-a_{\theta}\rt|}\rt)=
R^{-1}\lt(\frac{b_{\theta}-a_{\theta}}{\lt|b_{\theta}-a_{\theta}\rt|}\rt)
=R_0\lt(\frac{b_{\theta}-a_{\theta}}{\lt|b_{\theta}-a_{\theta}\rt|}\rt)$ hence (again recalling definition (\ref{eq751})),(\ref{eq104}) also follows in this case. $\Box$

%
%
%
%

\begin{a1}
\label{L5}
Let $\Omega$ be a convex body with $\mathrm{diam}(\Omega)=2$. Let
$u:W^{2,2}(\Omega)\rightarrow \R$ be a function satisfying (\ref{eq1}) and (\ref{eq2}) and in addition
$u$ satisfies $u=0$ on $\partial \Omega$ and $\na u(z)\cdot \eta_z=1$ on $\partial \Omega$ in the sense of trace where $\eta_z$ is the
inward pointing unit normal to $\partial\Omega$ at $z$.
Let $a,b\in\Omega$ be such that $\mathrm{diam}\lt(\Omega\rt)=\lt|a-b\rt|$. 
We will show there exists constant $\CI_{3}>1$ and $r_0\in (\CI_{3}^{-1}\beta^{\frac{1}{512}},
\CI_{3}\beta^{\frac{1}{512}})$ such that
\begin{equation}
\label{eq30}
u\lt(x\rt)\geq 1-\CI_{3}\beta^{\frac{1}{512}}\text{ for any }x\in \partial B_{r_0}\lt(\frac{a+b}{2}\rt).
\end{equation}
\end{a1}
\em Proof of Lemma 4. \rm Let $\UI$ be the convex set and $\ui$ be the function constructed in Lemma
\ref{L3}. To simplify our notation we will without loss of generality assume that $\frac{a+b}{2}=0$. 
It is easy to see we can chose $\ai,\bi\in\UI$ such that $\frac{\ai-\bi}{\lt|\ai-\bi\rt|}
=\frac{a-b}{\lt|a-b\rt|}$, $\lt|\ai-\bi\rt|=\mathrm{diam}\lt(\UI\rt)$ and $\frac{\ai+\bi}{2}=0$. Without loss of 
generality also assume $\frac{\ai-\bi}{\lt|\ai-\bi\rt|}=e_2$. For any $z\in \partial \UI$ let $\eta_z$ denote the inward pointing 
unit normal to $\partial \UI$ at $z$. Note that $\eta_{\ai}=-e_2$ since otherwise 
$\UI\not\subset B_{\lt|\ai-\bi\rt|}(\bi)$ and this contradicts the fact that $\lt|\ai-\bi\rt|=\mathrm{diam}(\UI)$. For the same reason 
$\eta_{\bi}=e_2$.\nl

\begin{figure}[h]
\centerline{\psfig{figure=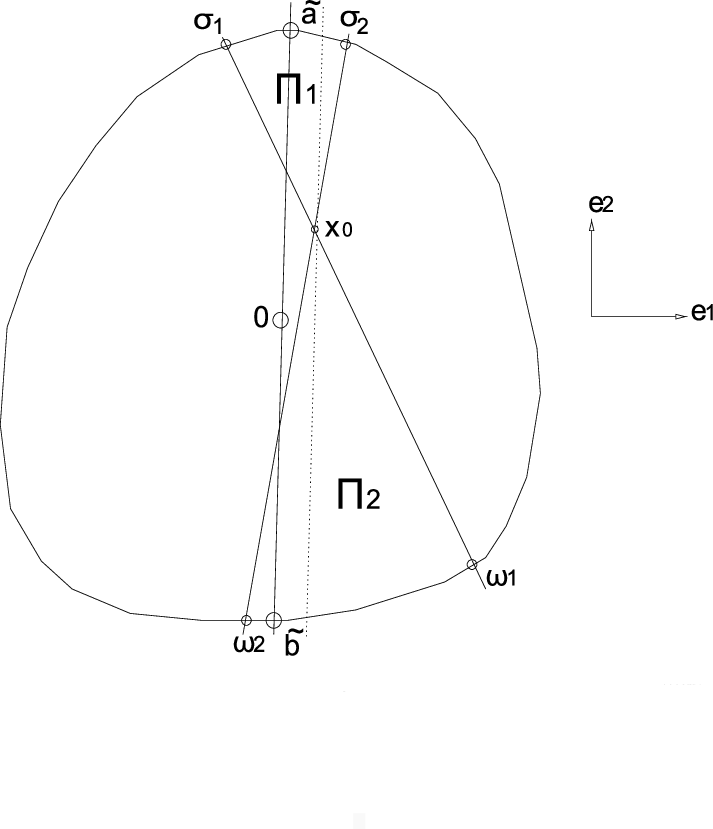,height=14cm,angle=0}}
\caption{}\label{fig4}
\end{figure}

\em Step 1. \rm Let $P:\lt[0,H^1(\partial \UI)\rt)\rightarrow \partial \UI$ be a 
`clockwise' parameterisation of $\partial \UI$ by arclength with $P(0)=\ti{a}$. For 
some $\gamma_1\in (H^1(\partial \UI)-2\beta^{\frac{1}{512}}, H^1(\partial \UI)-\beta^{\frac{1}{256}})$ 
and  $\gamma_2\in (\beta^{\frac{1}{256}}, 2\beta^{\frac{1}{512}})$ we have that for 
$\sigma_1=P(\gamma_1)$, $\sigma_2=P(\gamma_2)$, (see figure \ref{fig4}) 
the points $\sigma_1,\sigma_2$ satisfy the following properties. Firstly 
\begin{equation}
\label{ula700}
\eta_{\sigma_i}\in \wt{\Gamma}\text{ and }\eta_{\sigma_i}\cdot (-e_2)\geq 1-c\beta^{\frac{1}{128}}\text{ for }i=1,2.
\end{equation}
Secondly 
\begin{equation}
\label{ffgg19.7}
\lt|\sigma_1-\sigma_2\rt|\leq 40\beta^{\frac{1}{512}}.
\end{equation}
Thirdly 
\begin{equation}
\label{ula702}
\sigma_1\cdot (-e_1)\geq \frac{\beta^{\frac{1}{256}}}{2}\text{ and }
\sigma_2\cdot e_1\geq \frac{\beta^{\frac{1}{256}}}{2}.
\end{equation}
\em Proof of Step 1. \rm Recall $\UI=\Omega_{\frac{1}{10}}(\Omega)$, so for any 
$x\in \partial \UI$ let $z_x\in\partial \Omega$ be such that $d(x,\Omega)=\lt|x-z_x\rt|$, note 
that we can inscribe 
a ball $B_{\frac{1}{10}}(z_x)\subset \UI$ with $x\in \partial B_{\frac{1}{10}}(z_x)\cap 
\partial\UI$ and $B_{\frac{1}{10}}(z_x)\cap \partial\UI=\emptyset$. Thus the curvature of $\partial\UI$ is bounded 
above by $10$ and so 
\begin{equation}
\label{fff1}
\|\ddot{P}\|_{L^{\infty}(\partial\UI)}\leq 10. 
\end{equation}
Let $\wt{\Gamma}\subset S^1$ be the set constructed in Lemma \ref{L3}. We will show 
\begin{equation}
\label{fineqq50}
\inf\lt\{h\in\lt[\beta^{\frac{1}{256}},H^1(\partial \UI)\rt]:\eta_{P(h)}\in \wt{\Gamma}\rt\}
\leq 2\beta^{\frac{1}{512}}.
\end{equation}
Suppose this is not true, so for every $h\in \lt[\beta^{\frac{1}{256}},2\beta^{\frac{1}{512}}\rt]$, 
$\eta_{P(h)}\not\in \wt{\Gamma}$. Note that since $\partial \UI$ is $C^1$, 
$\lt\{\eta_{P(h)}:h\in \lt[\beta^{\frac{1}{256}},2\beta^{\frac{1}{512}}\rt]\rt\}$ is connected 
and since $H^1(S^1\backslash \wt{\Gamma})\leq 40\pi \beta^{\frac{1}{8}}$, so 
\begin{equation}
\label{fff1.5}
H^1\lt(\lt\{\eta_{P(h)}:h\in\lt[\beta^{\frac{1}{256}},2\beta^{\frac{1}{512}}\rt]\rt\}\rt)\leq 
40 \pi\beta^{\frac{1}{8}}.
\end{equation}

Note that as $P(0)=e_2$, $\dot{P}(0)=e_1$ and as generally for $x\in \lt[0,H^1(\partial \UI)\rt]$,  
$\dot{P}(x)=R(\eta_{P(x)})$ 
so for any $h\in \lt[0,2\beta^{\frac{1}{512}}\rt]$, 
\begin{eqnarray}
\label{ffd1}
\lt|\dot{P}(h)-e_1\rt|&\leq&\lt|\dot{P}(\beta^{\frac{1}{256}})-\dot{P}(0)\rt|+\lt|\dot{P}(h)-\dot{P}(\beta^{\frac{1}{256}})\rt|\nn\\
&\overset{(\ref{fff1.5}),(\ref{fff1})}{\leq}&20 \beta^{\frac{1}{256}}+40\pi \beta^{\frac{1}{8}}\leq 40\pi\beta^{\frac{1}{256}}.
\end{eqnarray}
So by the fundamental theorem of Calculus, $\lt|P(2\beta^{\frac{1}{512}})-(\ti{a}+2\beta^{\frac{1}{512}}e_1)\rt|\leq 80 
\pi\beta^{\frac{1}{256}}\beta^{\frac{1}{512}}$. Now  
\begin{eqnarray}
\lt|(\ti{a}+2\beta^{\frac{1}{512}}e_1)-\ti{b}\rt|&=&\sqrt{\lt|\ti{a}-\ti{b}\rt|^2+4\beta^{\frac{1}{256}}}\nn\\
&\geq&\lt|\ti{a}-\ti{b}\rt|+\frac{3}{4}\beta^{\frac{1}{256}}.\nn
\end{eqnarray} 
Thus $\lt|P(2\beta^{\frac{1}{512}})-\ti{b}\rt|\geq \lt|\ti{a}-\ti{b}\rt|+\frac{\beta^{\frac{1}{256}}}{2}$ which is a contradiction. Thus we have established (\ref{fineqq50}).

Hence (recalling the fact $H^1(S^1\backslash \wt{\Gamma})\leq 40\pi\beta^{\frac{1}{8}}$) we can 
pick $\gamma_2\in \lt[\beta^{\frac{1}{256}},2\beta^{\frac{1}{512}}\rt]\cap \wt{\Gamma}$ such that 
\begin{equation}
\label{fff3}
\lt|\eta_{P(\beta^{\frac{1}{256}})}-\eta_{P(\gamma_2)}\rt|\leq 50\pi\beta^{\frac{1}{8}}
\end{equation} 
and $\eta_{P(\gamma_2)}\in \wt{\Gamma}$. In the same way we can pick $\gamma_1\in \lt[H^1(\partial \UI)-2\beta^{\frac{1}{512}},
H^1(\partial \UI)-\beta^{\frac{1}{256}}\rt]$ such that 
$\lt|\eta_{P(H^1(\partial \UI)-\beta^{\frac{1}{256}} )}-\eta_{P(\gamma_1)}\rt|\leq 50\pi\beta^{\frac{1}{8}}$ 
and $\eta_{P(\gamma_1)}\in \wt{\Gamma}$. 

Define $\sigma_2=P(\gamma_2)$ and $\sigma_1=P(\gamma_1)$. Since $\dot{P}(0)=e_1$ and recalling again that 
$\eta_{P(s)}=R^{-1}(\dot{P}(s))$, 
\begin{eqnarray}
\lt|\dot{P}(0)-\dot{P}(\gamma_2)\rt|&\leq& \lt|\dot{P}(0)-\dot{P}(\beta^{\frac{1}{256}})\rt|+\lt|\dot{P}(\beta^{\frac{1}{256}})-\dot{P}(\gamma_2)\rt|\nn\\
&\overset{(\ref{fff1}),(\ref{fff3})}{\leq}&60 \pi\beta^{\frac{1}{256}}.\nn
\end{eqnarray}
Arguing in the same way we can establish $\lt|\dot{P}(0)-\dot{P}(\gamma_1)\rt|\leq 60 \pi\beta^{\frac{1}{256}}$. Thus as 
$\partial \UI$ is convex $\lt|\eta_{\sigma_i}+e_2\rt|\leq 60\pi\beta^{\frac{1}{256}}$ for $i=1,2$ which establishes 
(\ref{ula700}). Hence
$$
\sigma_2\cdot e_1=\lt(\sigma_2-\ti{a}\rt)\cdot e_1=\int_{0}^{\gamma_2} \dot{P}(s)\cdot e_1 ds
\overset{(\ref{ffd1})}{\geq} (1-40\pi\beta^{\frac{1}{256}})\gamma_2\geq \frac{\beta^{\frac{1}{256}}}{2}
$$
which establishes (\ref{ula702}) for $\sigma_2$. Inequality (\ref{ula702}) for $\sigma_1$ can be 
established in the same way. Finally note 
\begin{equation}
\label{ffgg19}
\lt|\sigma_1-\sigma_2\rt|=\lt|P(\gamma_2)-P(\gamma_1)\rt|\leq \int_{\gamma_2}^{H^1(\partial \UI)} \lt|\dot{P}(z)\rt| dz
+\int_{0}^{\gamma_2}  \lt|\dot{P}(z)\rt| dz\leq 40\beta^{\frac{1}{512}}
\end{equation}
which establishes (\ref{ffgg19.7}). \nl

\em Step 2. \rm For $y\in \R^2$, $\psi\in \R^2$, $\gamma>0$ define 
$X(y,\psi,\gamma):=\lt\{z:\lt|\frac{z-y}{\lt|z-y\rt|}\cdot \lt(\frac{\psi}{\lt|\psi\rt|}\rt)^{\perp}\rt|\leq \gamma\rt\}$. We will show there exists positive constant $\CI_4$ and 
$x_0\in N_{\CI_4\beta^{\frac{1}{512}}}\lt(\lt[\ai,\bi\rt]\rt)\cap \UI$ such
that for some $\psi_0\in B_{\CI_4\beta^{\frac{1}{256}}}(e_2)$ the following inequality holds
\begin{equation}
\label{eq120}
\lt|X\lt(x_0,\psi_0,\CI_4\beta^{\frac{1}{256}}\rt)\cap \UI\backslash \lt\{x:\lt|\na \ui\lt(x\rt)\cdot e_1\rt|<
\CI_4\beta^{\frac{1}{256}}\rt\}\rt|\leq \CI_4\beta^{\frac{1}{24}}.
\end{equation}
\em Proof of Step 2. \rm Recall we know $\sigma_1$ and $\sigma_2$ are chosen 
so that $\eta_{\sigma_1}\in \wt{\Gamma}$ and $\eta_{\sigma_2}\in \wt{\Gamma}$. We also 
know $\eta_{\ai}=-e_2$ and $\eta_{\bi}=e_2$. Let $\omega_1\in\partial \UI$ be the
unique point for which $-\eta_{\omega_1}=\eta_{\sigma_1}$ and let
$\omega_2\in\partial \UI$ be the
unique point for which $-\eta_{\omega_2}=\eta_{\sigma_2}$, see figure \ref{fig4}.

Define
\begin{equation}
\label{eq507}
\Pi_2:=H\lt(\frac{\sigma_2+\omega_2}{2}, R\lt(\frac{\omega_2-\sigma_2}{\lt|\omega_2-\sigma_2\rt|}\rt)\rt)\cap 
H\lt(\frac{\sigma_1+\omega_1}{2}, R^{-1}\lt(\frac{\omega_1-\sigma_1}{\lt|\omega_1-\sigma_1\rt|}\rt)\rt)
\end{equation}
and
\begin{equation}
\label{eq508}
\Pi_1:=H\lt(\frac{\sigma_2+\omega_2}{2}, R^{-1}\lt(\frac{\omega_2-\sigma_2}{\lt|\omega_2-\sigma_2\rt|}\rt)\rt)
\cap H\lt(\frac{\sigma_1+\omega_1}{2}, R\lt(\frac{\omega_1-\sigma_1}{\lt|\omega_1-\sigma_1\rt|}\rt)\rt)
\end{equation}
and let $\Pi=\Pi_1\cup \Pi_2$ and let $x_0:=\overline{\Pi_1}\cap \overline{\Pi_2}$,
see again figure \ref{fig4}.

Let us define $l_x^{\theta}:=x+\R_{+}\theta$ for any $x\in \R^2$, $\theta\in S^1$. 
First we will show $\lt(x_0+\R e_2\rt)\subset \Pi$
however this inclusion is relatively easy to see because firstly
$$
e_2\cdot R\lt(\frac{\omega_1-\sigma_1}{\lt|\omega_1-\sigma_1\rt|}\rt)=
e_1\cdot \lt(\frac{\omega_1-\sigma_1}{\lt|\omega_1-\sigma_1\rt|}\rt)
\overset{(\ref{ula702})}{\geq} \frac{10\beta^{\frac{1}{256}}}{44}
$$
thus $l^{e_2}_{0}\subset H\lt(0,R\lt(\frac{\omega_1-\sigma_1}{\lt|\omega_1-\sigma_1\rt|}\rt)\rt)$. And
secondly as $x_0\in\partial H\lt(\frac{\sigma_1+\omega_1}{2}, R\lt(\frac{\omega_1-\sigma_1}{\lt|\omega_1-\sigma_1\rt|}\rt)\rt)$
$$
l^{e_2}_{x_0}\subset H\lt(x_0,R\lt(\frac{\omega_1-\sigma_1}{\lt|\omega_1-\sigma_1\rt|}\rt)\rt)
=H\lt(\frac{\sigma_1+\omega_1}{2},R\lt(\frac{\omega_1-\sigma_1}{\lt|\omega_1-\sigma_1\rt|}\rt)\rt).
$$
In exactly the same way $l_{x_0}^{e_2}\subset H\lt(\frac{\sigma_2+\omega_2}{2}, R^{-1}\lt(\frac{\omega_2-\sigma_2}{\lt|\omega_2-\sigma_2\rt|}\rt)\rt)$. Hence
$l_{x_0}^{e_2}\subset \Pi_1$. Arguing in the same manner we have $l_{x_0}^{-e_2}\subset \Pi_2$
and thus we have established the claim.

Let $\gamma=l_{x_0}^{e_2}\cap \partial\UI$, by construction we have that $\gamma$ lies in the component of
$\partial\UI$ between $\sigma_1$ and $\sigma_2$ and hence by (\ref{ffgg19}) we know
$d\lt(\gamma,l^{e_2}_{0}\rt)\leq 40\beta^{\frac{1}{512}}$ and so it follows
$x_0\in N_{c\beta^{\frac{1}{512}}}\lt(\lt[\ai,\bi\rt]\rt)\cap\UI$.

Since $\eta_{\ai}=-e_2$, $\eta_{\bi}=e_2$ and $\UI$ is convex we know $\omega_2\in H\lt(0,-e_1\rt)$ and for
the same reasons $\omega_1\in H\lt(0,e_1\rt)$ see figure \ref{fig4}. So
$\lt(\sigma_2-\omega_2\rt)\cdot e_1\geq \sigma_2\cdot e_1\overset{(\ref{ula702})}{\geq} c\beta^{\frac{1}{256}}$
and for exactly the same reason $\lt(\sigma_1-\omega_1\rt)\cdot \lt(-e_1\rt)\geq \sigma_1\cdot
\lt(-e_1\rt)\overset{(\ref{ula702})}{\geq} c\beta^{\frac{1}{256}}$. Thus as $\lt|\sigma_1-\omega_1\rt|\leq 2\mathrm{diam}\lt(\UI\rt)$ and  $\lt|\sigma_2-\omega_2\rt|\leq 2\mathrm{diam}\lt(\UI\rt)$ we have $\frac{\sigma_2-\omega_2}{\lt|\sigma_2-\omega_2\rt|}
\cdot e_1\geq c\beta^{\frac{1}{256}}$ and
$\frac{\sigma_1-\omega_1}{\lt|\sigma_1-\omega_1\rt|}\cdot \lt(-e_1\rt)\geq c\beta^{\frac{1}{256}}$.
Hence
\begin{eqnarray}
\lt(\frac{\sigma_1-\omega_1}{\lt|\sigma_1-\omega_1\rt|}\rt)\cdot
\lt(\frac{\sigma_2-\omega_2}{\lt|\sigma_2-\omega_2\rt|}\rt)&=&
\lt(\frac{\sigma_1-\omega_1}{\lt|\sigma_1-\omega_1\rt|}\cdot e_1\rt)\lt(\frac{\sigma_2-\omega_2}{\lt|\sigma_2-\omega_2\rt|}\cdot e_1\rt)\nn\\
&~&
+\lt(\frac{\sigma_1-\omega_1}{\lt|\sigma_1-\omega_1\rt|}\cdot e_2\rt)\lt(\frac{\sigma_2-\omega_2}{\lt|\sigma_2-\omega_2\rt|}\cdot e_2\rt)\nn\\
&\leq& -c\beta^{\frac{1}{128}}+1.\nn
\end{eqnarray}

In other words the angle between $\frac{\sigma_1-\omega_1}{\lt|\sigma_1-\omega_1\rt|}$ and
 $\frac{\sigma_2-\omega_2}{\lt|\sigma_2-\omega_2\rt|}$ is greater than
$\CI_4 \beta^{\frac{1}{256}}$ for some positive constant $\CI_4$. Thus there exists 
$\psi_0\in B_{c\beta^{\frac{1}{256}}}(e_2)$
such that $X\lt(x_0,\psi_0,\CI_4 \beta^{\frac{1}{256}}\rt)\subset \Pi$. Now 
since $\eta_{\sigma_1},\eta_{\sigma_2}\in \wt{\Gamma}$ we can apply Lemma \ref{L3} so we know that
$$
\lt|\UI\cap H\lt(\frac{\sigma_2+\omega_2}{2},R^{-1}\lt(\frac{\omega_2-\sigma_2}
{\lt|\omega_2-\sigma_2\rt|}\rt)\rt)\backslash \GI^{R^{-1}}_{\eta_{\sigma_2}}\rt|\leq c\beta^{\frac{1}{24}}
$$
and
$$
\lt|\UI\cap H\lt(\frac{\sigma_1+\omega_1}{2},R\lt(\frac{\omega_1-\sigma_1}
{\lt|\omega_1-\sigma_1\rt|}\rt)\rt)\backslash \GI^{R}_{\eta_{\sigma_1}}\rt|\leq c\beta^{\frac{1}{24}}.
$$
Thus (recalling the definition of $\Pi_1$, (\ref{eq508}))
\begin{equation}
\label{eqa1}
\lt|\Pi_1\cap \UI\backslash \GI^{R^{-1}}_{\eta_{\sigma_2}}\cap \GI^{R}_{\eta_{\sigma_1}}\rt|\leq c\beta^{\frac{1}{24}}.
\end{equation}
In exactly the same way we have (recall (\ref{eq507}))
\begin{equation}
\label{eqa2}
\lt|\Pi_2\cap \UI\backslash \GI^{R^{-1}}_{\eta_{\sigma_1}}\cap \GI^{R}_{\eta_{\sigma_2}}\rt|\leq c\beta^{\frac{1}{24}}.
\end{equation}

Now for any $x\in \GI^{R^{-1}}_{\eta_{\sigma_2}}\cap \GI^{R}_{\eta_{\sigma_1}}$ we have 
$\na \ui\lt(x\rt)\cdot R\eta_{\sigma_2}\geq 0$ and $\na \ui\lt(x\rt)\cdot R^{-1}\eta_{\sigma_1}\geq 0$. 
Since from (\ref{ula700}) $\eta_{\sigma_i}\in X^{+}\lt(0,-e_2,c\beta^{\frac{1}{256}}\rt)$
for $i=1,2$ we know $R\eta_{\sigma_2}\in X^{+}\lt(0,e_1,c\beta^{\frac{1}{256}}\rt)$ and $R^{-1}\eta_{\sigma_1}\in X^{+}\lt(0,-e_1,\beta^{\frac{1}{256}}\rt)$, from this it is easy to see (assuming 
we chose $\CI_4$ large enough)$\lt|\na \ui\lt(x\rt)\cdot e_1\rt|\leq
\CI_4\beta^{\frac{1}{256}}$. And in the same way for any $x\in \GI^{R^{-1}}_{\eta_{\sigma_1}}\cap \GI^{R}_{\eta_{\sigma_2}}$ we also have
$\lt|\na \ui\lt(x\rt)\cdot e_1\rt|\leq \CI_4\beta^{\frac{1}{256}}$.
\begin{eqnarray}
&~&\lt|X\lt(x_0,\psi_0,\CI_4\beta^{\frac{1}{256}}\rt)\cap \UI\backslash \lt\{x:\lt|\na \ui\lt(x\rt)\cdot e_1\rt|
<\CI_4\beta^{\frac{1}{256}}\rt\}\rt|\nn\\
&~&\qd\qd\qd \leq c\lt|\Pi_1\cap \UI\backslash  \GI^{R}_{\eta_{\sigma_1}}\cap \GI^{R^{-1}}_{\eta_{\sigma_2}}\rt|
+ c \lt|\Pi_2\cap \UI\backslash  \GI^{R}_{\eta_{\sigma_2}}\cap \GI^{R^{-1}}_{\eta_{\sigma_1}}\rt|\nn\\
&~&\qd\qd\qd \leq \CI_4\beta^{\frac{1}{24}}\nn
\end{eqnarray}
which establishes (\ref{eq120}).\nl

\em Step 3. \rm There exists positive constant $\CI_5$ such that for some $v_1\in\lt\{e_2,-e_2\rt\}$ we have 
\begin{eqnarray}
\label{eq13.5}
&~&\lt|X\lt(x_0,\psi_0,\CI_4\beta^{\frac{1}{256}}\rt)\cap \UI\cap H\lt(\CI_5\beta^{\frac{1}{256}}v_1,v_1\rt)\backslash \Vb_{-v_1}\rt|\leq \CI_5\beta^{\frac{1}{24}}.
\end{eqnarray}
where 
\begin{equation}
\label{ulazz5}
\Vb_{-v_1}:=\lt\{x\in \UI:\na \ui\lt(x\rt)\in N_{\CI_5\beta^{\frac{1}{256}}}\lt(-v_1\rt)\rt\}.
\end{equation}

\em Proof of Step 3. \rm Let
$\wt{\varpi_0}=l_0^{-e_1}\cap \partial \UI$. Note since $\UI$ is convex $\eta_{\wt{\varpi_0}}\cdot e_1>0$. We claim
\begin{equation}
\label{eq400}
\eta_{\wt{\varpi_0}}\cdot e_1> \frac{1}{10}.
\end{equation}
Suppose this were not the case, then
$\eta_{\wt{\varpi_0}}\cdot e_1\leq \frac{1}{10}$. Since $\UI$ is convex (and recall $\UI=\Omega_{\frac{1}{10}}$) and
$\mathrm{diam}(\UI)=\frac{22}{10}$ we know
$\UI\subset \overline{H\lt(\wt{\varpi_0},\eta_{\wt{\varpi_0}}\rt)}\subset H\lt(-\frac{22}{10}e_1,\eta_{\wt{\varpi_0}}\rt)$ which
implies $(\bi+\frac{22}{10}e_1)\cdot \eta_{\wt{\varpi_0}}>0$ and thus
\begin{eqnarray}
\bi\cdot e_2 \sqrt{\frac{99}{100}}&\geq& \lt(\lt(\bi+\frac{22}{10}e_1\rt)\cdot e_2 \rt)\lt(\eta_{\wt{\varpi_0}}\cdot e_2\rt)
>-\lt(\lt(\bi+\frac{22}{10}e_1\rt)\cdot e_1\rt)\lt(\eta_{\wt{\varpi_0}}\cdot e_1\rt)\nn\\
&=&-\frac{22}{10}\eta_{\wt{\varpi_0}}\cdot e_1\geq -\frac{22}{100}
\end{eqnarray}
however as $\lt|\ai-\bi\rt|=\mathrm{diam}(\UI)=\frac{22}{10}$, $\frac{\ai+\bi}{2}=0$ and $\frac{\ai-\bi}{\lt|\ai-\bi\rt|}=e_2$ this is a contradiction, thus (\ref{eq400}) is established. 

Let 
$$
\alpha_0=\sup\lt\{\alpha>0:\lt\{\eta_x:x\in B_{\alpha}(\wt{\varpi_0})\cap \partial \UI\rt\}\cap \wt{\Gamma}=\emptyset\rt\},
$$
in the case where $\lt\{\alpha>0:\lt\{\eta_x:x\in B_{\alpha}(\wt{\varpi_0})\cap \partial \UI\rt\}\cap \wt{\Gamma}=\emptyset\rt\}=\emptyset$ 
let $\alpha_0=0$. 

Since $H^1(S^1\backslash \wt{\Gamma})\leq 40\pi\beta^{\frac{1}{8}}$ we know 
$\partial \UI\backslash B_{\alpha_0}(\wt{\varpi_0})\not =\emptyset$. Note also 
$$
\MI_0:=\lt\{\eta_x:x\in B_{\alpha_0}(\wt{\varpi_0})\cap \partial \UI\rt\}
$$ 
is a connected 
subset of $S^1$, so $H^1(\MI_0)\leq 40\pi\beta^{\frac{1}{8}}$ hence for every 
$z\in B_{\alpha_0}(\wt{\varpi_0})\cap \partial \UI$, $\lt|\eta_{z}-\eta_{\wt{\varpi_0}}\rt|\leq 40\pi\beta^{\frac{1}{8}}$. So we can pick $\alpha_1>\alpha_0$ such that some point 
$\varpi_0\in \partial B_{\alpha_1}(\wt{\varpi_0})\cap \partial \UI$ satisfies 
$\eta_{\varpi_0}\in \wt{\Gamma}$ and 
\begin{equation}
\label{fineqq60}
\lt|\eta_z-\eta_{\varpi_0}\rt|\leq 50\pi\beta^{\frac{1}{8}}\text{ for all }
z\in B_{\alpha_1}(\wt{\varpi_0}).
\end{equation}

Now since $B_{\frac{1}{10}}(0)\subset \UI$, we know $\wt{\varpi_0}\cdot (-e_1)\geq \frac{1}{10}$. Using 
again the fact that $\eta_{P(s)}=R^{-1}(\dot{P}(s))$ (where $P$ is the parameterisation of $\partial \UI$) it is easy to see 
by the fundamental theorem of Calculus that (\ref{fineqq60}) implies 
\begin{equation}
\label{ffgg2}
\varpi_0\cdot (-e_1)\geq \frac{1}{11}.
\end{equation} 
Also from (\ref{eq400}) and (\ref{fineqq60}) we know that 
\begin{equation}
\label{fff5}
\eta_{\varpi_0}\cdot e_1>\frac{1}{11}.
\end{equation}

Let $\varpi_1\in \partial\UI$ be the unique point for which $\eta_{\varpi_1}=-\eta_{\varpi_0}$. Note that by 
(\ref{fff5}) we know that $\eta_{\varpi_1}\cdot (-e_1)>\frac{1}{11}$ and as 
$\eta_{\ai}=-e_2$ and $\eta_{\bi}=e_2$ by convexity 
of $\UI$ this implies  
\begin{equation}
\label{fineq4}
\varpi_1\in \partial \UI\cap H(0,e_1).
\end{equation}

Now let $l\in\lt(\frac{\varpi_1-\varpi_0}{\lt|\varpi_1-\varpi_0\rt|}\rt)^{\perp}\cap S^1$ be such that
\begin{equation}
\label{eq129}
H^1\lt(\lt[a,b\rt]\cap H\lt(\frac{\varpi_1+\varpi_0}{2},l\rt)\rt)\geq \frac{\lt|a-b\rt|}{2}.
\end{equation}
Choose $S\in \lt\{R^{-1},R\rt\}$ so that $S\lt(\frac{\varpi_1-\varpi_0}{\lt|\varpi_1-\varpi_0\rt|}\rt)=l$, since $\eta_{\varpi_0}\in \wt{\Gamma}$ we can apply Lemma \ref{L3} and hence 
we have 
\begin{equation}
\label{eq401}
\lt|\UI\cap H\lt(\frac{\varpi_1+\varpi_0}{2},l\rt)\backslash \GI^S_{\eta_{\varpi_0}}\rt|\leq c\beta^{\frac{1}{24}}.
\end{equation}
From (\ref{eq2}) and (\ref{eq120}) we know
\begin{equation}
\label{eq740}
\lt|X\lt(x_0,\psi_0,\CI_4 \beta^{\frac{1}{256}}\rt)\cap \UI\backslash \lt\{x:\na \ui(x)\in
N_{100^{-1}}(\lt\{e_2,-e_2\rt\})\rt\}\rt|\leq c\beta^{\frac{1}{24}}.
\end{equation}
Since so $\lt|S^{-1}\eta_{\varpi_0}\cdot e_2\rt|\overset{(\ref{fff5})}{>}11^{-1}$ there exists 
some fixed vector $v_0\in\lt\{e_2,-e_2\rt\}$ such that if
$x\in \GI^S_{\eta_{\varpi_0}}\cap
\lt\{x:\na \ui\lt(x\rt)\in N_{100^{-1}}\lt(\lt\{e_2,-e_2\rt\}\rt)\rt\}$ then 
$\na \ui\lt(x\rt)\in B_{100^{-1}}\lt(v_0\rt)$. So using (\ref{eq401}) and (\ref{eq740})
\begin{equation}
\label{eq128}
\lt|X\lt(x_0,\psi_0,\CI_4 \beta^{\frac{1}{256}}\rt)\cap \UI\cap H\lt(\frac{\varpi_1+\varpi_0}{2},l\rt)
\backslash \lt\{x:\na \ui(x)\in B_{100^{-1}}\lt(v_0\rt)\rt\}\rt|\leq c\beta^{\frac{1}{24}}.
\end{equation}
Now for any $w\in H\lt(0,v_0\rt)$ we have the elementary inequality $\lt|w-v_0\rt|\leq 4 d(w,S^1)
+2\lt|w\cdot e_1\rt|$, so using (\ref{eq2}), (\ref{eq120}) and (\ref{eq128}) we have (assuming constant $\CI_5$ is large enough, recall definition 
(\ref{ulazz5}))
\begin{equation}
\label{eq37}
\lt|X\lt(x_0,\psi_0,\CI_4 \beta^{\frac{1}{256}}\rt)\cap \UI\cap H\lt(\frac{\varpi_1+\varpi_0}{2},l\rt)\backslash
\Vb_{v_0}\rt|\leq c\beta^{\frac{1}{24}}.
\end{equation}
By (\ref{fineq4}) $\varpi_1\cdot e_1\geq 0$ and so $\lt|\frac{\varpi_1-\varpi_0}{\lt|\varpi_1-\varpi_0\rt|}\cdot e_1\rt|\overset{(\ref{ffgg2})}{\geq} \frac{1}{44}$ and so $\lt|l\cdot e_2\rt|\geq \frac{1}{44}$. Thus by the fact that $\psi_0\in B_{\CI_4\beta^{\frac{1}{256}}}(e_2)$ and
that inequality (\ref{eq129}) implies $0\in\overline{H(\frac{\varpi_1+\varpi_0}{2},l)}$ there 
exists $v_1\in\lt\{e_2,-e_2\rt\}$ such that for some constant $\CI_5$ we have 
\begin{equation}
\label{uz11}
X\lt(x_0,\psi_0,\CI_4 \beta^{\frac{1}{256}}\rt)\cap H\lt(\CI_5\beta^{\frac{1}{256}}v_1,v_1\rt)\subset
H\lt(\frac{\varpi_1+\varpi_0}{2},l\rt).
\end{equation}
Putting (\ref{uz11}) together with (\ref{eq37}) gives
$$
\lt|X\lt(x_0,\psi_0,\CI_4\beta^{\frac{1}{256}}\rt)\cap\UI\cap H\lt(\CI_5\beta^{\frac{1}{256}}v_1,v_1\rt)\backslash
\Vb_{v_0}\rt|\leq c\beta^{\frac{1}{24}}.
$$

Let $x\in \UI\backslash\overline{\Omega}\cap X\lt(x_0,\psi_0,\CI_4\beta^{\frac{1}{256}}\rt)\cap
H(\CI_5\beta^{\frac{1}{256}}v_1,v_1)$ so as $\ui(x)=d(x,\partial \UI)$ (and since again 
$\psi_0\in B_{\CI_4\beta^{\frac{1}{256}}}\lt(e_2\rt)$) so $\na \ui(x)\in
N_{\CI_5\beta^{\frac{1}{256}}}(-v_1)$ thus we must have $v_0=-v_1$, this gives (\ref{eq13.5}).\nl

\begin{figure}[h]
\centerline{\psfig{figure=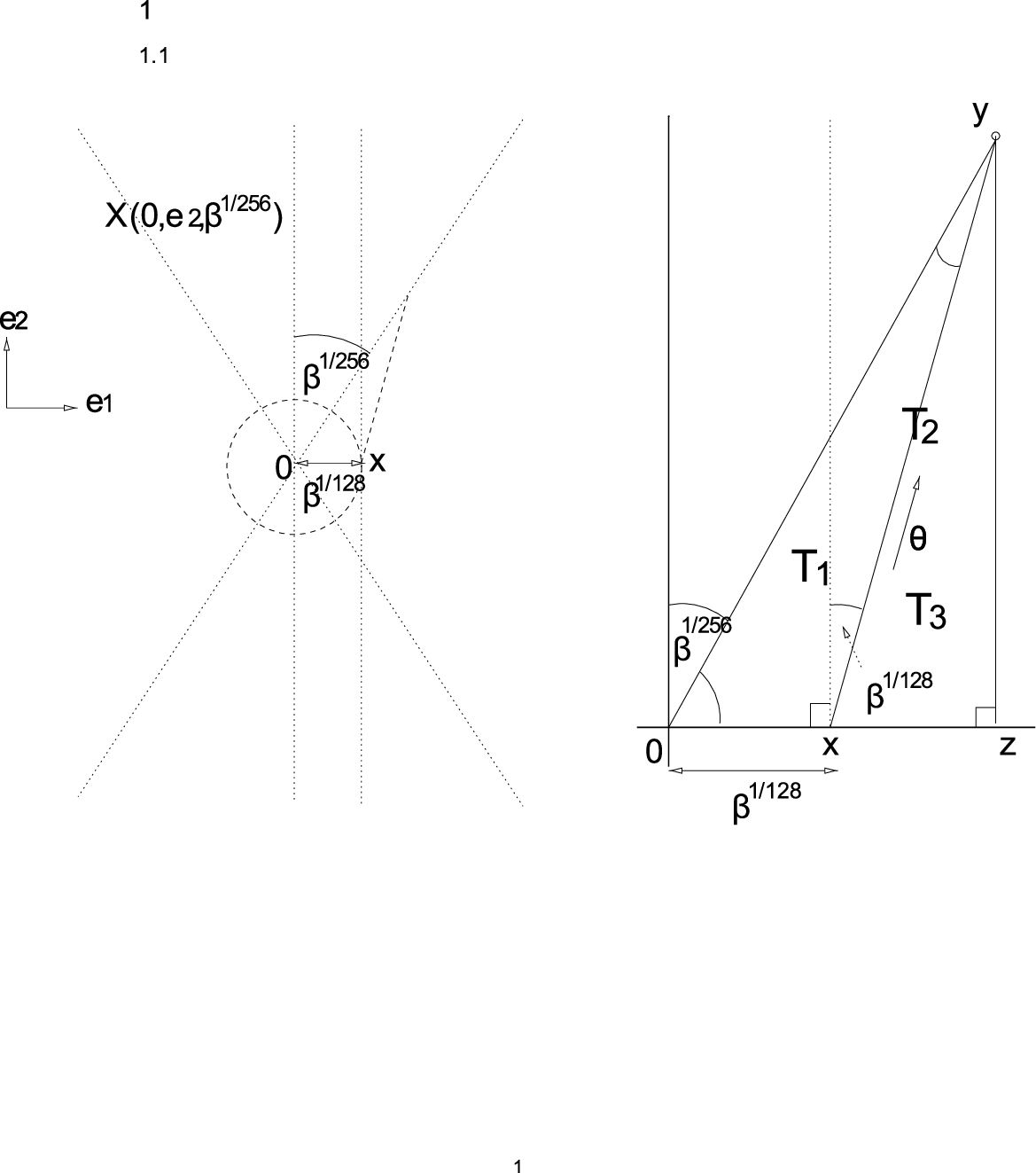,height=15cm,angle=0}}
\caption{}\label{fig6}
\end{figure}

\em Step 4. \rm We will show there exists a positive constant $\CI_6$ such that
\begin{equation}
\label{eq515}
\lt(l_x^{-\theta}\cup l_x^{\theta}\rt)\backslash B_{\CI_6 \beta^{\frac{1}{256}}}(x)\subset
X(x_0,\psi_0,\CI_4\beta^{\frac{1}{256}})\text{ for all } x\in B_{\beta^{\frac{1}{128}}}(x_0), \theta\in S^1\cap B_{\beta^{\frac{1}{128}}}(\psi_0).
\end{equation}
\em Proof of Step 4. \rm Without loss of generality we assume $x_0=0$, $\psi_0=e_2$ and
$\CI_4=1$. To begin with to take point
$x=\beta^{\frac{1}{128}}e_1$, we will show later the general case follows from this. See
figure \ref{fig6}.

Let
$\theta=\lt(\begin{matrix} \sin\beta^{\frac{1}{128}}\\ \cos\beta^{\frac{1}{128}}\end{matrix}\rt)$ and let
$y=\partial X(0,e_2,\beta^{\frac{1}{256}})\cap l_x^{\theta}$. We will get an upper bound on $\lt|y\rt|$.
Let $z=y\cdot e_1 e_1$. We have two triangles to calculate with, triangle $T_1$ with corners on
$0,x,y$ which is a subset of triangle $T_2$ with corners on $0,z,y$. Note that by applying the
law of sins we have $\lt|y\rt|^{-1}\sin(\frac{\pi}{2}+\beta^{\frac{1}{128}})=\lt|x-y\rt|^{-1}\sin(\frac{\pi}{2}-\beta^{\frac{1}{256}})$.
Note that $T_3=T_2\backslash T_1$ is also a right angle triangle and since
$\lt|z\rt|=\beta^{\frac{1}{128}}+\lt|x-z\rt|$ we have
$\lt|y\rt|\cos(\frac{\pi}{2}-\beta^{\frac{1}{256}})=\beta^{\frac{1}{128}}+\lt|y-x\rt|
\cos(\frac{\pi}{2}-\beta^{\frac{1}{128}})$.
Putting this together with the previous equation we have $\lt|y\rt|\sin \beta^{\frac{1}{256}}=
\beta^{\frac{1}{128}}+\lt|y\rt|\frac{\cos \beta^{\frac{1}{256}}}{\cos\beta^{\frac{1}{128}}}
\sin\beta^{\frac{1}{128}}$ which gives
$\lt|y\rt|\lt(\sin \beta^{\frac{1}{256}}-
\frac{\cos \beta^{\frac{1}{256}}}{\cos\beta^{\frac{1}{128}}}\sin\beta^{\frac{1}{128}}\rt)=\beta^{\frac{1}{128}}$.
Now by taking the Taylor series approximating $\sin$ and $\cos$ we have
$\lt|y\rt|\lt(\beta^{\frac{1}{256}}+O\lt(\beta^{\frac{1}{128}}\rt)\rt)=\beta^{\frac{1}{128}}$. Thus
$\lt|y\rt|\sim \beta^{\frac{1}{256}}$ and thus the existence of constant $\CI_6$ such that
(\ref{eq515}) holds follows instantly for the case $x=\beta^{\frac{1}{128}}e_1$.

In the general case where $x\not=\beta^{\frac{1}{128}}e_1$ suppose without loss of generality $x\cdot e_1>0$, 
define $\ti{x}=\lt(x+\la \theta\ra\rt)\cap\la e_1\ra$, since the angle between $\theta$ and $e_1$
is with $c\beta^{\frac{1}{256}}$ of $\frac{\pi}{2}$ it is easy to see
$\ti{x}\in B_{2\beta^{\frac{1}{128}}}(0)$ and of course $l^{\theta}_{\ti{x}}\cap \partial X(0,e_2,
\beta^{\frac{1}{256}})=l^{\theta}_{x}\cap \partial X(0,e_2,\beta^{\frac{1}{256}})$ so the
argument for the special case $x=\beta^{\frac{1}{128}}e_1$ can be applied to show the
existence of constant $\CI_6$ satisfying (\ref{eq515}).\nl

\em Step 5. \rm We will establish (\ref{eq30}). \nl
\em Proof of Step 5. \rm Let 
\begin{equation}
\label{ula901}
h(z):=\cha_{X\lt(x_0,\psi_0,\CI_4\beta^{\frac{1}{256}}\rt)\cap H\lt(\CI_5\beta^{\frac{1}{256}}v_1,v_1\rt)\cap \UI\backslash
\Vb_{-v_1}} 
\end{equation}
so we know $\int h\overset{(\ref{eq13.5})}{\leq} c\beta^{\frac{1}{24}}$.
So by the Fubini's Theorem
\begin{eqnarray}
\label{eq45}
&~&
\int_{\UI}\int_{\UI} \lt(h(z)+\beta^{-1}\lt|1-\lt|\na \ui(z)\rt|^2\rt|^2\rt)\lt|z-x\rt|^{-1} dz dx\nn\\
&~&\qd\qd\qd\qd\leq
\int_{\UI} \lt(h(z)+ \beta^{-1}\lt|1-\lt|\na \ui(z)\rt|^2\rt|^2\rt)\lt(\int_{\UI} \lt|z-x\rt|^{-1} dx\rt) dz\nn\\
&~&\qd\qd\qd\qd\leq    c\int_{\UI}  \lt(h(z)+ \beta^{-1}\lt|1-\lt|\na \ui(z)\rt|^2\rt|^2\rt) dz\nn\\
&~& \qd\qd\qd\qd\overset{(\ref{eq2})}{\leq}  c \beta^{\frac{1}{24}}.
\end{eqnarray}
Let
\begin{equation}
\label{ulzu1}
G:=\lt\{x\in B_{\beta^{\frac{1}{128}}}\lt(x_0\rt):\int_{\UI} 
\lt(h(z)+\beta^{-1}\lt|1-\lt|\na \ui(z)\rt|^2\rt|^2\rt)\lt|z-x\rt|^{-1} dz\leq \beta^{\frac{1}{48}}\rt\}
\end{equation}
so we know $\beta^{\frac{1}{48}}\lt|B_{\beta^{\frac{1}{128}}}(x_0)\backslash G\rt|\leq c\beta^{\frac{1}{24}}$,
thus $\lt| B_{\beta^{\frac{1}{128}}}(x_0)\backslash G\rt|\leq c\beta^{\frac{1}{48}}$, 
assuming $\beta$ is small enough $\lt| G\rt|
\geq 2^{-1}\beta^{\frac{1}{64}}$.
By Step 4, (\ref{eq515}) for any
$x\in B_{\beta^{\frac{1}{128}}}(x_0)$, $\theta\in B_{\beta^{\frac{1}{128}}}(\psi_0)\cap S^1$ we have
$(l^{-\theta}_x\cup l^{\theta}_x)\backslash B_{\CI_6\beta^{\frac{1}{256}}}(x)\subset X(x_0,\psi_0,\CI_4\beta^{\frac{1}{256}})$.

Since $X(x_0,\psi_0,\CI_4\beta^{\frac{1}{256}})=X(x_0,-\psi_0,\CI_4\beta^{\frac{1}{256}})$ we can assume without loss
generality that $\psi_0\cdot v_1>0$. Pick $x\in G$, by the Co-area formula we must be able to find
$\theta_1\in B_{\beta^{\frac{1}{128}}}\lt(\psi_0\rt)\cap S^1$
such that
\begin{equation}
\label{eq960}
\int_{(l_{x}^{-\theta_1}\cup l_{x}^{\theta_1})\cap \UI} h(z)+\beta^{-1}\lt|1-\lt|\na \ui(z)\rt|^2\rt|^2 dH^1 z\leq c\beta^{\frac{1}{48}}/\beta^{\frac{1}{128}}\leq c\beta^{\frac{1}{128}}.
\end{equation}

Let   
$\KI:=(l_{x}^{-\theta_1}\cup l_{x}^{\theta_1})\cap \UI\cap H(\CI_5 \beta^{\frac{1}{256}}v_1,v_1)$. Let $d,e$ be 
the endpoint of $\KI$ where we chose $d\in\partial H(\CI_5 \beta^{\frac{1}{256}}v_1,v_1)$ and $e\in\partial \UI$. As already 
noted, by Step 4 $\KI\backslash B_{\CI_6 \beta^{\frac{1}{256}}}(x)\overset{(\ref{eq515})}{\subset} X(x_0,\psi_0,\CI_4\beta^{\frac{1}{256}})\cap H(\CI_5\beta^{\frac{1}{256}} v_1,v_1)\cap \UI$, so for any $z\in \KI\backslash B_{\CI_6 \beta^{\frac{1}{256}}}(x)$ with $h(z)=0$ by definition (\ref{ula901}) we must have 
$z\in \Vb_{-v_1}$ so 
\begin{equation}
\label{ulzz1}
H^1\lt(\KI\backslash \Vb_{-v_1}\rt)\leq 4\CI_6 \beta^{\frac{1}{256}}+H^1\lt(\KI
\backslash\lt(B_{\CI_6 \beta^{\frac{1}{256}}}(x)\cup \Vb_{-v_1}\rt)\rt)\overset{(\ref{eq960})}{\leq} c\beta^{\frac{1}{256}}.
\end{equation}

Note also that if $z\in \Vb_{-v_1}$ so $\na\ui(z)\in B_{\CI_5\beta^{\frac{1}{256}}}(-v_1)$ and as 
(recall from Step 2, $\lt|\psi_0\cdot e_1\rt|<c\beta^{\frac{1}{256}}$ and we 
assumed without loss of generality $\psi_0\cdot v_1>0$) 
\begin{equation}
\label{ffe1}
\theta_1\in B_{\beta^{\frac{1}{128}}}(\psi_0)\subset B_{2\CI_4\beta^{\frac{1}{256}}}(v_1)
\end{equation}
thus 
$\na\ui(z)\cdot (-\theta_1)\geq 1+\frac{\lt|\na \ui(z)\rt|^2-1}{2}-c\beta^{\frac{1}{128}}$. Now for $z\in \KI$ let $t_z$ denote the tangent 
to $\KI$, since $t_z=-\theta_1$ so by the fundamental theorem of Calculus 
\begin{eqnarray}
\label{eq273}
\ui\lt(d\rt)-\ui\lt(e\rt)&\geq&\int_{\Vb_{-v_1}\cap \KI} \na \ui(z) \cdot(-\theta_1) dH^1 z-
\int_{\KI\backslash \Vb_{-v_1}} \lt|\na \ui(z)\rt| dH^1 z\nn\\
&\geq&\lt(1-c\beta^{\frac{1}{128}}\rt)H^1\lt(\Vb_{-v_1}\cap \KI\rt)
-H^1\lt(\KI\backslash\Vb_{-v_1}\rt)\nn\\
&~&-c\int_{\KI} \lt|1-\lt|\na \ui\rt|^2\rt| dH^1 z\nn\\
&\overset{(\ref{eq960}),(\ref{ulzz1})}{\geq}& \lt|d-e\rt|(1-c\beta^{\frac{1}{256}}).
\end{eqnarray}

Since the curvature of $\partial \UI$ is bounded above by $10$ and by 
(\ref{ffe1}) it is easy to see either $e$ is very close to $\ti{a}$ or $\ti{b}$, 
we will without loss of generality assume the former, so by (\ref{ffe1}) we have  
\begin{equation}
\label{eq274}
\lt|e-\ti{a}\rt|\leq c\beta^{\frac{1}{256}},
\end{equation}
it is also easy to see $\lt[e,\ti{a}\rt]\subset \UI\backslash \Omega$ and $\ui$ is $1$-Lipschitz
on $\UI\backslash\Omega$ so
\begin{equation}
\label{eq275}
\lt|\ui(e)-\ui(\ti{a})\rt|\leq c\beta^{\frac{1}{256}}.
\end{equation}
Note also as $d\in \partial H(\CI_5 \beta^{\frac{1}{256}}v_1,v_1)\cap \lt(l_x^{-\theta_1}\cup l_x^{\theta_1}\rt)$ by (\ref{ffe1}) and the fact that $x\in B_{\beta^{\frac{1}{128}}}(x_0)$ and from Step 2 we know 
$x_0\in N_{c_4 \beta^{\frac{1}{512}}}\lt(\lt[\ai,\bi\rt]\rt)$, thus $d\in B_{c\beta^{\frac{1}{512}}}(0)$. Thus we have
\begin{eqnarray}
\label{eq198}
\ui(d)&=&\ui(d)-\ui(\ti{a})\nn\\
&\overset{(\ref{eq273}),(\ref{eq274}),(\ref{eq275})}{\geq}& \lt|d-\ti{a}\rt|-c\beta^{\frac{1}{256}}\nn\\
&\geq& \lt|\ti{a}\rt|-c\beta^{\frac{1}{512}}=2^{-1}\mathrm{diam}(\UI)-c\beta^{\frac{1}{512}}.
\end{eqnarray}

Pick $r_0\in \lt[\lt|d\rt|+\beta^{\frac{1}{512}},\lt|d\rt|+2\beta^{\frac{1}{512}}\rt]$ such that
$\int_{\partial B_{r_0}(0)} \lt|1-\lt|\na \ui(z)\rt|^2\rt| dH^1 z\leq c\beta^{-\frac{1}{512}}\beta$. Now fix
$y\in \partial B_{r_0}(0)$, let $s=\KI\cap \partial B_{r_0}\lt(0\rt)$ and
$\Gamma_1$ denote a connected component of
$\partial B_{r_0}(0)\backslash\lt\{s,y\rt\}$.
So we know $\int_{\Gamma_1\cup \lt[d,s\rt]} \lt|\na \ui(z)\rt| dH^1 z\leq
c H^1(\Gamma_1\cup \lt[d,s\rt])\leq c\beta^{\frac{1}{512}}$ so we can apply the fundamental theorem of Calculus
we have that $\lt|u(y)-u(d)\rt|\leq c\beta^{\frac{1}{512}}$ and since $y$ was an arbitrary point in
$\partial B_{r_0}(0)$, using (\ref{eq198}) this gives
\begin{equation}
\label{eq48}
\inf\lt\{\ui(z):z\in \partial B_{r_0}(0)\rt\}\geq 2^{-1}\mathrm{diam}(\UI)-c\beta^{\frac{1}{512}}.
\end{equation}
By definition (see (\ref{eq27})) $\ui(z)=u(z)+10^{-1}$ for any $z\in \partial B_{r_0}(0)$. Since $\mathrm{diam}(\UI)=\frac{22}{10}$ putting this with (\ref{eq48}) we have
(\ref{eq30}).  $\Box$\nl\nl

%
%
%
%
%

\em Proof of Theorem \ref{T1}. \rm Let $r_0\in(\CI_{3}^{-1}\beta^{\frac{1}{512}},\CI_{3}\beta^{\frac{1}{512}})$ be a number
we obtain from Lemma \ref{L5} that satisfies (\ref{eq30}). By Fubini's Theorem we know
$\int_{\Omega}\int_{\Omega} \lt|1-\lt|\na u(z)\rt|^2\rt|^2\lt|z-y\rt|^{-1} dz dy\leq \CI_7 \beta^2$ for some
constant $\CI_7>0$.
Let
\begin{equation}
\label{eqa38}
G_0:=\lt\{y\in\Omega:\int_{\Omega} \lt|1-\lt|\na u(z)\rt|^2\rt|^2 \lt|z-y\rt|^{-1} dz\leq\beta\rt\}.
\end{equation}
Note that $\lt|\Omega\backslash G_0\rt|\leq \CI_7\beta$.

Let $a,b\in \overline{\Omega}$ be such that $\lt|a-b\rt|=\mathrm{diam}(\Omega)$. Let $\vthe=\frac{a+b}{2}$. 
Since $r_0>\CI_{3}^{-1}\beta^{\frac{1}{512}}$ we can
pick $x_0\in B_{\beta^{\frac{1}{4}}}(\vthe)\cap G_0\subset B_{r_0}(\vthe)$. So by the Co-area formula there exists
$\Psi\subset S^1$ such that $H^1(S^1\backslash \Psi)\leq \sqrt{\beta}$ and
\begin{equation}
\label{eq57}
\int_{l_{x_0}^{\theta}\cap \Omega} \lt|1-\lt|\na u\rt|^2\rt|^2 dH^1 z\leq c\sqrt{\beta}
\text{ for each }\theta\in \Psi.
\end{equation}
For any $\theta\in S^1$ define $P(\theta):=l^{\theta}_{x_0}\cap \partial\Omega$, we will show
\begin{equation}
\label{eq50.5}
\lt|P(\theta)-x_0\rt|\geq 1-c\beta^{\frac{1}{512}}\text{ for any }\theta\in\Psi.
\end{equation}
To see this we argue as follows
\begin{eqnarray}
\label{eq58}
u(x_0)&=&u(x_0)-u\lt(P\lt(\theta\rt)\rt)\nn\\
&=&\int_{\lt[x_0,P\lt(\theta\rt)\rt]} \na u(z)\cdot (-\theta) dH^1 z\nn\\
&\overset{(\ref{eq57})}{\leq}&\lt|x_0-P\lt(\theta\rt)\rt|+c\beta^{\frac{1}{4}}.
\end{eqnarray}
Let $y_{\theta}:=\lt[x_0,P(\theta)\rt]\cap \partial B_{r_0}(\vthe)$. In exactly the same
way we have
\begin{equation}
\label{eqa21}
\lt|u(y_{\theta})-u(x_0)\rt|\leq c\beta^{\frac{1}{512}}.
\end{equation}
So
\begin{equation}
\label{eq50}
u(x_0)\geq u(y_{\theta})-\lt|u(y_{\theta})-u(x_0)\rt|\overset{(\ref{eqa21})}{\geq}
u(y_{\theta})-c\beta^{\frac{1}{512}}\overset{(\ref{eq30})}{\geq} 1-c\beta^{\frac{1}{512}}
\end{equation}
this together with (\ref{eq58}) establishes (\ref{eq50.5}).

Let $N=\lt[2^{-1}\beta^{-\frac{1}{2}}\rt]$, we
can divide $S^1$ into $N$ disjoint pieces of equal length, denote them $I_1,I_2,\dots I_N$.
Formally; $\bigcup_{k=1}^N I_k=S^1$ and $H^1(I_k)=\frac{2\pi}{N}$ for each $k=1,2,\dots N$. We can pick
$\theta_k\in I_k\cap\Psi$ for each $k=1,2,\dots N$.

Let
\begin{equation}
\label{eq51}
h=\min\lt\{\lt|P(\theta_k)-x_0\rt|:k\in\lt\{1,2,\dots N\rt\}\rt\}.
\end{equation}
We define $\Pi$ to be the convex hull of the points $x_0+h\theta_1, x_0+h\theta_2,\dots x_0+h\theta_N$.
Now by the construction of $\Pi$, for any $y\in\partial\Pi$ we can find $k\in\lt\{1,2,\dots N\rt\}$ such
that $\lt|y-(x_0+h\theta_k)\rt|\leq c\sqrt{\beta}$ and thus $\lt|y-x_0\rt|\geq h-c\sqrt{\beta}$ and so
\begin{equation}
\label{eqa31}
B_{h-c\sqrt{\beta}}(x_0)\subset \Pi.
\end{equation}
Note that by using (\ref{eq50.5}) we know 
$h>1-c\beta^{\frac{1}{512}}$ and since $\lt|x_0-\vthe\rt|\leq \beta^{\frac{1}{4}}$ (recalling also that $\Omega$ is convex and so $\Pi\subset \Omega$) there exists positive constant $\CI_8$ such that
\begin{equation}
\label{eqa30}
B_{1-\CI_8\beta^{\frac{1}{512}}}(\vthe)\subset \Omega.
\end{equation}
We claim
\begin{equation}
\label{eqa40}
\Omega\subset B_{1+2\CI_8\beta^{\frac{1}{512}}}(\vthe).
\end{equation}
Suppose not, so there exists $y\in\partial \Omega$ such that
$\lt|y-\vthe\rt|\geq  1+2\CI_8\beta^{\frac{1}{512}}$.
By inequality (\ref{eqa30}) we know $-\frac{y-\vthe}{\lt|y-\vthe\rt|}\lt(1-\CI_8\beta^{\frac{1}{512}}\rt)+\vthe\subset \Omega$
and as by convexity of $\Omega$, 
$\lt[y,\vthe-\frac{y-\vthe}{\lt|y-\vthe\rt|}\lt(1-\CI_8\beta^{\frac{1}{512}}\rt)\rt]\subset \Omega$ thus
$$
H^1\lt(\lt[y, \vthe-\frac{y-\vthe}{\lt|y-\vthe\rt|}\lt(1-\CI_8\beta^{\frac{1}{512}}\rt)\rt]\rt)
\geq 2+\CI_8\beta^{\frac{1}{512}}
$$
which contradicts the fact $\mathrm{diam}(\Omega)=2$ hence (\ref{eqa40}) is
established. Since the center of mass of $\Omega$ is $0$, i.e.\ $\int_{\Omega} dx=0$, by (\ref{eqa30}), (\ref{eqa40}) 
we have that $\lt|\vthe\rt|\leq c\beta^{\frac{1}{512}}$. 
Recall $x_0\in B_{\beta^{\frac{1}{4}}}(\vthe)$ so $\lt|x_0-P(\theta)\rt|\leq \lt|P(\theta)\rt|+\lt|x_0\rt|\overset{(\ref{eqa40})}{\leq}
1+c\beta^{\frac{1}{512}}$ so putting this together with (\ref{eq50}) we have
\begin{equation}
\label{eqa41}
u(x_0)-u(P(\theta))=u(x_0)\geq \lt|x_0-P(\theta)\rt|-c\beta^{\frac{1}{512}}.
\end{equation}

Thus
\begin{eqnarray}
\label{eqa42}
\int_{\lt[x_0,P(\theta)\rt]} \lt|\na u(z)+\theta\rt|^2 dH^1 z&=&\int_{\lt[x_0,P(\theta)\rt]}
\lt(\lt|\na u(z)\rt|^2+2\na u(z)\cdot\theta+1\rt) dH^1 z\nn\\
&\overset{(\ref{eq57})}{\leq}& 2(1+c\beta^{\frac{1}{4}})\lt|x_0-P(\theta)\rt|+2\lt(u(P(\theta))-u(x_0)\rt)\nn\\
&\overset{(\ref{eqa41})}{\leq}&c\beta^{\frac{1}{512}}\text{ for any } \theta\in\Psi.
\end{eqnarray}

Now using the elementary fact that $\lt|\na u(z)+\frac{z-x_0}{\lt|z-x_0\rt|}\rt|^2\leq
\lt|\lt|\na u(z)\rt|^2-1\rt|^2+4$, since $x_0\in G_0$ we have
\begin{eqnarray}
\label{eqc57}
&~&\int_{\theta\in S^1\backslash \Psi} \int_{l_{x_0}^{\theta}} \lt|\na u(z)+\frac{z-x_0}{\lt|z-x_0\rt|}\rt|^2 dH^1 z dH^1 \theta\nn\\
&~&\qd\qd\qd\qd\qd\qd\leq 4 H^1(S^1\backslash \Psi)+ \int_{\theta\in S^1}
\int_{l_{x_0}^{\theta}} \lt|\lt|\na u(z)\rt|^2-1\rt|^2 dH^1 z dH^1 \theta\nn\\
&~&\qd\qd\qd\qd\qd\qd \overset{(\ref{eqa38})}{\leq} 5\sqrt{\beta}.
\end{eqnarray}
And thus
\begin{eqnarray}
\int_{\Omega} \lt|\na u(z)+\frac{z-x_0}{\lt|z-x_0\rt|}\rt|^2 dz &\leq&c\int_{\Omega} \lt|\na u(z)+\frac{z-x_0}{\lt|z-x_0\rt|}\rt|^2 \lt|z-x_0\rt|^{-1}dz\nn\\
&\leq&c\int_{\theta\in S^1} \int_{l_{x_0}^{\theta}} \lt|\na u(z)+\frac{z-x_0}{\lt|z-x_0\rt|}\rt|^2 dH^1 z dH^1 \theta\nn\\
&\overset{(\ref{eqc57}),(\ref{eqa42})}{\leq}& c\beta^{\frac{1}{512}}.\nn
\end{eqnarray}
By Holder's inequality this gives
\begin{equation}
\label{eq60}
\lt(\int_{\Omega} \lt|\na u(z)+\frac{z-x_0}{\lt|z-x_0\rt|}\rt|^2 dz\rt)^{\frac{1}{2}}\leq c\beta^{\frac{1}{1024}}.
\end{equation}
Note that as $x_0\in B_{\beta^{\frac{1}{4}}}(\vthe)$ and (\ref{eqa31}), (\ref{eqa30}) we established that $\lt|\vthe\rt|\leq c\beta^{\frac{1}{512}}$  so 
$\lt|x_0\rt|\leq c\beta^{\frac{1}{512}}$. 
Now for any $z\in \Omega\backslash B_{\beta^{\frac{1}{1024}}}(0)$
\begin{eqnarray}
\label{ulaz100}
\lt|\frac{z}{\lt|z\rt|}-\frac{z-x_0}{\lt|z-x_0\rt|}\rt|&=&	
\lt|\frac{z\lt|z-x_0\rt|-(z-x_0)\lt|z\rt|}{\lt|z\rt|\lt|z-x_0\rt|}\rt|\nn\\
&=&\lt|\frac{z(\lt|z-x_0\rt|-\lt|z\rt|)+x_0\lt|z\rt|}{\lt|z\rt|\lt|z-x_0\rt|}\rt|\nn\\
&\leq&  \lt|\frac{\lt|z-x_0\rt|-\lt|z\rt|}{\lt|z-x_0\rt|}\rt|+\frac{\lt|x_0\rt|}{\lt|z-x_0\rt|}\nn\\
&\leq&c\beta^{\frac{1}{1024}}.
\end{eqnarray}
So
$$
\lt(\int_{\Omega} \lt|\frac{z}{\lt|z\rt|}-\frac{z-x_0}{\lt|z-x_0\rt|}\rt|^2 dz\rt)^{\frac{1}{2}}\leq
c\beta^{\frac{1}{512}}+\lt(\int_{\Omega\backslash B_{\beta^{\frac{1}{1024}}}} \lt|\frac{z}{\lt|z\rt|}-\frac{z-x_0}{\lt|z-x_0\rt|}\rt|^2 dz\rt)^{\frac{1}{2}}\overset{(\ref{ulaz100})}{\leq} c \beta^{\frac{1}{1024}}.
$$
Putting this together with (\ref{eq60}) we have (\ref{eq70.6}). $\Box$

\section{Proof of Corollary \ref{CC3}}
\label{sec4}

We begin by establishing the following proposition.

\begin{a5}
\label{PP7}

Let $\Omega$ be a bounded convex domain with $C^2$ boundary 
and $\lt|\Omega\triangle B_1(0)\rt|\leq \beta$ there exists a 
sequence $u^{\ep}\in C^{\infty}(\overline{\Omega})$ such that 
$u^{\ep}(z)=0$, 
$\na u^{\ep}(z)\cdot \eta_z=1$ for $z\in \partial \Omega$ (where $\eta_z$ is the inward 
pointing unit normal to $\partial \Omega$ at $z$) and for which
\begin{equation}
\label{ula201}
\limsup_{\ep\rightarrow 0}\int_{\Omega} \ep^{-1}\lt|1-\lt|\na u^{\ep}\rt|^2\rt|^2+\ep\lt|\na^2 u^{\ep}\rt|^2 dz\leq c\beta^{\frac{3}{32}}.
\end{equation}
\end{a5}

\subsubsection{Proof of Proposition \ref{PP7}}

\begin{a1}
\label{LL9}
 Suppose $\Omega$ is a convex and 
$\lt|\Omega\triangle B_1(0)\rt|=\beta$. Let $a_{\theta}=\partial\Omega \cap l^{\theta}_{0}$ we have
\begin{equation}
\label{abb22}
\lt|\lt|a_{\theta}\rt|-1\rt|\leq c\sqrt{\beta} \text{ and so }\partial \Omega\subset N_{c\sqrt{\beta}}(\partial B_1(0)).
\end{equation}
In addition there exists constant $c$ such that 
\begin{equation}
\label{abb23}
\lt|\eta_{a_{\theta}}+\theta\rt|\leq c\beta^{\frac{1}{4}}\text{ for any }\theta\in S^1.
\end{equation}
\end{a1}
\em Proof of Lemma. \rm

\em Step 1. \rm We will show $B_{\frac{1}{2}}(0)\subset \Omega$.

\em Proof of Step 1. \rm Suppose not, so we can pick $x\in\partial \Omega\cap B_{\frac{1}{2}}(0)$. Let $\eta_x$ be an inward pointing unit normal to $\partial \Omega$ at $x$, by convexity of
$\Omega$ we have $\Omega\subset \overline{H(x,\eta_x)}$ and so $B_1(0)\cap H(x,-\eta_x)\cap \Omega=\emptyset$ which implies $\lt|B_1(0)\backslash \Omega\rt|\geq \lt|B_1(0)\cap H(x,-\eta_x)\rt|>\frac{1}{8}$ which
contradicts that $\lt|\Omega\triangle B_1\rt|\leq \beta$.

\em Step 2. \rm $a_{\theta}\in B_{1+c\sqrt{\beta}}(0)$.

\em Proof of Step 2. \rm Suppose not. Since $\Omega$ is convex we have $\mathrm{conv}\lt(\lt\{a_{\theta}\rt\}\cup B_{\frac{1}{2}}(0)\rt)\subset \Omega$ and
$$
\lt|\mathrm{conv}\lt(\lt\{a_{\theta}\rt\}\cup B_{\frac{1}{2}}(0)\rt)\backslash B_1(0)\rt|>c\beta,
$$
thus we have $\lt|\Omega\backslash B_1(0)\rt|>c\beta$ which contradicts the fact that
$\lt|\Omega\triangle B_1(0)\rt|=\beta$.

\em Step 3. \rm We will show $a_{\theta}\not\in B_{1-c\sqrt{\beta}}(0)$.

\em Proof of Step 3. \rm Suppose $a_{\theta}\in B_{1-c\sqrt{\beta}}(0)$ this implies 
$\lt|B_1(0)\backslash H(a_{\theta},\eta_{a_{\theta}})\rt|\geq c\beta^{\frac{3}{4}}$ and 
$\Omega\subset H(a_{\theta},\eta_{a_{\theta}})$ so $\lt|B_1(0)\backslash \Omega\rt|\geq c\beta^{\frac{3}{4}}$ which gives a contradiction.

\em Proof of Lemma completed. \rm Suppose (\ref{abb23}) is false, since $\lt|a_{\theta}-\theta\rt|\leq c\sqrt{\beta}$ we have
$$
\lt|B_1(0)\backslash H(a_{\theta},\eta_{a_{\theta}})\rt|\geq c\beta^{\frac{3}{4}},
$$
as before this implies $\lt|B_1(0)\backslash \Omega\rt|>c\beta^{\frac{3}{4}}$ which is a 
contradiction. $\Box$

%
%

\begin{a1}
\label{LL10}
Let $\Omega$ be convex and define $u(x):=d(z,\partial \Omega)$ for any $z\in \Omega$ then
function $u$ is concave.
\end{a1}
\em Proof of Lemma. \rm Let $a,b\in \Omega$. Since $\Omega$ is convex
$\mathrm{conv}\lt(B_{u(a)}(a)\cup B_{u(b)}(b)\rt)\subset \Omega$. Now suppose there exists
$\lm\in (0,1)$ such that
$$
u\lt(\lm a+(1-\lm)b\rt)<\lm u(a)+(1-\lm)u\lt(b\rt)
$$
then as this implies $B_{u\lt(\lm a+\lt(1-\lm\rt) b\rt)}\lt(\lm a+\lt(1-\lm\rt) b\rt)\subset
\mathrm{int}
\lt(\mathrm{conv}\lt( B_{u(a)}(a)\cup B_{u(b)}(b)\rt)\rt)$ we must be able to find $x\in \partial \Omega$ with 
$x\in \partial \Omega\cap \mathrm{conv}\lt( B_{u(a)}(a)\cup B_{u(b)}(b)\rt)$ which is a contradiction. $\Box$

%
%

\begin{a1}
\label{LL25}
Let $\beta>0$, suppose $\Omega$ is a convex set with $\lt|\Omega\triangle B_1(0)\rt|\leq \beta$. 
Let $u(z)=d(z,\partial \Omega)$. For any $x\in \Omega\backslash B_{\beta^{\frac{1}{8}}}(0)$ for which the approximate
derivative $\na u$ exists
\begin{equation}
\label{abb72.6}
\lt|\na u(x)+\frac{x}{\lt|x\rt|}\rt|\leq c\beta^{\frac{3}{16}}.
\end{equation}
\end{a1}
\em Proof. \rm For any $x\in \Omega\backslash B_{\beta^{\frac{1}{8}}}(0)$ let
$b_x\in \partial \Omega$ be such that $\lt|b_x-x\rt|=u\lt(x\rt)$. 

We begin by showing 
\begin{equation}
\label{ula41}
\lt|b_x-\frac{x}{\lt|x\rt|}\rt|\leq c\beta^{\frac{3}{16}}.
\end{equation}

Recall from Lemma \ref{LL9} 
$a_{\frac{x}{\lt|x\rt|}}=\partial \Omega\cap l_0^{\frac{x}{\lt|x\rt|}}$.
Using (\ref{abb22}) from Lemma \ref{LL9} and the fact 
$(x-a_{\frac{x}{\lt|x\rt|}})\lt|x-a_{\frac{x}{\lt|x\rt|}}\rt|^{-1}=\frac{x}{\lt|x\rt|}$,
\begin{equation}
\label{ula84}
\lt|x-b_x\rt|\leq \lt|x-a_{\frac{x}{\lt|x\rt|}}\rt|\overset{(\ref{abb22})}{\leq} 1-\lt|x\rt|+c\sqrt{\beta}.
\end{equation}
Hence 
\begin{equation}
\label{ula165}
\lt|x-b_x\rt|^2=\lt|x\rt|^2-2x\cdot b_x+\lt|b_x\rt|^2\overset{(\ref{ula84})}{\leq} 1-2\lt|x\rt|+\lt|x\rt|^2+c\sqrt{\beta}.
\end{equation}
Therefor 
\begin{eqnarray}
-2x\cdot b_x&\overset{(\ref{ula165})}{\leq}& 1-2\lt|x\rt|+c\sqrt{\beta}-\lt|b_x\rt|^2\nn\\
&\overset{(\ref{abb22})}{\leq}&-2\lt|x\rt|+c\sqrt{\beta}.\nn
\end{eqnarray}
Thus $2\lt|x\rt|
\leq 2x\cdot b_x+c\sqrt{\beta}$. Since $\lt|x\rt|>\beta^{\frac{1}{8}}$ we have
\begin{equation}
\label{axa6}
1-c\beta^{\frac{3}{8}}\leq 1-c\frac{\beta^{\frac{1}{2}}}{\lt|x\rt|}\leq \frac{x}{\lt|x\rt|}\cdot b_x.
\end{equation}
Hence 
$$
\lt|b_x-\frac{x}{\lt|x\rt|}\rt|^2=\lt|b_x\rt|^2+1-2\frac{x}{\lt|x\rt|}\cdot b_x
\overset{(\ref{axa6}),(\ref{abb22})}{\leq} c\beta^{\frac{3}{8}}
$$
which gives 
\begin{equation}
\label{ax5}
\lt|\frac{x}{\lt|x\rt|}-b_x\rt|\leq c\beta^{\frac{3}{16}}.
\end{equation}
Let $\theta_x=\frac{b_x}{\lt|b_x\rt|}$ so using Lemma \ref{LL9} $\lt|\eta_{b_x}+\frac{b_x}{\lt|b_x\rt|}\rt|=
\lt|\eta_{a_{\theta_x}}+\theta_x\rt|\overset{(\ref{abb23})}{\leq} c\beta^{\frac{1}{4}}$ and by (\ref{abb22}) this easily implies 
\begin{equation}
\label{axa9}
\lt|\eta_{b_x}+b_x\rt|\leq c\beta^{\frac{1}{4}}. 
\end{equation}
Now since $\na u(x)=\frac{x-b_x}{\lt|x-b_x\rt|}=\eta_{b_x}$ and so 
$$
\lt|\na u(x)+\frac{x}{\lt|x\rt|}\rt|\leq \lt|\eta_{b_x}+b_x\rt|+\lt|\frac{x}{\lt|x\rt|}-b_x\rt|
\overset{(\ref{ax5}),(\ref{axa9})}{\leq}c\beta^{\frac{3}{16}}
$$
thus we have established (\ref{abb72.6}). $\Box$

%
%
%

%
%
%

\begin{a1}
\label{LL12}
Let $\Omega$ be a convex set and $\lt|\Omega\triangle B_1(0)\rt|\leq\beta$. 
Define $u(x)=d\lt(x,\partial \Omega\rt)$, 
note that since $u$ is convex $\na u$ is BV. Let $V(\na u,\cdot)$ denotes the total variation of the 
measure $\na u$. Firstly we have 
\begin{equation}
\label{ula301}
V(\na u,\Omega\backslash \overline{B_{3\beta^{\frac{1}{8}}}(0)})\leq 16\pi.
\end{equation}
For any $\veps\in (0,\beta^{\frac{1}{2}}]$, for any 
$x\in \Omega\backslash \lt(N_{2\veps}(\partial\Omega)\cup B_{4\beta^{\frac{1}{8}}}(0)\rt)$ we have
\begin{equation}
\label{abb40}
V(\na u, B_{\veps}(x))\leq c\beta^{\frac{3}{16}}\veps.
\end{equation}

\end{a1}
\em Proof. \rm Let $\tau\in (0,\frac{\veps}{20})$ be some small number. 
For any $x\in \Omega\backslash (\overline{N_{4\tau}(\partial \Omega)\cup B_{\frac{3}{2}\beta^{\frac{1}{8}}}(0)})=:\Pi_{\tau}$. 
Let $w_{\tau}(x)=u * \rho_{\tau}(x)$ and $v^{\tau}=\frac{\na w_{\tau}}{\lt|\na w_{\tau}\rt|}$. Note from Lemma \ref{LL25} for any $x\in \Pi_{\tau}$ 
\begin{eqnarray}
\label{axa19}
\lt|\na w_{\tau}(x)+\frac{x}{\lt|x\rt|}\rt|&=&\lt|\int \lt(\na u(x-z)+\frac{x}{\lt|x\rt|}\rt)
\rho_{\tau}(z) dz\rt|\nn\\
&\leq&\int \lt|\lt(\na u(x-z)+\frac{x-z}{\lt|x-z\rt|}\rt)\rho_{\tau}(z)\rt| dz
+\int \lt| \frac{x-z}{\lt|x-z\rt|}-\frac{x}{\lt|x\rt|}\rt| \rho_{\tau}(z) dz\nn\\
&\overset{(\ref{abb72.6})}{\leq}&c\sup_{z\in B_{2\tau}(0)}
\lt|\frac{x-z}{\lt|x-z\rt|}-\frac{x}{\lt|x\rt|}\rt|+c\beta^{\frac{3}{16}}\nn\\
&\leq&c\beta^{\frac{3}{16}}.
\end{eqnarray}                          

From this it is easy to conclude that 
\begin{equation}
\label{zaz1.8}
\|w_{\tau}-dist(\cdot, \partial B_1(0)))\|_{L^{\infty}(\Pi_{\tau})}\leq c\beta^{\frac{3}{16}}.
\end{equation}

\em Step 1. \rm Let $\tau_0>0$ be a very small number. We will show 
\begin{equation}
\label{fg2}
\lim_{\tau\rightarrow 0} \| v^{\tau}-\na u\|_{L^1(\Pi_{\tau_0})}=0.
\end{equation}

\em Proof of Step 1. \rm Now 
\begin{eqnarray}
\label{fg.05}
\int_{\Pi_{\tau_0}} \lt|1-\lt|\na w_{\tau}\rt|\rt| dz&=&\int_{\Pi_{\tau_0}}  \lt|\lt|\na u\rt|-\lt|\na w_{\tau}\rt|\rt|  dz\nn\\
&\leq& \int_{\Pi_{\tau_0}}  \lt|\na u-\na w_{\tau}\rt|  dz\rightarrow 0\text{ as }\tau\rightarrow 0.
\end{eqnarray}

Now from (\ref{axa19}) we have 
\begin{equation}
\label{fg1}
\lt|\na w_{\tau}(x)\rt|\geq \frac{1}{2}\text{ for any }x\in \Pi_{\tau_0}, \tau\in (0,\tau_0).
\end{equation}

So 
\begin{eqnarray}
\label{fg1.5}
\|\frac{\na w_{\tau}}{\lt|\na w_{\tau}\rt|}-\na w_{\tau}\|_{L^1(\Pi_{\tau_0})}&=&
\|\na w_{\tau}\lt(\frac{1}{\lt|\na w_{\tau}\rt|}-1\rt)\|_{L^1(\Pi_{\tau_0})}\nn\\
&\overset{(\ref{fg.05}),(\ref{fg1})}{\leq}&2\|1-\lt|\na w_{\tau}\rt|\|_{L^1(\Pi_{\tau_0})}\rightarrow 0\text{ as }\tau\rightarrow 0.
\end{eqnarray}
Since $\|\na w_{\tau}-\na u\|_{L^1(\Pi_{\tau_0})}\rightarrow 0$ as $\tau\rightarrow 0$ putting this together 
with (\ref{fg1.5}) gives (\ref{fg2}). \nl

\em Step 2. \rm We will show that for any $G\subset \subset \Omega\backslash \overline{B_{\frac{3}{2}\beta^{\frac{1}{8}}}(0)}$
\begin{equation}
\label{hhj1}
V(\na u,G)\leq 2\lt|\mathrm{div}(\na u)\rt|(G)
\end{equation}
and 
\begin{equation}
\label{ff6}
\lt|\mathrm{div}(\na u)\rt|(G)\leq \liminf_{\tau\rightarrow 0} \int_{G} \lt|v^{\tau}_{1,1}+v^{\tau}_{2,2}\rt| dz,
\end{equation}
where $\lt|\mathrm{div}(\na u)\rt|$ denotes the total variation of measure $\mathrm{div}(\na u)$.

\em Proof of Step 2. \rm We can find $\tau_0>0$ such that $G\subset \Pi_{\tau_0}$. 
Now from \cite{amdel} $\na u\in SBV_{loc}$ so in particular $\mathrm{div}(\na u)$ is a signed measure defined by 
\begin{equation}
\label{ff34}
\int \mathrm{div}(\na u) \phi dz=\int \phi_{,1}u_{,1}+\phi_{,2}u_{,2} dz\text{ for all }\phi\in C^{\infty}_{c}(\Omega).
\end{equation}

So for $\phi\in C^{\infty}_{c}(\Omega)$ we have 
\begin{eqnarray}
\int \mathrm{div}(\na u) \phi dz&\overset{(\ref{ff34}),(\ref{fg2})}{=}& \lim_{\tau\rightarrow 0} 
\int \phi_{,1}v^{\tau}_{1}+\phi_{,2}v^{\tau}_{2} dz\nn\\
&=&\lim_{\tau\rightarrow 0} \int ( v^{\tau}_{1,1}+v^{\tau}_{2,2})\phi dz. \nn
\end{eqnarray}
Now given open set $G\subset \Pi_{\tau_0}$ if $\phi\in C^{\infty}_{c}(G)$ then 
\begin{eqnarray}
\lt|\int \mathrm{div}(\na u)\phi dz\rt|&=& \lt|  \lim_{\tau\rightarrow 0} 
\int (v^{\tau}_{1,1}+v^{\tau}_{2,2})\phi dz\rt|\nn\\
&\leq& \|\phi\|_{L^{\infty}(\Pi_{\tau_0})}\int_{G} \lt|v^{\tau}_{1,1}+v^{\tau}_{2,2}\rt| dz.\nn
\end{eqnarray}
So this in particular by Proposition 1.47 \cite{amb} implies (\ref{ff6}).

Now since $\na u\in SBV_{loc}(\Omega)$ we know by Theorem 3.78 \cite{amb} that there exists a rectifiable set $J_{\na u}\subset S_{\na u}$ 
(where $S_{\na u}$ denotes the set of approximate jump points of $\na u$) with $H^{n-1}(S_{\na u}\backslash J_{\na u})=0$ and 
$D \na u\lfloor J_{\na u}=\lt(\na u^{+}-\na u^{-}\rt)\otimes \nu H^{n-1}\lfloor J_{\na u}$ where $\nu(x)$ is the normal to the 
approximate tangent of the rectifiable set $J_{\na u}$ at point $x$. Following \cite{amb} Definition 3.67 we assume 
that the triple $(\na u^{+}(x),\na u^{-}(x),\nu(x))$ satisfies (3.69) of \cite{amb}. By Theorem 3.94 \cite{amb} we have 
that $(\na u^{+}(x)-\na u^{-}(x))\otimes \nu(x)$ is a rank-$1$ matrix for $\lt|D \na u\rt|$ a.e.\ $x\in J_{\na u}$. Now 
$D \na u$ is a matrix valued measure and indeed letting $\partial_{i} u_{,j}$ denote the individual `component' measures, just from the 
definition we know that $\partial_{i} u_{,j}=\partial_j u_{,i}$ so $D \na u$ is a symmetric matrix valued measure. Specifically 
by differentiation of measures (see Theorem 2.2 \cite{amb}) $M(x):=\lim_{r\rightarrow 0} \frac{D \na u(B_r(x))}{\lt|D \na u\rt|(B_r(x))}$ 
exists for $\lt|D \na u\rt|$ a.e.\ $x$ and $M(x)$ will be a symmetric $2\times 2$ matrix. So for $H^{n-1}$ a.e.\ $x\in J_{\na u}$, 
$(\na u^{+}(x)-\na u^{-}(x))\otimes \nu(x)$ is a symmetric rank-$1$ matrix, this is easily seen to imply 
$\frac{\na u^{+}(x)-\na u^{-}(x)}{\lt|\na u^{+}(x)-\na u^{-}(x)\rt|}=\nu(x)$. So 
$(\na u^{+}(x)-\na u^{-}(x))\otimes \nu(x)=\lt|\na u^{+}(x)-\na u^{-}(x)\rt| \nu(x)\otimes \nu(x)$. Thus we can decompose $D(\na u)$ into 
absolutely continuous and singular parts we have  
\begin{equation}
\label{ff51}
D(\na u)(S)=\int_{S} D(\na u) dx+\int_{S\cap J_{\na u}} \lt|\na u^{+}-\na u^{-}\rt| \nu(x)\otimes \nu(x) dH^{1}\text{ for any set }S\subset \R^n.
\end{equation}
 
Obviously this is a matrix valued Radon measure and the signed Radon measure $\Delta u$ is given by the sum of diagonal elements of the 
matrix defined by (\ref{ff51}) and so is given by 
\begin{eqnarray}
\Delta u(S)&=& \int_{S} \mathrm{div}_{a}(\na u) dx+\int_{S\cap J_{\na u}} \lt|\na u^{+}-\na u^{-}\rt| \nu\cdot \nu dH^1 \nn\\
&=&\int_{S} \mathrm{div}_{a}(\na u) dx+\int_{S\cap J_{\na u}} \lt|\na u^{+}-\na u^{-}\rt|  dH^1\text{ for any }S\subset \R^n.\nn
\end{eqnarray}
Now recall $\lt|\na u(x)\rt|=1$ for a.e.\ $x\in \Omega$. So by Volpert chain rule (see Theorem 3.94 \cite{amb}) 
we have that the function $x\rightarrow \lt|\na u(x)\rt|^2$ is BV and the standard chain rule holds so 
\begin{eqnarray}
\label{fgg1.2}
&~&u_{,11}(x)u_{,1}(x)+u_{,12}(x)u_{,2}(x)=0\text{ and }\nn\\
&~&u_{,12}(x)u_{,1}(x)+u_{,22}(x)u_{,2}(x)=0\text{ for }a.e.\ x\in \Omega.
\end{eqnarray}
Since $u_{,21}=u_{,12}$ we have
$$
\lt(\begin{matrix} u_{,11} &    u_{,12}\nn\\
u_{,21} &    u_{,22} \end{matrix}\rt)\lt(\begin{matrix} u_{,1} \nn\\
u_{,2} \end{matrix}\rt)\overset{(\ref{fgg1.2})}{=}\lt(\begin{matrix} 0 \nn\\
 0 \end{matrix}\rt)\text{ and }\lt(\begin{matrix} u_{,11} &    u_{,12}\nn\\
u_{,21} &    u_{,22} \end{matrix}\rt)\lt(\begin{matrix} -u_{,2} \nn\\
u_{,1} \end{matrix}\rt)\overset{(\ref{fgg1.2})}{=}
(u_{,11}+ u_{,22})
\lt(\begin{matrix} -u_{,2} \nn\\
u_{,1} \end{matrix}\rt).\nn
$$
Letting $\|\cdot\|$ denote the operator norm of a matrix, since
$\lt(\begin{matrix} u_{,1} &    -u_{,2} \nn\\
u_{,2} & u_{,1}\end{matrix}\rt)\in O(2)$ thus
\begin{eqnarray}
\label{abb55}
\lt|\lt|\lt(\begin{matrix} u_{,11} &    u_{,12}\nn\\
u_{,21} &    u_{,22} \end{matrix}\rt)\rt|\rt|&=&
\lt|\lt| \lt(\begin{matrix} u_{,11} &    u_{,12}\nn\\
u_{,21} &    u_{,22} \end{matrix}\rt)
\lt(\begin{matrix} u_{,1} & -u_{,2} \nn\\
u_{,2} & u_{,1}\end{matrix}\rt)\rt|\rt|\nn\\
&=&\lt|\lt| \lt(\begin{matrix} 0 &    -(u_{,11}+u_{,22})u_{,2}\nn\\
0 &   (u_{,11}+u_{,22})u_{,1} \end{matrix}\rt)\rt|\rt|\nn\\
&\leq&2\lt|u_{,11}+u_{,22}\rt|.\nn
\end{eqnarray}
So 
$$
\lt|D_a(\na u(x))\rt|\leq 2\lt|\mathrm{div}_a (\na u(x))\rt|\text{ for }a.e.\ x\in \Omega. 
$$
Thus 
\begin{eqnarray}
\label{ff12}
V(\na u,G)&=&\int_{G} \lt|D_a (\na u)\rt| dz+\int_{G\cap J_{\na u}} \lt|\na u^{+}-\na u^{-}\rt| dH^1 \nn\\
&\leq& 2\int_{G} \lt|\mathrm{div}_a (\na u)\rt| dz+\int_{G\cap J_{\na u}} \lt|\na u^{+}-\na u^{-}\rt| dH^1 \nn\\
&\leq& 2\lt|\mathrm{div}(\na u)\rt|(G),\nn
\end{eqnarray}
thus establishing (\ref{hhj1}).\nl

\em Step 3. \rm We will show that for any $t\in (8\tau,1-2\beta^{\frac{1}{8}})$
\begin{equation}
\label{abc4}
\int_{w^{-1}_{\tau}(t)}
\lt|v^{\tau}_{1,1}(z)+v^{\tau}_{2,2}(z)\rt| dH^1 z\leq 2\pi.
\end{equation}

\em Proof of Step 3. \rm We define the `angle' function by
\begin{equation}
\label{acb10}
A(x):=\lt\{\begin{array}{ll} \arccos\lt(\frac{x_1}{\lt|x\rt|}\rt)&\text{ for } x_2\geq 0\\
2\pi- \arccos\lt(\frac{x_1}{\lt|x\rt|}\rt)&\text{ for } x_2< 0\end{array}\rt.
\end{equation}
Note that $A$ is smooth expect at the half line $\lt\{(x_1,x_2):x_2=0, x_1>0\rt\}$. For $x\in \Pi_{\tau}$ we 
have $\lt|v^{\tau}(x)\rt|^2=1$, so as before
\begin{equation}
\label{abbb24}
\partial_1(\lt|v^{\tau}(x)\rt|^2)=v^{\tau}_1(x)v^{\tau}_{1,1}(x)+v^{\tau}_2(x)v^{\tau}_{2,1}(x)=0.
\end{equation}

Since $u$ is the $1$-Lipschitz, 
\begin{equation}
\label{ulz3}
\|w_{\tau}-u\|_{L^{\infty}(\Pi_{\tau})}\leq 2\tau,
\end{equation}
and so from this and (\ref{zaz1.8}) we have that for any $t\in (8\tau,1-2\beta^{\frac{1}{8}})$, 
$w_{\tau}^{-1}(t)\subset \Pi_{\tau}$ and hence by (\ref{axa19}) $v^{\tau}$ 
is well defined along this level set. 
We also know that for any $x\in w^{-1}_{\tau}(t)$ the tangent to curve
$w^{-1}_{\tau}(t)$ is given by
$\lt(\begin{matrix} -v_2^{\tau}(x)\nn\\ v_1^{\tau}(x) \end{matrix}\rt)$. 
Note that $w_{\tau}^{-1}(t)$ is the boundary of a smooth convex set so there 
exists a point $x_t\in w_{\tau}^{-1}(t)$ such that 
$A\lt(\begin{matrix} -v_2^{\tau}(x_t)\nn\\ v_1^{\tau}(x_t) \end{matrix}\rt)=0$. There 
must also exist $y_t\in w_{\tau}^{-1}(t)$ such that 
\begin{equation}
\label{kband01}
A\lt(\begin{matrix} -v_2^{\tau}(y_t)\\ v_1^{\tau}(y_t) \end{matrix}\rt)=\pi.
\end{equation}
Let $\Phi^t:\lt[0,H^1(w_{\tau}^{-1}(t))\rt)\rightarrow w_{\tau}^{-1}(t)$ denote the clockwise parameterization of $w_{\tau}^{-1}(t)$ by arc-length with $\Phi^t(0)=x_t$. So $\dot{\Phi^t}(s)=\lt(\begin{matrix} -v_2^{\tau}(\Phi^t(s))\nn\\ v_1^{\tau}(\Phi^t(s)) \end{matrix}\rt)$. Define
$\Theta_t:[0,H^1(w_{\tau}^{-1}(t)))\rightarrow \R$ by $\Theta_t(s)=A(\dot{\Phi}^t(s))$. Now pick
$s\in\lt(0,H^1(w_{\tau}^{-1}(t))\rt)$, suppose $v_1^{\tau}\lt(\Phi^t(s)\rt)> 0$, then
\begin{eqnarray}
\label{axx14}
\dot{\Theta}_t(s)&=&\dot{\arccos}\lt(-v_2^{\tau}\lt(\Phi^t(s)\rt)\rt)\frac{\partial}{\partial t}
\lt( -v_2^{\tau}\lt(\Phi^t(s)\rt)\rt)\nn\\
&=&\dot{\arccos}\lt(-v_2^{\tau}\lt(\Phi^t(s)\rt)\rt)
\lt( -v_{2,1}^{\tau}\lt(\Phi^t(s)\rt)\dot{\Phi}^t_1(t)-v_{2,2}^{\tau}\lt(\Phi^t(s)\rt)\dot{\Phi}^t_2(t)\rt)\nn\\
&=&\dot{\arccos}\lt(-v_2^{\tau}\lt(\Phi^t(s)\rt)\rt)
\lt(v_{2,1}^{\tau}\lt(\Phi^t(s)\rt)v_2^{\tau}\lt(\Phi^t(s)\rt)  -v_{2,2}^{\tau}\lt(\Phi^t(s)\rt)v_1^{\tau}\lt(\Phi^t(s)\rt)  \rt)\nn\\
&\overset{(\ref{abbb24})}{=}&\dot{\arccos}\lt(-v_2^{\tau}\lt(\Phi^t(s)\rt)\rt)
\lt( -v_{1,1}^{\tau}\lt(\Phi^t(s)\rt)v_1^{\tau}\lt(\Phi^t(s)\rt)  -v_{2,2}^{\tau}\lt(\Phi^t(s)\rt)v_1^{\tau}\lt(\Phi^t(s)\rt)  \rt)\nn\\
&=&-\dot{\arccos}\lt(-v_2^{\tau}\lt(\Phi^t(s)\rt)\rt)v_1^{\tau}\lt(\Phi^t(s)\rt)
\lt(v_{1,1}^{\tau}\lt(\Phi^t(s)\rt)+ v_{2,2}^{\tau}\lt(\Phi^t(s)\rt)\rt).
\end{eqnarray}
Now for any $w\in (-1,1)$, $\dot{\arccos}(w)=-(\sin(\arccos(w)))^{-1}$ so
\begin{equation}
\label{axx12}
\dot{\Theta}_t(t)=\frac{v_1^{\tau}\lt(\Phi^t(s)\rt)}{\sin(\arccos(-v_2^{\tau}\lt(\Phi^t(s)\rt)))}
\lt(v_{1,1}^{\tau}\lt(\Phi^t(s)\rt)+ v_{2,2}^{\tau}\lt(\Phi^t(s)\rt)\rt).
\end{equation}
Recall $\lt|\lt(\begin{matrix} -v_2^{\tau}\lt(\Phi^t(s)\rt)\nn\\ v_1^{\tau}\lt(\Phi^t(s)\rt) \end{matrix}\rt)\rt|=1$ and we supposed $v_1^{\tau}\lt(\Phi^t(s)\rt)> 0$, so
\begin{eqnarray}
\label{axa12}
v_1^{\tau}\lt(\Phi^t(s)\rt)&=&\sqrt{1-\lt(v_2^{\tau}\lt(\Phi^t(s)\rt)\rt)^2}\nn\\
&=&\sqrt{1-\lt(\cos\lt(\arccos\lt(-v_2^{\tau}\lt(\Phi^t(s)\rt)\rt)\rt)\rt)^2}\nn\\
&=&\sin\lt(\arccos\lt(-v_2^{\tau}\lt(\Phi^t(s)\rt)\rt)\rt).
\end{eqnarray}

Thus from (\ref{axx12})
\begin{equation}
\label{abb27}
\dot{\Theta}_t(s)=\lt(v_{1,1}^{\tau}\lt(\Phi^t(s)\rt)+ v_{2,2}^{\tau}\lt(\Phi^t(s)\rt)\rt)\text{ for any }
s\in\lt(0,H^1(w_{\tau}^{-1}(t))\rt)\text{ with }v_1^{\tau}\lt(\Phi^t(s)\rt)> 0.
\end{equation}

Suppose we have $s\in\lt(0,H^1(w_{\tau}^{-1}(t))\rt)$ with $v_{1}^{\tau}\lt(\Phi^t(s)\rt)<0$, then 
in the same way as (\ref{axa12}) we have 
\begin{equation}
\label{axa13}
v_1^{\tau}\lt(\Phi^t(s)\rt)=-\sqrt{1-\lt(\cos\lt(\arccos\lt(-v_2^{\tau}\lt(\Phi^t(s)\rt)\rt)\rt)\rt)^2}
=-\sin\lt(\arccos\lt(-v_2^{\tau}\lt(\Phi^t(s)\rt)\rt)\rt).
\end{equation}
And since $v_1^{\tau}(\Phi^t(s))<0$, by definition of $A$ (see (\ref{acb10})) arguing as in 
(\ref{axx12}) we have
\begin{eqnarray}
\label{abb26}
\dot{\Theta}_t(s)&=&\frac{-v_1^{\tau}\lt(\Phi^t(s)\rt)}{\sin\lt(\arccos\lt(-v_2^{\tau}\lt(\Phi^t(s)\rt)\rt)\rt)}
\lt(v_{1,1}^{\tau}\lt(\Phi^t(s)\rt)+ v_{2,2}^{\tau}\lt(\Phi^t(s)\rt)\rt)\nn\\
&\overset{(\ref{axa13})}{=}&v_{1,1}^{\tau}\lt(\Phi^t(s)\rt)+ v_{2,2}^{\tau}\lt(\Phi^t(s)\rt)\text{ for }
s\in\lt(0,H^1(w_{\tau}^{-1}(t))\rt)\text{ with }v_1^{\tau}\lt(\Phi^t(s)\rt)< 0.\nn
\end{eqnarray}
Without loss of generality we can assume 
$\lt|\lt\{s\in [0,H^1(w_{\tau}^{-1}(t))]: v_1^{\tau}\lt(\Phi^t(s)\rt)=0\rt\}\rt|=0$. 
Thus by continuity of $\dot{\Theta}_t(\cdot)$, $v_{1,1}^{\tau}(\Phi^{t}(\cdot))$ and 
$v_{2,2}^{\tau}(\Phi^{t}(\cdot))$ we have
\begin{equation}
\label{abb65}
\dot{\Theta}_t(s)=v_{1,1}^{\tau}\lt(\Phi^t(s)\rt)+v_{2,2}^{\tau}\lt(\Phi^t(s)\rt)\text{ for }
s\in\lt[0,H^1(w^{-1}_{\tau}(t))\rt).
\end{equation}
Now since $u$ is concave, $w_{\tau}$ is concave and so the set
$w^{-1}_{\tau}(\lt[t,\infty\rt))$ is a convex set,
hence
\begin{equation}
\label{abb60}
v_{1,1}^{\tau}\lt(\Phi^t(s)\rt)+ v_{2,2}^{\tau}\lt(\Phi^t(s)\rt)=
\dot{\Theta}_t(s)\geq 0\text{ for any }s\in\lt[0,H^1(w^{-1}_{\tau}(t))\rt).
\end{equation}
Hence
\begin{equation}
\label{lband0205}
\int_{w_{\tau}^{-1}(t)} \lt|v_{1,1}^{\tau}(z)+v_{2,2}^{\tau}(z)\rt| dH^1 z=\int_{0}^{H^1(w_{\tau}^{-1}(t))}
\dot{\Theta}_t(s) ds\leq 2\pi.
\end{equation}

\em Step 4. \rm Let $x\in \Pi_{\tau}\backslash N_{2\veps}(\partial \Omega)$ and define 
\begin{equation}
\label{axa15}
t_1=\inf\lt\{s\in \R:w^{-1}_{\tau}(s)\cap B_{\veps}(x)\not=\emptyset\rt\}\text{ and }
t_2=\sup\lt\{s\in \R:w^{-1}_{\tau}(s)\cap B_{\veps}(x)\not=\emptyset\rt\}.
\end{equation}

Recall $y_t\in w_{\tau}^{-1}(t)$ was chosen so that 
(\ref{kband01}) holds true, let 
$\pi_t:=(\Phi^t)^{-1}(y_t)$. We have for any $t\in\lt(t_1,t_2\rt)$
\begin{equation}
\label{abb99}
\sup\lt\{\lt|\Theta_t(s_1)-\Theta_t(s_2)\rt|:s_1, s_2\in (\Phi^t)^{-1}\lt(w_{\tau}^{-1}(t)\cap B_{\veps}(x)\rt)\cap \lt[0,\pi_t\rt]\rt\}\leq c\beta^{\frac{3}{16}}
\end{equation}
and 
\begin{equation}
\label{abb99.6}
\sup\lt\{\lt|\Theta_t(s_1)-\Theta_t(s_2)\rt|:s_1, s_2\in (\Phi^t)^{-1}\lt(w_{\tau}^{-1}(t)\cap B_{\veps}(x)\rt)\cap \lt[\pi_t, H^1(w^{-1}_{\tau}(t))  \rt)\rt\}\leq c\beta^{\frac{3}{16}}.
\end{equation}

\em Proof of Step 4. \rm Let $s_1,s_2\in \lt[0,\pi_t\rt]$ such that 
$\Phi^t(s_1),\Phi^t(s_2)\in B_{\veps}(x)$, since $\Phi^t$ is parameterization of $w_{\tau}^{-1}(t)$
by arclength $\dot{\Phi}^t(s)$ is the unit tangent to $w^{-1}_{\tau}(t)$ at $\Phi^t(s)$. Thus
$$
R\lt(\frac{\na w_{\tau}\lt(\Phi^t(s_i)\rt)}{\lt|\na w_{\tau}\lt(\Phi^t(s_i)\rt)\rt|}\rt)
=\dot{\Phi}^t(s_i)\text{ for }i=1,2.
$$

However by Lemma \ref{LL25} (recalling the fact that $\lt|\Phi^t(s_1)\rt|>\frac{3\beta^{\frac{1}{8}}}{2}$ and  
$\lt|\Phi^t(s_2)\rt|>\frac{3\beta^{\frac{1}{8}}}{2}$ in order to apply the lemma)
\begin{eqnarray}
\label{ula51.5}
\lt|\na w_{\tau}\lt(\Phi^t(s_1)\rt)-\na w_{\tau}\lt(\Phi^t(s_2)\rt)\rt|&=&
\lt|\int \lt(\na u\lt(\Phi^t\lt(s_1\rt)-z\rt)- \na u\lt(\Phi^t\lt(s_2\rt)-z\rt)\rt)\rho_{\tau}(z) dz
\rt|\nn\\
&\overset{(\ref{abb72.6})}{\leq}&c\int_{B_{\tau}(0)} \lt|\frac{\Phi^t(s_1)-z}{\lt|\Phi^t(s_1)-z\rt|}-
\frac{\Phi^t(s_2)-z}{\lt|\Phi^t(s_2)-z\rt|}\rt|\rho_{\tau}(z) dz+c\beta^{\frac{3}{16}}.\nn\\
\end{eqnarray}
Note $z\in B_{\tau}(0)\subset B_{\frac{\beta^{\frac{1}{2}}}{20}}(0)$ so as $\lt|\Phi^t(s_1)\rt|>\frac{3\beta^{\frac{1}{8}}}{2}$ we have $\lt|\Phi^t(s_1)-z\rt|\geq \lt|\Phi^t(s_1)\rt|-\lt|z\rt|\geq \beta^{\frac{1}{8}}.$
Recall the elementary inequality inequality 
\begin{equation}
\label{df1}
\lt|\frac{z}{\lt|z\rt|}-\frac{y}{\lt|y\rt|}\rt|\leq 2\lt|z-y\rt|\text{ for any }z,y\text{ with }\lt|z\rt|\geq 1, \lt|y\rt|\geq 1.
\end{equation}
So in particular we have  
\begin{equation}
\label{ewq1}
\lt|\lt|\frac{\Phi^t(s_1)-z}{\lt|\Phi^t(s_1)-z\rt|}-
\frac{\Phi^t(s_2)-z}{\lt|\Phi^t(s_2)-z\rt|}\rt|\rt|
\leq \frac{2}{\beta^{\frac{1}{8}}}\lt|\Phi^t(s_1)-\Phi^t(s_2)\rt|\leq 2\beta^{\frac{3}{8}}.
\end{equation}
Thus with (\ref{ula51.5}) this gives
\begin{equation}
\label{ula51}
\lt|\na w_{\tau}\lt(\Phi^t(s_1)\rt)-\na w_{\tau}\lt(\Phi^t(s_2)\rt)\rt|\leq c\beta^{\frac{3}{16}}.
\end{equation}
As a consequence of (\ref{axa19}) we know 
\begin{equation}
\label{ula91}
\lt|\lt|\na w_{\tau}(x)\rt|-1\rt|\leq c\beta^{\frac{3}{16}}\text{ for any }x\in \Pi_{\tau}
\end{equation}
so
\begin{eqnarray}
\label{ula53}
\lt|\dot{\Phi}(s_1)-\dot{\Phi}(s_2)\rt|&\overset{(\ref{ula51})}{\leq}&
\lt|R\lt( \frac{\na w_{\tau}\lt(\Phi^t(s_1)\rt)}{\lt|\na w_{\tau}\lt(\Phi^t(s_1)\rt)\rt|}\rt)
-R\lt(\na w_{\tau}\lt(\Phi^t(s_1)\rt)\rt)\rt|\nn\\
&~&
+\lt|R\lt( \frac{\na w_{\tau}\lt(\Phi^t(s_2)\rt)}{\lt|\na w_{\tau}\lt(\Phi^t(s_2)\rt)\rt|}\rt)
-R\lt(\na w_{\tau}\lt(\Phi^t(s_2)\rt)\rt)\rt|+c\beta^{\frac{3}{16}}\nn\\
&\overset{(\ref{ula91})}{\leq}&c\beta^{\frac{3}{16}}.
\end{eqnarray}
Now as $s_1,s_2\in \lt[0,\pi_t\rt]$, since $w_{\tau}^{-1}(t)$ is the boundary of a convex set so we know  
$\dot{\Phi}^t(s_1), \dot{\Phi}^t(s_2)\in 
\lt\{v\in S^1: v\cdot e_2\leq 0\rt\}$. Now as $A$ is Lipschitz on 
$\lt\{v\in S^1: v\cdot e_2\leq 0\rt\}$,
\begin{equation}
\label{ula175}
\lt|\Theta_t(s_1)-\Theta_t(s_2)\rt|=\lt|A\lt(\dot{\Phi}^t(s_1)\rt)-A\lt(\dot{\Phi}^t(s_2)\rt)\rt|
\overset{(\ref{ula53})}{\leq} c\beta^{\frac{3}{16}}
\end{equation}
and so (\ref{abb99}) is established. Inequality (\ref{abb99.6}) follows in exactly 
the same way. \nl

\em Step 5. \rm We will show 
\begin{equation}
\label{ffg14}
V(\na u,B_{\veps}(x))\leq c\veps \beta^{\frac{3}{16}}\text{ for all }x\in \Omega\backslash \lt(N_{2\veps}(\partial \Omega)\cup 
B_{4\beta^{\frac{1}{4}}}(0)\rt)
\end{equation}

\em Proof of Step 5. \rm Let $x\in \Omega\backslash \lt(N_{2\veps}(\partial \Omega)\cup 
B_{4\beta^{\frac{1}{4}}}(0)\rt)$. Let $t\in (t_1,t_2)$. The 
most non-trivial case is where 
$$
\lt\{s\in \lt[0,\pi_t\rt]:\Phi^t(s)\in B_{\veps}(x)\rt\}\not=\emptyset\text{ and }
\lt\{s\in \lt[\pi_t,H^1(w_{\tau}^{-1}(t))\rt]:\Phi^t(s)\in B_{\veps}(x)\rt\}\not=\emptyset.
$$
When either of these sets is empty the proof follow in a very similar way.

Let $s_1^t=\inf\lt\{s\in \lt[0,\pi_t\rt]:\Phi^t(s)\in B_{\veps}(x)\rt\}$,
$s_2^t=\sup\lt\{s\in \lt[0,\pi_t\rt]:\Phi^t(s)\in B_{\veps}(x)\rt\}$. So
$\lt[s^t_1,s^t_2\rt]=\lt\{s\in \lt[0,\pi_t\rt]:\Phi^t(s)\in B_{\veps}(x)\rt\}$. Now
\begin{eqnarray}
\label{axz40}
\int_{\lt[s^t_1,s^t_2\rt]} \lt|v^{\tau}_{1,1}(\Phi^t(s))+ v^{\tau}_{2,2}(\Phi_t(s))\rt| ds
&\overset{(\ref{abb60})}{=}&\int_{\lt[s^t_1,s^t_2\rt]} \dot{\Theta}_t(s) ds\nn\\
&\overset{(\ref{abb99})}{\leq}&c\beta^{\frac{3}{16}}.
\end{eqnarray}
In the same way of we let 
$$
r_1^t=\inf\lt\{s\in \lt[\pi_t,H^1(w_{\tau}^{-1}(t))\rt]:\Phi^t(s)\in B_{\veps}(x)\rt\}, 
r_2^t=\sup\lt\{s\in \lt[\pi_t,H^1(w_{\tau}^{-1}(t))\rt]:\Phi^t(s)\in B_{\veps}(x)\rt\}
$$
then 
\begin{equation}
\label{kband02}
\int_{\lt[r^t_1,r^t_2\rt]} \lt|v^{\tau}_{1,1}(\Phi^t(s))+ v^{\tau}_{2,2}(\dot{\Phi}_t(s))\rt| ds\leq c\beta^{\frac{3}{16}}.
\end{equation}
Thus
\begin{eqnarray}
&~&
\int_{B_{\veps}(x)}  \lt|v^{\tau}_{1,1}(z)+ v^{\tau}_{2,2}(z)\rt| \lt|\na w_{\tau}(z)\rt| dz\nn\\
&~&\qd\qd\qd=\int_{t_1}^{t_2} \int_{w_{\tau}^{-1}(t)} \lt|v^{\tau}_{1,1}(z)+ v^{\tau}_{2,2}(z)\rt| dH^1 z dt\nn\\
&~&\qd\qd\qd=\int_{t_1}^{t_2} \int_{\lt[s^t_1,s^t_2\rt]\cup \lt[r^t_1,r^t_2\rt]} \lt|v^{\tau}_{1,1}(\Phi^t(s))+ v^{\tau}_{1,1}(\Phi^t(s))\rt| ds dt\nn\\
&~&\qd\qd\qd\overset{(\ref{axz40}),(\ref{kband02})}{\leq} c\lt|t_1-t_2\rt|\beta^{\frac{3}{16}}.\nn
\end{eqnarray}
By using (\ref{axa19}) and recalling the definition (\ref{axa15}) of Step 2 we must have $\lt|t_1-t_2\rt|\leq c\veps$. Also 
from (\ref{axa19}) we 
know $\lt|\na w_{\tau}(z)\rt|\geq 1-c\beta^{\frac{3}{16}}$ for all $z\in B_{\veps}(x)$, so putting these things together
we have
\begin{equation}
\label{acb1}
\int_{B_{\veps}(x)}  \lt|v^{\tau}_{1,1}(z)+ v^{\tau}_{2,2}(z)\rt| dz\leq c\veps\beta^{\frac{3}{16}}\text{ for all }
x\in \Omega\backslash \lt(N_{2\veps}(\partial \Omega)\cup B_{4\beta^{\frac{1}{8}}}(0)\rt)
\end{equation}
So for any $x\in\Omega\backslash \lt(N_{2\veps}(\partial \Omega)\cup B_{4\beta^{\frac{1}{8}}}(0)\rt)$ we know 
$B_{\veps}(x)\subset \Pi_{\frac{\veps}{4}}$ so by Step 2  
\begin{eqnarray}
V(\na u,B_{\veps}(x))&\overset{(\ref{hhj1})}{\leq}&2\lt|\mathrm{div}(\na u)\rt|(B_{\veps}(x))\nn\\
&\overset{(\ref{ff6})}{\leq}&2\liminf_{\tau\rightarrow 0} \int_{B_{\veps}(x)} \lt|v_{1,1}^{\tau}+v_{2,2}^{\tau}\rt| dz\nn\\
&\leq&c\veps \beta^{\frac{3}{16}},\nn
\end{eqnarray}
and so we have established (\ref{ffg14}). 

\em Proof of Lemma completed. \rm Note that by (\ref{zaz1.8}) and (\ref{ulz3}) we have 
$$
\Pi_{16\tau}\backslash \overline{B_{3\beta^{\frac{1}{8}}}(0)}\subset w_{\tau}^{-1}\lt(\lt[8\tau,1-2\beta^{\frac{1}{8}}\rt]\rt)
$$ 
by using the Co-area formula 
$$
\int_{\Pi_{16\tau}\backslash \overline{B_{3\beta^{\frac{1}{8}}}(0)}} \lt|v_{1,1}^{\tau}+v_{2,2}^{\tau}\rt|\lt|\na w^{\tau}\rt| dz\leq 
\int_{8\tau}^{1-2\beta^{\frac{1}{8}}} \int_{w_{\tau}^{-1}(s)} \lt|v_{1,1}^{\tau}+v_{2,2}^{\tau}\rt| dH^1 z ds\leq 4\pi.
$$
Thus using (\ref{axa19}) 
\begin{equation}
\label{ula280}
\int_{\Pi_{16\tau}\backslash  \overline{B_{3\beta^{\frac{1}{8}}}(0)}} \lt|v_{1,1}^{\tau}+v_{2,2}^{\tau}\rt| dz\leq 8\pi.
\end{equation}
By Step 2 this implies $V(\na u, \Pi_{16\tau}\backslash \overline{B_{3\beta^{\frac{1}{8}}}(0)})\leq 16\pi$ and as 
$\tau$ is arbitrary  $V(\na u,\Omega\backslash \overline{B_{3\beta^{\frac{1}{8}}}(0)})\leq 16\pi$. $\Box$

%
%

\begin{a1}
\label{LL34}
Let $\Omega$ be a convex domain and $\lt|\Omega\triangle B_1(0)\rt|\leq \beta$.

Let $u(x)=d(x,\partial \Omega)$ and $\eta(x)=1-8\beta^{\frac{3}{32}}+\lt|x\rt|$. 
Define $\Gamma:=\lt\{x:u(x)=\eta(x)\rt\}$, we will show $\Gamma$ is the boundary of a convex set 
with $H^1(\Gamma)\leq c\beta^{\frac{3}{32}}$,
\begin{equation}
\label{azz7.5}
\Gamma\subset N_{c\beta^{\frac{3}{16}}}(\partial B_{4\beta^{\frac{3}{32}}}(0))
\end{equation}
and for any $\veps\in (0,\beta^{\frac{3}{16}}]$
\begin{equation}
\label{azz32}
\lt|N_{2\veps}(\Gamma)\rt|\leq c\veps \beta^{\frac{3}{32}}.
\end{equation}

\end{a1}
%
%
\em Proof of Lemma. \rm 

\em Step 1. \rm We will show $\Pi:=\lt\{x\in\Omega:\eta(x)\leq u(x)\rt\}$ is convex. 

\em Proof of Step 1. \rm Take $a,b\in \Pi$ and pick $\lm\in\lt[0,1\rt]$. Since $u$ is concave  
$u(\lm a+(1-\lm)b)\geq \lm u(a)+(1-\lm)u(b)$ and since $\eta$ is convex 
$\eta(\lm a+(1-\lm)b)\leq \lm \eta(a)+(1-\lm)\eta(b)$. 
Hence as $a,b\in\Pi$, $u(\lm a+(1-\lm)b)\geq \eta(\lm a+(1-\lm)b)$. Thus 
$\lt[a,b\rt]\subset \Pi$ and thus the set $\Pi$ is convex.\nl

\em Step 2. \rm We will establish (\ref{azz7.5}). 

\em Proof of Step 2. \rm Let $x\in \Gamma$ and let $b_x\in \partial \Omega$ be such that 
$\lt|x-b_x\rt|=u(x)$. So 
\begin{equation}
\label{zaz1}
1-8\beta^{\frac{3}{32}}+\lt|x\rt|=\lt|b_x-x\rt|.
\end{equation}
And thus $1-8\beta^{\frac{3}{32}}+\lt|x\rt|\geq\lt|b_x\rt|-\lt|x\rt|$, so using (\ref{abb22}) 
\begin{equation}
\label{ula341}
2\lt|x\rt|\geq  \lt|b_x\rt|-1+8\beta^{\frac{3}{32}}\geq 8\beta^{\frac{3}{32}}-c\sqrt{\beta}.
\end{equation}
Also from (\ref{zaz1}) we have 
\begin{equation}
\label{ula340}
\lt|x\rt|=\lt|b_x-x\rt|-(1-8\beta^{\frac{3}{32}})\overset{(\ref{abb22})}{\leq} 8\beta^{\frac{3}{32}}+\sqrt{\beta}.
\end{equation}
Now using Lemma \ref{LL25}, since $\na u(x)=\frac{x-b_x}{\lt|x-b_x\rt|}$ so 
\begin{eqnarray}
\label{abb72.3}
\lt|\frac{x}{\lt|x\rt|}-\frac{b_x}{\lt|b_x\rt|}\rt|&\leq& 
\lt|\frac{b_x-x}{\lt|b_x-x\rt|}-\frac{b_x}{\lt|b_x\rt|}\rt|
+\lt|\frac{x-b_x}{\lt|x-b_x\rt|}+\frac{x}{\lt|x\rt|}\rt|\nn\\
&\overset{(\ref{ula340}),(\ref{abb72.6})}{\leq}& c\beta^{\frac{3}{32}}
\end{eqnarray}
 so 
\begin{equation}
\label{abz45}
\lt|1-\frac{b_x}{\lt|b_x\rt|}\cdot \frac{x}{\lt|x\rt|}\rt|= 2^{-1}\lt|\frac{b_x}{\lt|b_x\rt|}-\frac{x}{\lt|x\rt|}\rt|^2\leq c\beta^{\frac{3}{16}}.
\end{equation}

Again by Lemma \ref{LL25} we have 
\begin{eqnarray}
\label{uu1}
\lt|\lt|b_x-x\rt|+(x-b_x)\cdot \frac{x}{\lt|x\rt|}\rt|&\leq&\lt|x-b_x\rt|\lt|\frac{x-b_x}{\lt|x-b_x\rt|}+\frac{x}{\lt|x\rt|}\rt|\nn\\
&\leq&2\lt|\na u(x)+\frac{x}{\lt|x\rt|}\rt|\nn\\
&\overset{(\ref{abb72.6})}{\leq}& c\beta^{\frac{3}{16}}
\end{eqnarray}
and thus 
\begin{eqnarray}
\lt|2x\cdot\frac{x}{\lt|x\rt|}-8\beta^{\frac{3}{32}}\rt|&\overset{(\ref{abz45})}{\leq}&
\lt|-8\beta^{\frac{3}{32}}+1-\frac{b_x}{\lt|b_x\rt|}\cdot \frac{x}{\lt|x\rt|}+2x\cdot \frac{x}{\lt|x\rt|}\rt|
+c\beta^{\frac{3}{16}}\nn\\
&=&\lt|1-8\beta^{\frac{3}{32}}+\lt|x\rt|-\lt(\frac{b_x}{\lt|b_x\rt|}-x\rt)
\cdot\frac{x}{\lt|x\rt|}\rt|+c\beta^{\frac{3}{16}}\nn\\
&\overset{(\ref{zaz1})}{=}&\lt|\lt|b_x-x\rt|-\lt(\frac{b_x}{\lt|b_x\rt|}-x\rt)\cdot \frac{x}{\lt|x\rt|}\rt|+
c\beta^{\frac{3}{16}}\nn\\
&\overset{(\ref{abb22})}{\leq}&\lt|\lt|b_x-x\rt|+\lt(x-b_x\rt)\cdot \frac{x}{\lt|x\rt|}\rt|+
c\beta^{\frac{3}{16}}\nn\\
&\overset{(\ref{uu1})}{\leq}&c\beta^{\frac{3}{16}} \nn
\end{eqnarray}
hence $\lt|2\lt|x\rt|-8\beta^{\frac{3}{32}}\rt|\leq c\beta^{\frac{3}{16}}$ for any $x\in \Gamma$, so 
(\ref{azz7.5}) is established.

Since (\ref{azz7.5}) implies the diameter of $\Pi$ is bounded by $c\beta^{\frac{3}{32}}$ and since $\Pi$ is a 
convex set it follows immediately that $H^1(\Gamma)\leq c\beta^{\frac{3}{32}}$. 

Now the set $\Gamma$ equipped with the Euclidean norm is a boundly compact metric space. So by applying the 
5r Covering Theorem (Theorem 2.1 \cite{mat}) we can find a disjoint collection of balls $B_{2\veps}(x_1), B_{2\veps}(x_2), \dots B_{2\veps}(x_M)$ with $x_1,x_2,\dots x_M\in \Gamma$ such that $\Gamma\subset \bigcup_{i=1}^{n} B_{10\veps}(x_i)$. This 
implies $N_{2\veps}(\Gamma)\subset \bigcup_{i=1}^{n} B_{20\veps}(x_i)$. 
Since $H^1(\Gamma)\leq c\beta^{\frac{3}{32}}$ so $M\leq c\veps^{-1}\beta^{\frac{3}{32}}$ and thus  
$\lt|N_{2\veps}(\Gamma)\rt|\leq c\veps \beta^{\frac{3}{32}}$ which establishes (\ref{azz32}). $\Box$

%
%

\begin{a1} 
\label{LLA1}

Let $\Omega$ be a convex set. Let $\beta=\lt|\Omega\triangle B_1(0)\rt|$. 
Let 
$$
w(z):=\min\lt\{d(z,\partial \Omega),1-8\beta^{\frac{3}{32}}+\lt|z\rt|\rt\}.
$$
We will show $\na w\in SBV(\Omega:S^1)$ and 
\begin{equation}
\label{qband03}
\int_{J_{\na w}\cap \Omega} \lt|\na w^{+}-\na w^{-}\rt|^3 dH^1 \leq 
c\beta^{\frac{3}{32}}.
\end{equation}
\end{a1}
\em Proof. \rm By Lemma \ref{LL12} we know $\na u\in BV(\Omega\backslash B_{3\beta^{\frac{1}{3}}}(0))$ 
and $V(\na u,\Omega\backslash B_{3\beta^{\frac{1}{3}}}(0))\leq 8\pi$. This implies 
\begin{equation}
\label{qband01}
\int_{(\Omega\backslash B_{3\beta^{\frac{1}{8}}}(0))\cap S_{\na u}} 
\lt|\na u^{+}-\na u^{-}\rt| dH^1 \leq 8\pi.
\end{equation}
Now by Lemma \ref{LL25} (\ref{abb72.6}) for any 
$x\in  (\Omega\backslash B_{3\beta^{\frac{1}{8}}}(0))\cap S_{\na u}$ we have 
$\lt|\na u^{+}(x)-\na u^{-}(x)\rt|\leq c\beta^{\frac{3}{16}}$. So 
\begin{eqnarray}
\label{gband011}
\int_{(\Omega\backslash B_{3\beta^{\frac{1}{8}}}(0))\cap S_{\na u}} 
\lt|\na u^{+}-\na u^{-}\rt|^3 dH^1 &\leq&
c\beta^{\frac{3}{8}}\int_{(\Omega\backslash B_{3\beta^{\frac{1}{8}}}(0))\cap S_{\na u}} 
\lt|\na u^{+}-\na u^{-}\rt| dH^1\nn\\
&\overset{(\ref{qband01})}{\leq}&c\beta^{\frac{3}{8}}.
\end{eqnarray}

As in Lemma \ref{LL34} let 
$\Pi:=\lt\{x:u(x)\leq 1-8\beta^{\frac{3}{32}}+\lt|x\rt|\rt\}$ and 
$\Gamma:=\partial \Pi$. Since $\Pi$ is convex it is also a 
set of finite perimeter. Let $\eta(z)=1-8\beta^{\frac{3}{32}}+\lt|x\rt|$, 
it is clear $w(z)=\cha_{\Pi} \eta(z)+\cha_{\Omega\backslash \Pi} u(z)$.  By 
Theorem 3.83 \cite{amb}, $\na w\in BV(\Omega:S^1)$. We know 
by Lemma \ref{LL34}, $H^1(\Gamma)\leq c\beta^{\frac{3}{32}}$. Now for any 
$x\in \Gamma$, since 
$\na w^{+}(x), \na w^{-}(x)\in S^1$, $\lt|\na w^{+}(x)-\na w^{-}(x)\rt|\leq 2$. So 
\begin{eqnarray}
\int_{J_{\na w}} \lt|\na w^{+}-\na w^{-}\rt|^3 dH^1 
&=&\int_{J_{\na w}\cap (\Omega\backslash \Pi)} \lt|\na w^{+}-\na w^{-}\rt|^3 dH^1 
+\int_{J_{\na w}\cap \Pi} \lt|\na w^{+}-\na w^{-}\rt|^3 dH^1 \nn\\
&\overset{(\ref{gband011})}{\leq}&c\beta^{\frac{3}{8}}+8H^1(\Gamma)\nn\\
&\leq&c\beta^{\frac{3}{32}}.\;\;\;\Box\nn 
\end{eqnarray}

\subsection{Proof of Proposition \ref{PP7} completed} By Lemma \ref{LLA1} we know 
that $w\in BV(\Omega,S^1)$ we can apply Theorem 1 of \cite{conti-del} or Corollary 1.1 \cite{poli} to find a 
sequence $u^{\ep}$ that satisfies $u^{\ep}(z)=0$ and $\na u^{\ep}(z)\cdot \eta_z=1$ for $z\in \partial \Omega$ (where $\eta_z$ is the inward 
pointing unit normal to $\partial \Omega$ at $z$) such that 
\begin{eqnarray}
\limsup_{\ep\rightarrow 0}\int_{\Omega} \ep^{-1}\lt|1-\lt|\na u^{\ep}\rt|^2\rt|^2+\ep\lt|\na^2 u^{\ep}\rt|^2 dz&\leq& 
\int_{J_{\na w}\cap \Omega} \lt|\na w^{+}-\na w^{-}\rt|^3 dH^1 \nn\\
&\overset{(\ref{qband03})}{\leq}& c\beta^{\frac{3}{32}}.\;\;\;\;\Box\nn
\end{eqnarray}

\subsection{Proof of Corollary \ref{CC3}} Let $\beta=\inf_{a\in \Omega} \lt|\Omega\triangle B_1(a)\rt|$. Without loss of generality we can 
assume $\lt|\Omega\triangle B_{1}(0)\rt|\leq 2\beta$. 
So by Proposition \ref{PP7} we can find $\ep_0\in (0,1)$ such that for $\ep\in (0,\ep_0)$, any minimiser 
$u^{\ep}$ of $I_{\ep}$ defined on $\Omega$ satisfies 
\begin{equation}
\label{opl1}
\int_{\Omega} \ep^{-1}\lt|1-\lt|\na u^{\ep}\rt|^2\rt|^2+\ep\lt|\na^2 u^{\ep}\rt|^2 dz\leq c\beta^{\frac{3}{32}}.
\end{equation}
So we can apply Theorem \ref{T0} to conclude that 
\begin{equation}
\label{kband010}
\int_{\Omega} \lt|\na u^{\ep}(z)+\frac{z}{\lt|z\rt|}\rt|^2 dz\leq c\beta^{\frac{1}{5462}}.
\end{equation}
Applying Lemma \ref{LL25} we have 
\begin{equation}
\label{uza1}
\int_{\Omega\backslash B_{\beta^{\frac{1}{8}}}(0)} \lt|\na u^{\ep}-\na \zeta\rt|^2\leq c\beta^{\frac{1}{5462}}.
\end{equation}
Now 
\begin{eqnarray}
\int_{B_{\beta^{\frac{1}{8}}}(0)} \lt|\na u^{\ep}-\na \zeta\rt|^2 dz&\leq&\int_{B_{\beta^{\frac{1}{8}}}(0)} 
\lt|\na u^{\ep}\rt|^2+2\lt|\na u^{\ep}\rt|+1 dz\nn\\
&\leq&\int_{B_{\beta^{\frac{1}{8}}}(0)} \lt(\lt|1-\lt|\na u^{\ep}\rt|^2\rt|+c\rt) dz\nn\\
&\overset{(\ref{opl1})}{\leq}&c\beta^{\frac{3}{32}}\nn
\end{eqnarray}
together with (\ref{uza1}) this gives $\|u^{\ep}-\zeta\|_{W^{1,2}(\Omega)}\leq 
c\beta^{\frac{1}{5462}}$. $\Box$

%
%

\section{Proof of Corollary \ref{CC1}}
\label{sec5}

In this section we will show that given a convex domain $\Omega$ with $C^2$ boundary with curvature
bounded above by $\ep^{-\frac{1}{2}}$ and that satisfies $\lt|B_1(0)\triangle \Omega\rt|\leq \beta$ we will construct a
function $u$ with $I_{\ep}(u)\leq \beta^{\frac{3}{16}}$, this is the contents of Proposition \ref{PP1} below. 
The proof of Corollary \ref{CC1} will follows easily from this.

\begin{a5}
\label{PP1}

Let $\Omega$ be a convex body with $C^2$ boundary and with curvature
bounded above by $\ep^{-\frac{1}{2}}$ and $\lt|\Omega\triangle B_1(0)\rt|\leq \beta$. 
Let $\ep\in (0,\frac{\beta^{\frac{1}{2}}}{4}]$, there exists a function $C^{2}$ function
$\xi:\Omega\rightarrow \R$ which satisfies $\na \xi(z)\cdot \eta_z=1$ (where $\eta_z$ is the inward 
pointing unit normal to $\partial \Omega$ at $z$), $\xi(z)=0$ for 
$z\in \partial \Omega$ and for which
\begin{equation}
\label{ula201}
\int_{\Omega} \ep^{-1}\lt|1-\lt|\na \xi\rt|^2\rt|^2+\ep\lt|\na^2 \xi\rt|^2 dz\leq c\beta^{\frac{3}{32}}.
\end{equation}
\end{a5}

%
%

\subsection{Proof of Proposition \ref{PP1}}

We begin with a preliminary lemma.
\begin{a1}
\label{LL6}
Let $\phi:\R_{+}\rightarrow \R_{+}$ be a continuous function. Let $\rho$ denote the standard convolution kernel, 
i.e.\ $\int \rho =1$ and $\spt\rho \subset B_{\frac{3}{2}}(0)$ and define $\rho_h(z):=h^{-2}\rho(h^{-1}z)$.

Suppose $f:\R^n\rightarrow \R$ be an affine function, let $g(x)=f*\rho_{\phi\lt(x\rt)}(x)$ then
\begin{equation}
\label{abb12}
g\lt(x\rt)=f\lt(x\rt) \text{ for all }x\in \R^n.
\end{equation}
\end{a1}
\em Proof of Lemma. \rm Let $\eta=\na f$. As $f$ is affine $f(x-y)=f(x)-\eta\cdot y$
\begin{eqnarray}
\label{aab14}
g(x)&=&\int f(x-y)(\phi(x))^{-2}\rho(\phi(x)^{-1} y) dy \nn\\
&=&\int (f(x)-\eta\cdot y)(\phi(x))^{-2}\rho(\phi(x)^{-1} y) dy \nn\\
&=&f(x).\;\;\Box
\end{eqnarray}

%
%

\begin{a1}
\label{LL7}
Let $\ep>0$, suppose $\Omega$ is a convex body with $C^2$ boundary and with curvature
bounded above by $\ep^{-\frac{1}{2}}$. Let $u(x)=d(x,\partial \Omega)$. Let 
$\rho$ be the standard convolution kernel and $\rho_{\ep}(z):=\rho\lt(\frac{z}{\ep}\rt)\ep^{-2}$. 
We will construct a function $\psi:\Omega\cap N_{8\ep}(\partial \Omega)\rightarrow \R$ with $\psi=0$  on $\partial \Omega$ which satisfies the following properties
\begin{equation}
\label{ula351}
\int_{\Omega\cap N_{8\ep}(\partial \Omega)} \lt|1-\lt|\na \psi\rt|^2\rt|^2 dz\leq c \ep^{2},
\end{equation}
\begin{equation}
\label{ula357}
\int_{\Omega\cap N_{8\ep}(\partial \Omega)}   \lt|\na^2 \psi\rt|^2 dz\leq c,
\end{equation}
\begin{equation}
\label{ulazbc1}
\psi(z)=\lt[u*\rho_{\ep}\rt](z)\text{ for any }z\in \Omega\backslash N_{8\ep}(\partial \Omega)
\end{equation}
and 
\begin{equation}
\label{uz22}
\na \psi(z)=\eta_z\text{ for each }z\in \partial \Omega. 
\end{equation}

\end{a1}
\em Proof. \rm Let $w:\R_{+}\rightarrow \R_{+}$ be a smooth monotonic function with the following properties 
\begin{equation}
\label{ula70}
w\lt(z\rt)=
\begin{cases} z & \text{for }z\in \lt[0,\frac{\ep}{3}\rt)\\
 \ep &\text{ for } z\geq \ep
\end{cases}
\end{equation}
and $\sup \lt|\ddot{w}\rt|\leq c\ep^{-1}$. 

For any $x\in \Omega\cap N_{8\ep}(\partial \Omega)$ define 
\begin{equation}
\label{kband0202}
\phi(x)=w(u(x)). 
\end{equation}
We will 
convolve the function $u$ with convolution kernel 
$\rho_{\phi(x)}(z):=\rho\lt(\frac{z}{\phi(x)}\rt)/\lt(\phi(x)\rt)^{2}$. Since the convulsion 
kernel varies with $x$, when we differentiate $u*\rho_{\phi(x)}$, the derivative will involve a 
term with the derivative of $\rho_{\phi(x)}$. For this reason we need to calculate various partial 
derivatives of $\rho_{\phi(x)}$.

Since the curvature of $\partial \Omega$ is bounded above by 
$\ep^{-\frac{1}{2}}$, for any $x\in \Omega\cap N_{8\ep}(\partial \Omega)$ we have that there is 
one unique $b_x\in \partial \Omega$ such that $\lt|x-b_x\rt|=u(x)$. 
We define $\vs_x=\frac{x-b_x}{\lt|x-b_x\rt|}$, let 
$R=\lt(\begin{matrix} 0 & -1\\ 1 & 0 \end{matrix}\rt)$ and define $\omega_x=R \vs_x$.

Note $\vs_x=\eta_{b_x}$, i.e.\ the inward pointing unit normal to $\partial \Omega$ at $b_x$. 
Note also that for all small enough $h$, $b_x=b_{x+h\vs_x}$ so $u(x+h\vs_x)=h+u(x)$. Thus
\begin{eqnarray}
\phi_{,\vs_x}(x)&=&\lim_{h\rightarrow 0}
\frac{\phi(x+h\vs_x)-\phi(x)}{h}\nn\\
&=&\lim_{h\rightarrow 0} \frac{w(u(x)+h)-w(u(x))}{h}\nn\\
&=&\dot{w}(u(x)).\nn
\end{eqnarray} 

Note also that since $\lt|\na u(x)\rt|=1$ and 
$u_{,\vs_x}(x)=\lim_{h\rightarrow 0} 
\frac{u(x+h\vs_x)-u(x)}{h}=1$
so 
$$
u_{,\omega_x}(x)=\lim_{h\rightarrow 0} 
\frac{u(x+h\omega_x)-u(x)}{h}=0.
$$
Thus 
\begin{equation}
\label{ula81}
\phi_{,\omega_x}(x)=\dot{w}(u(x))u_{,\omega_x}(x)=0.
\end{equation}
So 
\begin{eqnarray}
\label{ula120}
\frac{\partial}{\partial \vs_x}\lt(\rho_{\phi(x)}(z)\rt)&=&
\frac{\partial}{\partial \vs_x}\lt(\rho\lt(\frac{z}{\phi(x)}\rt)\lt(\phi(x)\rt)^{-2}\rt)\nn\\
&=&-\na\rho\lt(\frac{z}{\phi(x)}\rt)\cdot z \frac{\phi_{,\vs_x}(x)}{(\phi(x))^4}-
2\rho\lt(\frac{z}{\phi(x)}\rt)\frac{\phi_{,\vs_x}(x)}{\lt(\phi(x)\rt)^{3}}
\end{eqnarray}
and 
\begin{equation}
\label{ula83}
\frac{\partial}{\partial \omega_x}\lt(\rho_{\phi(x)}(z)\rt)=0.
\end{equation}
Define 
\begin{equation}
\label{qwer1}
\psi(x):=\lt\{\begin{array}{ll} \int u(x-z)\rho_{\phi\lt(x\rt)}(z) dz 
&\text{ for } x\in \Omega\\
0 &\text{ for }  x\not\in \Omega \end{array} \rt..
\end{equation}
Now 
\begin{eqnarray}
\label{aqa1}
\psi_{,\vs_x}(x)&\overset{(\ref{ula120})}{=}&\int u_{,\vs_x}(x-z) \rho_{\phi(x)}(z) dz\nn\\
&~&-\int u(x-z)\lt(\na\rho\lt(\frac{z}{\phi(x)}\rt)\cdot z  \frac{\phi_{,\vs_x}(x)}{(\phi(x))^4} +
2\rho\lt(\frac{z}{\phi(x)}\rt)\frac{\phi_{,\vs_x}(x)}{(\phi(x))^3} \rt) dz
\end{eqnarray}
In the same way it is easy to see $\psi_{,\omega_x}(x)=\int u_{,\omega_x}(x-z)\rho_{\phi(x)}(z) dz$ and so 
\begin{eqnarray}
\label{ula21}
\psi_{,\vs_x \omega_x}(x)&=&\int u_{,\omega_x \vs_x}(x-z)\rho_{\phi(x)}(z) dz
+\int u_{,\omega_x}(x-z)\frac{\partial}{\partial \vs_x}\lt(\rho_{\phi(x)}(z)\rt) dz.
\end{eqnarray}
And 
\begin{equation}
\label{ula21.5}
\psi_{,\omega_x \omega_x}(x)=\int u_{,\omega_x \omega_x}(x-z)\rho_{\phi(x)}(z) dz.
\end{equation}
Finally 
\begin{eqnarray}
\label{fineq35}
\psi_{,\vs_x \vs_x}(x)&=&\int u_{,\vs_x \vs_x}(x-z)\rho_{\phi(x)}(z)
+2u_{,\vs_x}(x-z)\frac{\partial }{\partial \vs_x}\lt(\rho_{\phi(x)}(z)\rt) dz\nn\\
&~&\qd+\int u(x-y)\frac{\partial^2 }{\partial^2 \vs_x}\lt(\rho_{\phi(x)}(z)\rt) dz
\end{eqnarray}
each term will be estimated later in Step 4. \nl

%
%

\em Step 1. \rm We will show 
\begin{equation}
\label{fineq60}
\lt|\na ^2 u(x)\rt|\leq c\ep^{-\frac{1}{2}}\text{ for any }x\in N_{\frac{\sqrt{\ep}}{3}}(\partial \Omega).
\end{equation}

\em Proof of Step 1. \rm Let $b_x\in \partial \Omega$ be such that 
$\mathrm{dist}(x,\partial \Omega)=\lt|x-b_x\rt|$. We start by showing 
\begin{equation}
\label{fineq30}
\lt|\na u(x)-\na u(y)\rt|\leq c\ep^{-\frac{1}{2}}\lt|x-y\rt|\text{ for any }x\in N_{\frac{\sqrt{\ep}}{3}}(\partial \Omega), 
y\in B_{\frac{\ep}{6}}(x).
\end{equation}
Now recall $\frac{y-b_y}{\lt|y-b_y\rt|}=\eta_{b_y}$, 
$\frac{x-b_x}{\lt|x-b_x\rt|}=\eta_{b_x}$. We have two cases to consider. Firstly the case that 
$(b_x+\R_{+}\eta_{b_x})\cap (b_y+\R_{+}\eta_{b_y})=\emptyset$. In this case since 
$\Omega$ is convex this implies $\eta_{b_x}=\eta_{b_y}$. Thus as 
$\lt|\na u(x)-\na u(y)\rt|=\lt|\frac{y-b_y}{\lt|y-b_y\rt|}-\frac{x-b_x}{\lt|x-b_x\rt|}\rt|=\lt|\eta_{b_x}-\eta_{b_y}\rt|=0$ so 
(\ref{fineq30}) is established. 

Now suppose we have the case that $\pi:=(b_x+\R_{+}\eta_{b_x})\cap (b_y+\R_{+}\eta_{b_y})\not=\emptyset$. Then let 
\begin{equation}
\label{wwwe1}
\theta=\arccos\lt(\frac{b_y-y}{\lt|b_y-y\rt|}\cdot \frac{b_x-x}{\lt|b_x-x\rt|}\rt). 
\end{equation}
Since the curvature of $\partial \Omega$ is bounded by $\ep^{-\frac{1}{2}}$ we know that 
$\pi\not \in N_{\sqrt{\ep}}(\partial \Omega)$. Consider the triangle whose corners are 
$x,y,\pi$, which we denote by $T(x,y,\pi)$. The 
angle at corner $\pi$ is $\theta$. Now since $\lt|x-y\rt|\leq \frac{\ep}{6}$,  
$\lt|x-\pi\rt|\geq \frac{\sqrt{\ep}}{2}$,  $\lt|y-\pi\rt|\geq \frac{\sqrt{\ep}}{2}$. 
So as $\lt|\lt|x-\pi\rt|-\lt|y-\pi\rt|\rt|\leq \lt|x-y\rt|\leq \frac{\ep}{6}$. Thus 
$$
\frac{\ep^2}{36}\geq \lt|\lt|x-\pi\rt|-\lt|y-\pi\rt|\rt|^2=
\lt|2\lt|x-\pi\rt|\lt|y-\pi\rt|-\lt|x-\pi\rt|^2-\lt|y-\pi\rt|^2\rt|.
$$
Thus by the law of cosines 
\begin{eqnarray}
2\lt|x-\pi\rt|\lt|y-\pi\rt|\cos\theta&=&\lt|x-\pi\rt|^2+\lt|y-\pi\rt|^2-\lt|x-y\rt|^2\nn\\
&\geq& \lt|x-\pi\rt|^2+\lt|y-\pi\rt|^2-\frac{\ep^2}{36}\nn\\
&\geq&2\lt|x-\pi\rt|\lt|y-\pi\rt|-\frac{\ep^2}{36}.\nn
\end{eqnarray}
Which implies $\cos \theta\geq 1-c\ep$ and so $\lt|\theta\rt|\leq c\sqrt{\ep}$. 

Let $\ti{y}:=\lt[b_y,\pi\rt]\cap \partial B_{\lt|x-\pi\rt|}(x)$, since $\lt|\theta\rt|\leq c\sqrt{\ep}$ we 
have $\lt|x-\ti{y}\rt|\leq \frac{11}{10}\lt|x-y\rt|$. Consider the triangle $T(x,\ti{y},\pi)$. Note the 
angle of this triangle at $\pi$ is $\theta$ and denoting the angle at $x$ by $\psi$ we have $\psi\sim \frac{\pi}{2}$.

Then by the law of sins, 
$$
\frac{\lt|x-\ti{y}\rt|}{\sin\theta}=\frac{\lt|\ti{y}-\pi\rt|}{\sin\psi}\geq 
\frac{\lt|\ti{y}-\pi\rt|}{2}\geq \frac{\sqrt{\ep}}{4}.
$$
So $4\frac{\lt|x-\ti{y}\rt|}{\sqrt{\ep}}\geq \sin\theta$ which 
gives $\lt|\theta\rt|\leq \frac{c\lt|x-\ti{y}\rt|}{\sqrt{\ep}}\leq \frac{c\lt|x-y\rt|}{\sqrt{\ep}}$. So as 
$\na u(x)=\frac{x-b_x}{\lt|x-b_x\rt|}$ and $\na u(y)=\frac{y-b_y}{\lt|y-b_y\rt|}$, (recalling the definition of $\theta$ from 
(\ref{wwwe1})) $\lt|\na u(x)-\na u(y)\rt|\leq c\arccos\lt(\na u(x)\cdot \na u(y)\rt)\leq \frac{c\lt|x-y\rt|}{\sqrt{\ep}}$. So 
(\ref{fineq30}) is established. Thus letting $y\rightarrow x$ we have that $\lt|\na ^2 u(x)\rt|\leq c\ep^{-\frac{1}{2}}$ and this 
completes the proof of Step 1. \nl


\em Step 2. \rm For any $x\in N_{16\ep}(\partial\Omega)\cap \Omega$ we have
\begin{equation}
\label{abb5}
\sup\lt\{\lt|\na u(z)-\vs_x\rt|:z\in B_{16 u\lt(x\rt)}(x)\cap\Omega\rt\}\leq c\ep^{-\frac{1}{2}}u(x).
\end{equation}
\em Proof of Step 2. \rm Since $\partial \Omega$ has curvature less than $\ep^{-\frac{1}{2}}$ for any
$x_1,x_2\in \partial \Omega$, $\lt[x_1,x_1+\ep^{\frac{1}{2}}\eta_{x_1}\rt]\cap \lt[x_2,x_2+\ep^{\frac{1}{2}}\eta_{x_2}\rt]=\emptyset$. So for any
$x_1,x_2\in B_{32 u(x)}(x)\cap \partial\Omega$, $\lt|\eta_{x_1}-\eta_{x_2}\rt|\leq \ep^{-\frac{1}{2}}
H^1(B_{32 u(x)}(x)\cap \partial\Omega)$. Note as $\Omega\cap B_{32 u(x)}(x)$ is convex
and $\partial \Omega \cap B_{32 u(x)}(x)\subset \partial ( \Omega\cap B_{32 u(x)}(x))$ so
$H^1(\partial \Omega \cap B_{32 u(x)}(x))\leq c u(x)$. 
Hence $\lt|\eta_{x_1}-\eta_{x_2}\rt|\leq c\ep^{-\frac{1}{2}} u(x)\leq c\sqrt{\ep}$ so it is clear that 
\begin{equation}
\label{aza1}
B_{16 u(x)}(x)\cap \Omega\subset \bigcup_{z\in  \partial \Omega \cap B_{32 u(x)}(x)}
\lt[z,z+\sqrt{\ep}\eta_z\rt].
\end{equation}

For any $z\in B_{16 u(x)}(x)\cap \Omega$ we have $\na u(z)=\frac{z-b_z}{\lt|z-b_z\rt|}=\eta_{b_z}$ where $b_z$ 
is such that $\lt|z-b_z\rt|=d\lt(z,\partial\Omega\rt)$. So for any 
$z_1,z_2\in B_{16u(x)}(x)\cap \Omega$ by (\ref{aza1}) we have that $b_{z_1}, b_{z_2}\in \partial \Omega\cap 
B_{32 u(x)}(x)$, so 
$\lt|\na u(z_1)-\na u(z_2)\rt|=\lt|\eta_{b_{z_1}}-\eta_{b_{z_2}}\rt|\leq c\ep^{-\frac{1}{2}}u(x)$.\nl


\em Step 3. \rm For any $x\in N_{8\ep}(\partial\Omega)\cap \Omega$ we have
\begin{equation}
\label{ula400.6}
\lt|\lt|\na \psi(x)\rt|-1\rt|\leq c\sqrt{\ep}.
\end{equation}
And 
\begin{equation}
\label{uz1}
\lim_{y\rightarrow z} \na \psi(y)=\eta_z.
\end{equation}

\em Proof of Step 3. \rm From (\ref{aqa1}) we have
\begin{eqnarray}
\label{ula81}
&~&\lt|\psi_{,\vs_x}(x)-1\rt|\nn\\
&~&\qd\qd\leq\overbrace{\lt|\int (u_{,\vs_x}(x-z)-1)\rho\lt(\frac{z}{\phi(x)}\rt)\lt(\phi(x)\rt)^{-2} dz\rt|}^B\nn\\
&~&\qd\qd\qd
+\overbrace{\lt|\int \frac{-u(x-z)\phi_{,\vs_x}(x)}{\lt(\phi(x)\rt)^{3}}
\lt(\na \rho\lt(\frac{z}{\phi(x)}\rt)\cdot \frac{z}{\phi(x)}+2\rho\lt(\frac{z}{\phi(x}\rt)\rt) dz\rt|}^C.
\end{eqnarray}

Now for any $z\in \spt\rho_{\phi(x)}$ we have that
$\na u(x-z)=u_{,\vs_x}(x-z)\vs_x+u_{,\omega_x}(x-z) \omega_x$ now since 
$\spt\rho_{\phi(x)}\subset B_{2\phi(x)}(0)\subset B_{2 u(x)}(0)$ so for any 
$z\in \spt\rho_{\phi(x)}$ by (\ref{abb5}) from Step 2 we have 
$\lt|\na u(x-z)-\vs_x\rt|\leq c\ep^{-\frac{1}{2}}u(x)$ and thus 
\begin{equation}
\label{abb16}
\lt|u_{,\vs_x}(x-z)-1\rt|\leq c\ep^{-\frac{1}{2}} u(x)\text{ for any }z\in \spt\rho_{\phi(x)}
\end{equation}
so (noting $u(x)\leq c\phi(x)$ for any $x\in N_{8\ep}(\partial\Omega)\cap \Omega$)
\begin{equation}
\label{abb17}
B\leq c u(x) \ep^{-\frac{1}{2}}<c\phi(x)\ep^{-\frac{1}{2}}.
\end{equation}

Also defining $w=\Xint{-}_{B_{\phi\lt(x\rt)}(x)} \na u$ 
\begin{equation}
\label{aza2}
\lt|w-\vs_x\rt|=\lt|\Xint{-}_{B_{\phi\lt(x\rt)}(x)} \lt(\na u(z)-\vs_x\rt) dz \rt|
\overset{(\ref{abb5})}{\leq}  c\ep^{-\frac{1}{2}} \phi\lt(x\rt).
\end{equation}
So by Poincare's inequality there exists affine function $l_w$ with $\na l_w=w$
\begin{eqnarray}
\label{aza}
\Xint{-}_{B_{\phi\lt(x\rt)}\lt(x\rt)} \lt|u\lt(z\rt)-l_w\lt(z\rt)\rt| dz&\leq&
c \phi\lt(x\rt)\Xint{-}_{B_{\phi\lt(x\rt)}(x)} \lt|\na u(z)-w\rt| dz\nn\\
&\leq&c \phi\lt(x\rt)
\lt(\Xint{-}_{B_{\phi\lt(x\rt)}(x)} \lt|\na u(z)-\vs_x\rt| dz+c\lt|w-\vs_x\rt|\rt)\nn\\
&\overset{(\ref{abb5}),(\ref{aza2})}{\leq}&
c  \ep^{-\frac{1}{2}} (\phi\lt(x\rt))^2.
\end{eqnarray}

Now using (\ref{aza2}), again for the appropriate choice of affine function $l_{\vs_x}$ with $\na l_{\vs_x}=\vs_x$ 
we have by Poincare's inequality 
$$
\Xint{-}_{B_{\phi\lt(x\rt)}\lt(x\rt)}  \lt| l_{\vs_x}\lt(z\rt)-l_{w}(z)\rt| dz\leq 
c\phi\lt(x\rt)\Xint{-}_{B_{\phi\lt(x\rt)}} \lt|w-\vs_x\rt| dz\overset{(\ref{aza2})}{\leq} 
c\ep^{-\frac{1}{2}} (\phi\lt(x\rt))^2
$$
with (\ref{aza}) gives 
\begin{equation}
\label{abb3}
\Xint{-}_{B_{\phi\lt(x\rt)}\lt(x\rt)}  \lt| l_{\vs_x}\lt(z\rt)-u(z)\rt| dz\leq
c\ep^{-\frac{1}{2}} (\phi\lt(x\rt))^2.
\end{equation}
Let $g$ be defined by $g(y)=l_{\vs_x}*\rho_{\phi(y)}(y)$, note by Lemma \ref{LL6} we have 
$\na g(y)=\vs_x$ for any $y\in \Omega$ and hence $g_{,\vs_x}(x)=1$ 
and as 
\begin{eqnarray}
g_{,\vs_x}(x)&=&\int \rho\lt(\frac{z}{\phi(x)}\rt)
\lt(\phi(x)\rt)^{-2} dz\nn\\
&~&
-\int\frac{l_{\vs_x}(x-z)}{\lt(\phi\lt(x\rt)\rt)^3} 
\phi_{,\vs_x}\lt(x\rt) \lt(\na \rho\lt(\frac{z}{\phi(x)}\rt)\cdot z
\lt(\phi(x)\rt)^{-1}+2 \rho\lt(\frac{z}{\phi(x)}\rt)\rt) dz\nn\\
&=&1-\int \frac{l_{\vs_x}(x-z)}{\lt(\phi\lt(x\rt)\rt)^3}
\phi_{,\vs_x}\lt(x\rt) \lt(\na \rho\lt(\frac{z}{\phi(x)}\rt)\cdot z
\lt(\phi(x)\rt)^{-1}+2 \rho\lt(\frac{z}{\phi(x)}\rt)\rt) dz.\nn
\end{eqnarray} 

Thus
\begin{eqnarray}
0&=&\int \frac{l_{\vs_x}(x-z)}{\lt(\phi\lt(x\rt)\rt)^3} 
\phi_{,\vs_x}\lt(x\rt) \lt(\na \rho\lt(\frac{z}{\phi(x)}\rt)\cdot z
\lt(\phi(x)\rt)^{-1}+2 \rho\lt(\frac{z}{\phi(x)}\rt)\rt) dz 
\end{eqnarray}

So
\begin{eqnarray}
\label{ula82}
C&\leq& \int \frac{\lt|l_{\vs_x}(x-z)-u(x-z)\rt|}{\lt(\phi\lt(x\rt)\rt)^3}
\lt|\phi_{,\vs_x}\lt(x\rt) \lt(\na \rho\lt(\frac{z}{\phi(x)}\rt)\cdot z
\lt(\phi(x)\rt)^{-1}+2 \rho\lt(\frac{z}{\phi(x)}\rt)\rt)\rt| dz \nn\\
&\leq& c(\phi(x))^{-3}\int_{B_{\phi(x)}(x)} \lt|l_{\vs_x}(z)-u(z)\rt| dz\nn\\
&\overset{(\ref{abb3})}{\leq}&c\ep^{-\frac{1}{2}}\phi(x).
\end{eqnarray}

Since $x\in N_{8\ep}(\partial\Omega)\cap \Omega$ we know $\phi(x)\leq c\ep$ 
applying (\ref{ula82}) and (\ref{abb17}) to (\ref{ula81}) gives 
\begin{equation}
\label{uz2}
\lt|\psi_{,\vs_x}(x)-1\rt|\leq c\ep^{-\frac{1}{2}}\phi\lt(x\rt)\leq c\sqrt{\ep}. 
\end{equation}
Now using that $u_{,\omega_x}(x)=0$ we have that
\begin{eqnarray}
\label{ffbb25}
\lt|\psi_{,\omega_x}(x)\rt|&\leq&\lt|\int u_{,\omega_x}(x-z)\rho_{\phi(x)}(z) dz\rt|\nn\\
&\leq&\int \lt|u_{,\omega_x}(x-z)-u_{,\omega_x}(x)\rt| \rho_{\phi(x)}(z) dz\nn\\
&\overset{(\ref{fineq60})}{\leq}&c\ep^{-\frac{1}{2}}\phi(x)\int \rho_{\phi(x)}(z) dz\nn\\
&\leq& c\ep^{\frac{1}{2}}\phi(x).
\end{eqnarray}
Thus $\lt|\na \psi(x)-\vs_x\rt|\leq c\sqrt{\ep}$ and 
(\ref{ula400.6}) follows easily. Also for (\ref{uz2}), (\ref{ffbb25}) we know 
$\lt|\na \psi(x)-\eta_{b_x}\rt|\leq c\ep^{-\frac{1}{2}}\phi(x)$ and (\ref{uz1}) follows. 
This completes the proof of Step 3.\nl

\em Step 4. \rm We will show 
\begin{equation}
\label{dd4}
\lt|\na^2 \psi(x)\rt|\leq c\ep^{-\frac{1}{2}}\text{ for any }x\in N_{8\ep}(\partial\Omega)\cap \Omega.
\end{equation}

\em Proof Step 4. \rm 
We will estimate the terms in (\ref{fineq35}) one by one.  First note 
\begin{eqnarray}
&~& \int u(x-y)\frac{\partial^2 }{\partial^2 \vs_x}\lt(\rho_{\phi(x)}(z)\rt) dz\nn\\
&~&\qd=\int u(x-z)\partial_{\vs_x}\lt(\sum_{k=1}^2
-\rho_{,k}\lt(\frac{z}{\phi\lt(x\rt)}\rt)
\frac{z_k \phi_{,\vs_x}(x) }{\lt(\phi\lt(x\rt)\rt)^4}-2\rho\lt(\frac{z}{\phi\lt(x\rt)}\rt)
\frac{\phi_{,\vs_x}(x)}{\lt(\phi\lt(x\rt)\rt)^3}\rt)dz\nn\\
&~&\qd=\int u(x-z)\lt(\sum_{k,l=1}^2 \rho_{,kl}\lt(\frac{z}{\phi\lt(x\rt)}\rt)
\frac{\lt(\phi_{,\vs_x}(x)\rt)^2}{(\phi(x))^6} z_k z_l-\sum_{k=1}^2 \rho_{,k}\lt(\frac{z}{\phi\lt(x\rt)}\rt) z_k
\partial_{\vs_x}\lt(\frac{\phi_{,\vs_x}(x)}{(\phi(x))^4}\rt)\rt.\nn\\
&~&\qd\qd\qd\qd\qd\qd\qd\qd \lt.+2\sum_{m=1}^2 \rho_{,m}\lt(\frac{z}{\phi\lt(x\rt)}\rt) z_m
\frac{(\phi_{,\vs_x}(x))^2}{(\phi(x))^5}-2\rho\lt(\frac{z}{\phi\lt(x\rt)}\rt)
\partial_{\vs_x}\lt(\frac{\phi_{,\vs_x}(x)}{\lt(\phi\lt(x\rt)\rt)^3} \rt)\rt) dz\nn
\end{eqnarray}
Note
\begin{equation}
\partial_{\vs_x}\lt(\frac{\phi_{,\vs_x}(x)}{\lt(\phi\lt(x\rt)\rt)^3}\rt)=
\frac{-3(\phi_{,\vs_x}(x))^2}{\lt(\phi\lt(x\rt)\rt)^4}+ \frac{\phi_{,\vs_x \vs_x}(x)}{\lt(\phi\lt(x\rt)\rt)^3}\nn
\end{equation}
and
\begin{equation}
\partial_{\vs_x}\lt(\frac{\phi_{,\vs_x}(x)}{(\phi(x))^4}\rt)=
\frac{-4(\phi_{,\vs_x}(x))^2}{\lt(\phi\lt(x\rt)\rt)^5}
+ \frac{\phi_{,\vs_x \vs_x}(x)}{\lt(\phi\lt(x\rt)\rt)^4}.
\end{equation}
So
\begin{eqnarray}
\label{abb9}
&~&\int u(x-y)\frac{\partial^2 }{\partial^2 \vs_x}\lt(\rho_{\phi(x)}(z)\rt) dz\nn\\
&~&\qd\qd\qd
=\int u(x-z)\lt(\lt(\na^2 \rho\lt(\frac{z}{\phi\lt(x\rt)}\rt):z\otimes z\rt)\frac{(\phi_{,\vs_x}(x))^2}{\lt(\phi\lt(x\rt)\rt)^6}\rt.\nn\\
&~&\qd\qd\qd\qd
\lt.
+\lt(-\frac{\phi_{,\vs_x \vs_x}(x)}{\lt(\phi\lt(x\rt)\rt)^4}
+\frac{6(\phi_{,\vs_x}(x))^2}{\lt(\phi\lt(x\rt)\rt)^5}\rt) \na \rho\lt(\frac{z}{\phi\lt(x\rt)}\rt)\cdot z \rt.\nn\\
&~&\qd\qd\qd\qd +\lt.\lt(\frac{6(\phi_{,\vs_x}(x))^2}{\lt(\phi\lt(x\rt)\rt)^4}
-\frac{2\phi_{,\vs_x \vs_x}(x)}{\lt(\phi\lt(x\rt)\rt)^3}\rt) \rho\lt(\frac{z}{\phi\lt(x\rt)}\rt)\rt)  dz.
\end{eqnarray}

From Step 2 (\ref{abb3}) we know the existence of an affine function $l_{\vs_x}$ with
$\na l_{\vs_x}=\vs_x$ with $\Xint{-}_{B_{\phi\lt(x\rt)}\lt(x\rt)}  \lt|u-l_{\vs_x}\rt| dz\leq
c\ep^{\frac{1}{2}} \phi\lt(x\rt)$. 
Let $g(x):=l_{\vs_x}* \rho_{\phi\lt(x\rt)}\lt(x\rt)$ so by Lemma \ref{LL6} we know
$g_{,\vs_x \vs_x}(x)=0$. By following through the same calculation that gave (\ref{abb9}) we have
\begin{eqnarray}
\label{abb10}
0&=&\int l_{\vs_x}(x-z)\lt(\lt(\na^2 \rho\lt(\frac{z}{\phi\lt(x\rt)}\rt):z\otimes z\rt)\frac{(\phi_{,\vs_x}(x))^2}{\lt(\phi\lt(x\rt)\rt)^6}\rt.\nn\\
&~&\qd\qd\qd\qd\qd\qd
\lt.
+\lt(-\frac{\phi_{,\vs_x \vs_x}(x)}{\lt(\phi\lt(x\rt)\rt)^4}
+\frac{6(\phi_{,\vs_x}(x))^4}{\lt(\phi\lt(x\rt)\rt)^5}\rt) \na \rho\lt(\frac{z}{\phi\lt(x\rt)}\rt)\cdot z \rt.\nn\\
&~&\qd\qd\qd\qd\qd\qd\qd +\lt.\lt(\frac{6(\phi_{,\vs_x}(x))^2}{\lt(\phi\lt(x\rt)\rt)^4}
-\frac{2\phi_{,\vs_x \vs_x}(x)}{\lt(\phi\lt(x\rt)\rt)^3}\rt) \rho\lt(\frac{z}{\phi\lt(x\rt)}\rt)\rt)  dz.
\end{eqnarray}

Note for $x\in N_{8\ep}(\partial \Omega)\cap \Omega$, 
$\lt|\phi_{,\vs_x}(x)\rt|\leq c$ and $\lt|\phi_{,\vs_x \vs_x}(x)\rt|\leq c\ep^{-1}\leq c(\phi(x))^{-1}$. So applying (\ref{abb10}) to (\ref{abb9})
\begin{eqnarray}
\label{fineq99}
&~&\lt|\int u(x-z)\frac{\partial^2}{\partial \vs_x^2}\lt(\rho_{\phi(x)}(z)\rt)\rt| \nn\\
&~&\qd\qd\leq \int \lt|u\lt(x-z\rt)-l_{\vs_x}(x-z)\rt|\lt|\lt(\na^2 \rho\lt(\frac{z}{\phi\lt(x\rt)}\rt):z\otimes z\rt)\frac{(\phi_{,\vs_x}(x))^2}{\lt(\phi\lt(x\rt)\rt)^6}\rt.\nn\\
&~&\qd\qd\qd\qd
\lt.
+\lt(\frac{-\phi_{,\vs_x \vs_x}(x)}{\lt(\phi\lt(x\rt)\rt)^4}
+\frac{6(\phi_{,\vs_x}(x))^4}{\lt(\phi\lt(x\rt)\rt)^5}\rt) \na \rho\lt(\frac{z}{\phi\lt(x\rt)}\rt)\cdot z \rt.\nn\\
&~&\qd\qd\qd\qd
+\lt.\lt(\frac{6(\phi_{,\vs_x}(x))^2}{\lt(\phi\lt(x\rt)\rt)^4}-\frac{2\phi_{,\vs_x \vs_x}(x)}{\lt(\phi\lt(x\rt)\rt)^3}\rt) \rho\lt(\frac{z}{\phi\lt(x\rt)}\rt)\rt|  dz\nn\\
&~&\qd\qd\leq c\int_{B_{\phi\lt(x\rt)}(0)} \frac{\lt|u\lt(x-z\rt)-l_{\vs_x}(x-z)\rt|}
{(\phi\lt(x\rt))^{4}} dz
\lt(\|\na^2 \rho\|_{\infty}+\|\na \rho\|_{\infty}+\|\rho\|_{\infty}\rt)\nn\\
&~&\qd\qd \leq c\int_{B_{\phi\lt(x\rt)}(x)} \lt|u\lt(z\rt)-l_{\vs_x}(z)\rt|\lt(\phi\lt(x\rt)\rt)^{-4} dz\nn\\
&~&\qd\qd \overset{(\ref{abb3})}{\leq}c \ep^{-\frac{1}{2}}.
\end{eqnarray}

Define $h(x):=\int \rho_{\phi(x)}(z) dz$, so note that $h\equiv 1$ and so 
$\frac{\partial h}{\partial \vs_x}(x)=\int \frac{\partial}{\partial \vs_x}\lt(\rho_{\phi(x)}(z)\rt) dz=0$. So 
\begin{eqnarray}
\label{fineq100}
&~&\lt|\int u_{,\vs_x}(x-z)\frac{\partial}{\partial \vs_x}\lt(\rho_{\phi(x)}(z)\rt) dz\rt|\nn\\
&~&\qd\qd\qd=\lt|\int \lt(u_{,\vs_x}(x-z)-1\rt)\frac{\partial}{\partial \vs_x}\lt(\rho_{\phi(x)}(z)\rt) dz\rt|\nn\\
&~&\qd\qd\qd\overset{(\ref{ula120}),(\ref{abb5})}{\leq} c \ep^{-\frac{1}{2}}u(x)\lt|\int \phi_{,\vs_x}\lt(x\rt) \lt(\na \rho\lt(\frac{z}{\phi(x)}\rt)\cdot z
\lt(\phi(x)\rt)^{-4}+2 \rho\lt(\frac{z}{\phi(x)}\rt)\lt(\phi(x)\rt)^{-3}\rt)\rt|\nn\\
&~&\qd\qd\qd \leq c\ep^{-\frac{1}{2}}.
\end{eqnarray}
Finally we estimate the first term from (\ref{fineq35}) 
\begin{eqnarray}
\label{fineq101}
\lt|\int u_{,\vs_x \vs_x}(x-z)\rho_{\phi(x)}(z) dz\rt|\leq 
\|\na^2 u\|_{L^{\infty}(B_{4\rho_{\phi(x)}}(x))}\lt|\int\rho_{\phi(x)}(z) dz\rt|
\overset{(\ref{fineq60})}{\leq} c\ep^{-\frac{1}{2}}.
\end{eqnarray}
Putting (\ref{fineq99}), (\ref{fineq100}) and (\ref{fineq101}) together and applying 
this to (\ref{fineq35}) we have 
\begin{equation}
\label{dd1}
\lt|\psi_{,\vs_x \vs_x}(x)\rt|\leq c\ep^{-\frac{1}{2}}\text{ for any }x\in 
N_{8\ep}(\partial \Omega)\cap \Omega.
\end{equation}
Now by (\ref{ula21}) for any $x\in N_{8\ep}(\partial \Omega)\cap \Omega$, 
\begin{eqnarray}
\label{dd2}
\lt|\psi_{,\omega_x \vs_x}(x)\rt|&\leq&\int \lt|\na^2 u(x-z)\rt|\rho_{\phi(x)}(z) dz
+\int \lt|u_{,\omega_x}(z-x)\frac{\partial }{\partial \vs_x}(\rho_{\phi(x)}(z))\rt| dz\nn\\
&\overset{(\ref{abb5}),(\ref{fineq60})}{\leq}&c\ep^{-\frac{1}{2}}
+c\ep^{\frac{1}{2}}\int \lt|{\frac{\partial}{\partial\vs_x}}(\rho_{\phi(x)}(z))\rt| dz\nn\\
&\overset{(\ref{ula120})}{\leq}&c\ep^{-\frac{1}{2}}.
\end{eqnarray}
And by (\ref{ula21.5})
\begin{eqnarray}
\label{dd3}
\lt|\psi_{\omega_x \omega_x}(x)\rt|&\leq&\lt|\int u_{,\omega_x \omega_x}(x-z)\rho_{\phi(x)}(z) dz\rt|\nn\\
&\overset{(\ref{fineq60})}{\leq}&c\ep^{-\frac{1}{2}}.
\end{eqnarray}
Putting (\ref{dd1}), (\ref{dd2}), (\ref{dd3}) together establishes (\ref{dd4}).

\em Proof of Lemma completed. \rm From Step 2, (\ref{ula400.6}), for any $x\in N_{8\ep}(\partial\Omega)\cap \Omega$ we have 
$$
\lt|\lt|\na \psi(x)\rt|^2-1\rt|^2\leq c\ep
$$ 
so (\ref{ula351}) follows. In the same way from Step 4 (\ref{dd4}), (\ref{ula357}) follows. 

Since for any $x\in \Omega\backslash N_{8\ep}(\partial\Omega)$ we know $u(x)\geq \ep$ and so $\phi(x)=w(u(x))=\ep$ and thus 
$\rho_{\phi(x)}(z)=\rho\lt(\frac{z}{\ep}\rt)\ep^{-1}$ and there for $\psi(x)=\int u(x-z)\rho_{\ep}(z) dz$. Thus 
(\ref{ulazbc1}) is established. Finally by (\ref{uz1}), (\ref{uz22}) follows. $\Box$

%
%

\begin{a1}
\label{LL14}
Let $\Omega$ be a convex domain and $\lt|\Omega\triangle B_1(0)\rt|\leq \beta$.
Let $u(x)=d(x,\partial \Omega)$ and for $\veps>0$
define $u _{\veps}:=u *\rho_{\veps}$. 
For any $a\in\Omega\backslash N_{4\veps}(\partial \Omega)$ we have
\begin{equation}
\label{acb25}
\lt|\lt|\na u_{\veps}(x)\rt|-1\rt|\leq c\veps^{-1}V(\na u,B_{4\veps}(a))\text{ for any }x\in B_{2\veps}(a).
\end{equation}
\end{a1}

\em Proof. \rm Firstly recall that since $u$ is concave and hence $\na u$ is BV. Let $w=\Xint{-}_{B_{4\veps}(a)} \na u$.
By Poincare's inequality (see Remark 3.45 \cite{amb})
\begin{equation}
\label{ula60}
\int_{B_{4\veps}(a)} \lt|\na u-w\rt| dz\leq c\veps V\lt(\na u, B_{4\veps}(a)\rt).
\end{equation}

Now 
\begin{eqnarray}
\pi 16 \veps^2 \lt|1-\lt|w\rt|\rt|&=&\int_{B_{4\veps}(a)} \lt|1-\lt|w\rt|\rt| dz\nn\\
&=&\int_{B_{4\veps}(a)} \lt|\lt|\na u\rt|-\lt|w\rt|\rt| dz\nn\\
&\overset{(\ref{ula60})}{\leq}&c\veps V(\na u, B_{4\veps}(a)).\nn
\end{eqnarray}
Thus $\lt|1-\lt|w\rt|\rt|\leq c\frac{V(\na u, B_{4\veps}(a))}{\veps}$ and so there must exists $v\in S^1$ such that 
$\lt|v-w\rt|\leq\lt|1-\lt|w\rt|\rt|$ hence putting this together with (\ref{ula60}) we have 
\begin{equation}
\label{ulaz90}
\Xint{-}_{B_{4\veps}(a)} \lt|\na u-v\rt| dz\leq c\frac{V(\na u,B_{4\veps}(a))}{\veps}.  
\end{equation}
Hence for any $w\in B_{2\veps}(a)$
\begin{eqnarray}
\lt|\na u_{\veps}(w)-v\rt|&=&\lt|\int \lt(\na u(z)-v\rt)\rho_{\veps}(w-z) dz \rt|\nn\\
&\leq&c\veps^{-2}\lt|\int \lt(\na u(z)-v\rt)\rho(\veps^{-1}(z-w)) dz \rt|\nn\\
&\leq& c\veps^{-2}\int_{B_{2\veps}(w)}\lt|\na u(z)-v\rt| dz \nn\\
&\overset{(\ref{ulaz90})}{\leq}&c\frac{V(\na u,B_{4\veps}(a))}{\veps}.  
\end{eqnarray}
This completes the proof of Lemma \ref{LL14}. $\Box$

%
%
%

\begin{a1}
\label{LL26}
Let $\Omega$ be a convex domain and $\lt|\Omega\triangle B_1(0)\rt|\leq \beta$.
Let $u(x)=d(x,\partial \Omega)$ and define $u_{\veps}:=u*\rho_{\veps}$. Define
$\Lambda:=\Omega\backslash \lt(N_{8\veps}(\partial \Omega)\cup B_{4\beta^{\frac{1}{8}}}(0)\rt)$, we will
show that for any $\veps\in (0,\frac{\beta^{\frac{1}{2}}}{4}]$
\begin{equation}
\label{ula350}
\int_{\Lambda} \veps^{-1}\lt|1-\lt|\na u_{\veps}\rt|^2\rt|^2+\veps \lt|\na^2 u_{\veps}\rt|^2 dz\leq c\beta^{\frac{3}{16}}.
\end{equation}

\end{a1}
\em Proof of Lemma. \rm By the 5r Covering Theorem (\cite{mat}, Theorem 2.1) them we can find a finite collection
of balls $J:=\lt\{B_{\frac{2\veps}{5}}(x_i):i=1,2,\dots m\rt\}$ that are piecewise disjoint and
$\Lambda\subset \bigcup_{i=1}^m B_{2\veps}(x_i)$.

Note that for any $i=1,2,\dots n$ since the set of ball in $J$ are pairwise disjoint, for some 
constant $C_1$ there are at most $C_1$ balls from the set $\lt\{B_{5\veps}(x_k):k=1,\dots m\rt\}$ intersecting 
$B_{5\veps}(x_i)$. Thus $\|\sum_{i=1}^m \cha_{B_{5\veps}(x_i)}\|_{L^{\infty}(\Omega)}\leq C_1$ and this obviously implies
$\|\sum_{i=1}^m \cha_{B_{2\veps}(x_i)}\|_{L^{\infty}(\Omega)}\leq C_1$.

For $x,y\in \R^2$ let $x\otimes y:=\lt(\begin{smallmatrix} x_1 y_1 & x_1 y_2 \\ x_2 y_1 & x_2 y_2\end{smallmatrix}\rt)$. 
Now given $a\in \Lambda$ if $x\in B_{2\veps}(a)$, let $w=\Xint{-}_{B_{\veps}(x)} \na u$
\begin{eqnarray}
\label{uz12}
\lt|\na^2 u_{\veps}(x)\rt|&=&\lt|\int \na u(z)\otimes \na\rho_{\veps}(x-z) dz\rt|\nn\\
&\leq&\lt|\int (\na u(z)-w)\otimes \na\rho\lt(\frac{x-z}{\veps}\rt)\veps^{-3} dz\rt|\nn\\
&\leq&c\veps^{-3}\lt|\int_{B_{2\veps}(x)} (\na u-w) dz\rt|\nn\\
&\overset{(\ref{ula60})}{\leq}&c\veps^{-2}V(\na u,B_{4\veps}(a)).
\end{eqnarray}

So 
\begin{eqnarray}
\label{ula984}
\int_{\Lambda} \lt|\na^2 u_{\veps}\rt|^2 dz&\leq&\sum_{i=1}^m c\int_{B_{2\veps}(x_i)} \lt|\na^2 u_{\veps}\rt|^2 dz\nn\\
&\leq&c\sum_{i=1}^m \veps^2\|\na^2 u_{\veps} \|^2_{L^{\infty}(B_{2\veps}(x_i))}\nn\\
&\overset{(\ref{uz12})}{\leq}&c\veps^2 \lt(\sum_{i=1}^m \veps^{-4}\lt(V(\na u,B_{4\veps}(x_i))\rt)^2\rt)\nn\\
&\overset{(\ref{abb40})}{\leq}&c\beta^{\frac{3}{16}}\veps^{-1}
\lt(\sum_{i=1}^m  V(\na u,B_{4\veps}(x_i))\rt)\nn\\
&\leq&c\beta^{\frac{3}{16}}\veps^{-1}V\lt(\na u,\Lambda\rt)\nn\\
&\overset{(\ref{ula301})}{\leq}&c\veps^{-1}\beta^{\frac{3}{16}}.
\end{eqnarray}

Now 
\begin{eqnarray}
\label{uz9.6}
\int_{\Lambda} \lt|1-\lt|\na u_{\veps}\rt|^2\rt|^2 dz
&\leq&c\sum_{i=1}^m \int_{B_{2\veps}(x_i)} \lt|1-\lt|\na u_{\veps}\rt|\rt|^2 dz\nn\\
&\overset{(\ref{abb40}),(\ref{acb25})}{\leq}& \sum_{i=1}^m c\veps^2\beta^{\frac{3}{16}}\|
\lt|1-\lt|\na u_{\veps}\rt|\rt|\|_{L^{\infty}(B_{2\veps}(x_i))}\nn\\
&\overset{(\ref{acb25})}{\leq}&\sum_{i=1}^m c \veps \beta^{\frac{3}{16}} V(\na u, B_{4\veps}(x_i))\nn\\
&\leq&c\veps \beta^{\frac{3}{16}} V(\na u,\Omega\backslash B_{2\beta^{\frac{1}{8}}}(0))\nn\\
&\overset{(\ref{ula301})}{\leq}&c \beta^{\frac{3}{16}}\veps.
\end{eqnarray}
Putting (\ref{uz9.6}) together with (\ref{ula984}) establishes (\ref{ula350}). $\Box$

%
%
%

\begin{a1}
\label{LL17}
Let $\eta(x)=\lt|x\rt|$, $\veps>0$ and define $\eta_{\veps}(x):=\int \eta(z)\rho_{\veps}(x-z) dz$. Then  
\begin{equation}
\label{acb30}
\int_{B_1(0)} \lt|1-\lt|\na \eta_{\veps}\rt|^2\rt|^2 dz\leq c\log(\veps^{-1})\veps^2
\end{equation}
and 
\begin{equation}
\label{acb31}
\int_{B_1(0)} \lt|\na^2 \eta_{\veps}\rt|^2 dz\leq c\log(\veps^{-1}).
\end{equation}
\end{a1}
%
%
\em Proof of Lemma. \rm Note for $x\not \in B_{2\veps}(0)$, $z\in B_{\veps}(x)$
\begin{eqnarray}
\label{ula99}
\lt|\frac{z}{\lt|z\rt|}-\frac{x}{\lt|x\rt|}\rt|&\leq&\lt|\frac{z\lt|x\rt|-x\lt|z\rt|}{\lt|z\rt|\lt|x\rt|}\rt|\nn\\
&\leq&\lt|\frac{z\lt|x\rt|-x\lt|x\rt|}{\lt|z\rt|\lt|x\rt|}\rt|
+\lt|\frac{x\lt|x\rt|-x\lt|z\rt|}{\lt|z\rt|\lt|x\rt|}\rt|\nn\\
&\leq&\frac{c\veps}{\lt|x\rt|-\veps}.
\end{eqnarray}
So for $x\not\in B_{4\veps}(0)$
\begin{eqnarray}
\label{ula101}
\lt|\na \eta_{\veps}(x)-\frac{x}{\lt|x\rt|}\rt|&=&
\lt| \int \rho_{\veps}(x-z)\lt(\frac{x}{\lt|x\rt|}-\frac{z}{\lt|z\rt|}\rt) dz \rt|\nn\\
&\overset{(\ref{ula99})}{\leq}&\frac{c\veps}{\lt|x\rt|-\veps}.
\end{eqnarray}
Since $\int \frac{x}{\lt|x\rt|}\otimes \na \rho_{\veps}(x-z) dz=0$, for any $x\not\in B_{4\veps}(0)$
\begin{eqnarray}
\label{aaz2}
\na^2 \eta_{\veps}(x)&=&\lt|\int \na \eta_{\veps}(z)\otimes \na\rho_{\veps}(x-z) dz\rt|\nn\\
&=&\lt|\int \lt(\na\eta_{\veps}(z)-\frac{z}{\lt|z\rt|}\rt) \otimes \na\rho_{\veps}(x-z) dz \rt|
+\lt|\int \lt(\frac{x}{\lt|x\rt|}-\frac{z}{\lt|z\rt|}\rt) \otimes \na\rho_{\veps}(x-z) dz \rt|\nn\\
&\overset{(\ref{ula99}),(\ref{ula101})}{\leq}& \frac{c\veps}{\lt|x\rt|-\veps}\lt|\int \na\rho_{\veps}(x-z) dz \rt|\nn\\
&\leq& \frac{c}{\lt|x\rt|-\veps}.
\end{eqnarray}

Hence 
\begin{eqnarray}
\label{ull2}
\int_{B_1(0)\backslash B_{4\veps}(0)} \lt|\na^2 \eta_{\veps}(x)\rt|^2 dx&\overset{(\ref{aaz2})}{=}&c\int_{4\veps}^1 
\int_{\partial B_h(0)} \lt(\frac{1}{\lt|z\rt|-\veps}\rt)^2  dH^1 z dr\nn\\
&\leq& c\int_{\veps}^1 \frac{1}{r} dr\nn\\
&\leq& c\log(\veps^{-1})
\end{eqnarray}
which establish. Now as $\lt|\na \eta_{\ep}(x)\rt|\leq c$ and 
$\lt|\na^2 \eta_{\ep}(x)\rt|\leq c\ep^{-1}$ for any $x\in B_{1-\ep}(0)$ so 
$$
\int_{B_{4\ep}(0)} \lt|\na^2 \eta_{\ep}\rt|^2 dz\leq c\ep.
$$
Thus putting this together with (\ref{ull2}) establishes (\ref{acb31}). 

Note $\lt|\lt|\na \eta_{\veps}(x)\rt|-1\rt|^2\leq \lt|\na \eta_{\veps}(x)-\frac{x}{\lt|x\rt|}\rt|^2\overset{(\ref{ula101})}{\leq} c\frac{\ep^2}{(\lt|x\rt|-\ep)^2}$ so arguing in the same way as in (\ref{ull2}) we have (\ref{acb30}). $\Box$

%
%

\subsection{Proof of Proposition \ref{PP1}}

Let $u(x)=d(x,\partial \Omega)$,  
let $w:\R_{+}\rightarrow \R_{+}$ be the smooth monotonic function from the proof of Lemma \ref{LL7}, 
so $w$ satisfies (\ref{ula70}) and $\sup \lt|\ddot{w}\rt|\leq c\ep^{-1}$ as in  
Lemma \ref{LL7} for 
$x\in N_{\ep}(\partial \Omega)\cap \Omega$ define 
\begin{equation}
\label{kband0203}
\phi(x)=w(u(x)). 
\end{equation}
Let  
\begin{equation}
\label{ula400}
v(x):=\min\lt\{u(x),1-8\beta^{\frac{3}{32}}+\lt|x\rt|\rt\}.
\end{equation}
and define 
\begin{equation}
\label{ula401}
\xi(x)=\int v(x-z)\rho_{\phi(x)}(z) dz.
\end{equation}

Let $\Pi:=\lt\{x:u(x)>1-8\beta^{\frac{3}{32}}+\lt|x\rt|\rt\}$, and 
define $\Lambda_0:=\Omega\backslash (N_{8\ep}(\partial \Omega)\cup N_{\ep}(\Pi))$, note that 
$\xi(x)=u_{\ep}(x)$ for any $x\in \Lambda_0$.  

Recall from (\ref{qwer1}) the function $\psi$ defined in Lemma \ref{LL7}. Note that for any 
$x\in N_{8\ep}(\partial \Omega)\cap \Omega$ function $\phi$ we defined by (\ref{kband0203}) 
is identical to $\phi$ defined by (\ref{kband0202}) in Lemma \ref{LL7}. Hence as $u=v$ in 
$N_{8\ep}(\partial \Omega)\cap \Omega$ we have $\xi(x)=\psi(x)$ for any 
$x\in N_{8\ep}(\partial\Omega)\cap \Omega$ thus 
\begin{equation}
\label{ula802}
\int_{N_{8\ep}(\partial\Omega)\cap \Omega} \ep^{-1}\lt|1-\lt|\na \xi\rt|^2\rt|^2+\ep\lt|\na^2 \xi\rt|^2 dx\nn\\
\overset{(\ref{ula351}),(\ref{ula357})}{\leq} c\ep
\end{equation}
Since $\psi=u_{\ep}$ in $\Lambda_0$, from (\ref{ula350}) we have 
$\int_{\Lambda_0} \ep^{-1}\lt|1-\lt|\na\xi\rt|^2\rt|^2+\ep\lt|\na^2 \xi\rt|^2 dx\leq 
c\beta^{\frac{3}{16}}$ and so putting this two inequalities together we have 
\begin{equation}
\label{ula370}
\int_{\Omega\backslash N_{\ep}(\Pi)} \ep^{-1}\lt|1-\lt|\na\xi\rt|^2\rt|^2+\ep\lt|\na^2 \xi\rt|^2 dx
\leq c\beta^{\frac{3}{16}}
\end{equation}

Now as for any $x\in \Pi\backslash N_{\ep}(\partial \Pi)$, $w(x)=1-8\beta^{\frac{3}{32}}+\lt|x\rt|$ and so 
$u_{\ep}(x)=\eta_{\ep}(x)+(1-8\beta^{\frac{3}{32}})$ where $\eta(x)=\lt|x\rt|$ and $\eta_{\ep}=\eta*\rho_{\ep}$. 
So $\na \xi(x)=\na \eta_{\ep}(x)$ and $\na^2 \xi(x)=\na^2 \eta_{\ep}(x)$ thus applying Lemma \ref{LL17} we have 
\begin{equation}
\label{ula501}
\int_{\Pi\backslash N_{\ep}(\partial \Pi)} \ep^{-1}\lt|1-\lt|\na \xi\rt|^2\rt|^2+\ep\lt|\na^2 \xi\rt|^2 dx
\overset{(\ref{acb30}),(\ref{acb31})}{\leq} c\ep\log(\ep^{-1}).
\end{equation}

Since $w$ is Lipschitz, so $\xi$ is Lipschitz and so from (\ref{azz32}) we have 
\begin{equation}
\label{ula800}
\int_{N_{\ep}(\partial \Pi)} \ep^{-1}\lt|1-\lt|\na \xi\rt|^2\rt|^2 dx\leq c\beta^{\frac{3}{32}}. 
\end{equation}
And note for any $x\in\Omega\backslash N_{\ep}(\partial \Omega)$
$$
\lt|\na^2 \xi(x)\rt|=\ep^{-3}\lt|\int \na v(z)\cdot \na\rho\lt(\frac{x-z}{\ep}\rt) dz\rt|\leq c\ep^{-1}
$$
so 
\begin{eqnarray}
\label{ula988}
\int_{N_{\ep}(\partial \Pi)}  \ep\lt|\na^2 \xi\rt|^2 dx
&\leq&c\ep^{-1}\lt|N_{\ep}(\partial \Pi)\rt|\nn\\
&\overset{(\ref{azz32})}{\leq}&c\beta^{\frac{3}{32}}.
\end{eqnarray}

Putting these inequalities together we have 
\begin{equation}
\label{uz21}
\int_{N_{\ep}(\partial \Pi)} \ep^{-1}\lt|1-\lt|\na \xi\rt|^2\rt|^2+\ep\lt|\na^2 \xi\rt|^2 dx
\leq c\beta^{\frac{3}{32}}.
\end{equation}
Now inequalities (\ref{ula370}), (\ref{ula501}) and (\ref{uz21}) give us that $\xi$ satisfies 
(\ref{ula201}). And since $\xi(x)=\psi(x)$ on $N_{\ep}(\partial\Omega)\cap \Omega$ from (\ref{uz22}) 
satisfies $\na \xi(x)\cdot \eta_x=1$ for any $x\in \partial \Omega$. This completes the proof of Proposition 
\ref{PP1}. $\Box$

\subsection{Proof Corollary \ref{CC1}}

Let $\alpha=\inf_{y\in \Omega} \lt|\Omega\triangle B_1(y)\rt|$.  Let $\beta=4(\alpha+\ep)$, 
note that since we can assume without loss of generality that $\alpha+\ep\leq \frac{1}{4}$ so $\beta\leq 1$ which gives 
$\beta\leq \beta^{\frac{1}{2}}$ and so $\ep\leq \frac{\beta^{\frac{1}{2}}}{4}$. Now we can also assume without loss of generality that $\lt|\Omega\triangle B_1(0)\rt|\leq \beta$. So 
we can apply Proposition \ref{PP1} which gives us the existence of $\xi\in \Lambda(\Omega)$ such that(\ref{ula201}) hold true. Hence we have 
that $\inf_{u\in \Lambda(\Omega)} I_{\ep}(u)\leq c\beta^{\frac{3}{32}}$. Let $v\in \Lambda(\Omega)$ 
be the minimiser of $I_{\ep}$ and since $v$ satisfies 
$$
\int_{\Omega} \lt|1-\lt|\na v\rt|^2\rt|\lt|\na^2 v\rt| dz\leq 
\int_{\Omega} \ep^{-1}\lt|1-\lt|\na v\rt|^2\rt|^2+
\ep\lt|\na^2 v\rt|^2 dz\leq c\beta^{\frac{3}{32}} 
$$
and as $\ep\in (0,\frac{\beta^{\frac{1}{2}}}{4})$
\begin{equation}
\label{uu6}
\int_{\Omega} \lt|1-\lt|\na v\rt|^2\rt|^2 dz\leq c\beta^{\frac{19}{32}}.
\end{equation}
So we have that (\ref{eq1}), (\ref{eq2}) are satisfied 
and hence by Theorem \ref{T1}
$$
\int_{\Omega} \lt|\na v(z)+\frac{z}{\lt|z\rt|}\rt|^2 dz\leq c\beta^{\frac{1}{5462}}.
$$
Applying Lemma \ref{LL25} we have 
$\int_{\Omega\backslash B_{\beta^{\frac{1}{8}}}(0)} \lt|\na v-\na \zeta\rt|^2\leq c\beta^{\frac{1}{5462}}$. So arguing is the same way as the proof of Corollary \ref{CC3} 
we have $\|v-\zeta\|_{W^{1,2}(\Omega)}\leq 
c\beta^{\frac{1}{5462}}\leq c(\ep+\alpha)^{\frac{1}{5462}}$. $\Box$

\end{document}